\documentclass[10pt]{amsart}
\usepackage{amssymb}
\usepackage{graphicx}
\usepackage{xcolor} % A package to add color.
\usepackage{tensor}
\usepackage{fullpage} % Sets all margins to 1 in.  
\usepackage{amsmath}
\usepackage{amsthm}
\usepackage{verbatim}
\usepackage{hyperref}
\usepackage{enumitem}
\usepackage{comment}
\usepackage{booktabs,tabularx}
\setlist[enumerate]{leftmargin=1.5em}
\setlist[itemize]{leftmargin=1.5em}

\providecommand{\MR}{\relax\ifhmode\unskip\space\fi MR }
\usepackage{xcolor}
\definecolor{purple}{rgb}{0.5, 0, 1}
\definecolor{orange}{rgb}{1,.5,0}

\providecommand{\MR}{\relax\ifhmode\unskip\space\fi MR }
\usepackage{seqsplit} 
\definecolor{green}{rgb}{0,0.8,0} % Redefines the color green.
\newtheorem{theorem}{Theorem}[section]
\newtheorem{corollary}[theorem]{Corollary}
\newtheorem{lemma}[theorem]{Lemma}
\newtheorem{proposition}[theorem]{Proposition}

\theoremstyle{definition}
\newtheorem{definition}[theorem]{Definition}

\theoremstyle{definition}
\newtheorem{remark}[theorem]{Remark}

\numberwithin{equation}{section}
\newcommand{\nrm}[1]{\Vert#1\Vert}

\newcommand{\tld}[1]{\widetilde{#1}}

\newcommand{\nnrm}[1]{{\vert\kern-0.25ex\vert\kern-0.25ex\vert #1 
		\vert\kern-0.25ex\vert\kern-0.25ex\vert}}

\newcommand{\dist}{\mathrm{dist}\,}
\newcommand{\sgn}{{\mathrm{sgn}}\,}
\newcommand{\supp}{{\mathrm{supp}}\,}

\newcommand{\lap}{\Delta}

\newcommand{\nb}{\nabla}

\newcommand{\gmm}{\gamma}

\newcommand{\dlt}{\delta}

\newcommand{\kpp}{\kappa}

\newcommand{\tht}{\theta}

\newcommand{\omg}{\omega}
\newcommand{\Omg}{\Omega}

\newcommand{\bbR}{\mathbb R}

\newcommand{\bbT}{\mathbb T}

\newcommand{\calA}{\mathcal A}

\newcommand{\calC}{\mathcal C}

\newcommand{\calD}{\mathcal D}
\newcommand{\calE}{\mathcal E}
\newcommand{\calF}{\mathcal F}

\newcommand{\calO}{\mathcal O}
\newcommand{\calP}{\mathcal P}

\newcommand{\calR}{\mathcal R}
\newcommand{\calS}{\mathcal S}
\newcommand{\calT}{\mathcal T}

\newcommand{\varep}{\varepsilon}

\begin{document}
	
	\title{A variational theory for von K\'arm\'an vortex streets
    %\\2. Variational construction of von K\'arm\'an vortex street\\3. Variational von K\'arm\'an vortex street\\    4. Energy maximizing von K\'arm\'an vortex streets\\    5. Energy maximization in periodic Euler flow: von K\'arm\'an vortex streets\\    6. ?
    } 
	
	\date{\today}
	
	\renewcommand{\thefootnote}{\fnsymbol{footnote}}
	\footnotetext{\emph{Key words: Incompressible Euler equation, von K\'arm\'an vortex street, variational principle, orbital stability.}  \\
		\emph{2020 AMS Mathematics Subject Classification:} Primary 76B47; Secondary 35Q35, 35B40. \\
		\emph{$^*$corresponding author}}

	\renewcommand{\thefootnote}{\arabic{footnote}}

	%\date{\today}%
	%\dedicatory{}%
	%\commby{}%
	% ----------------------------------------------------------------
	
	\author{In-Jee Jeong}
	\address{School of Mathematics, Korea Institute for Advanced Study,  Seoul 02455, Republic of Korea}
	\email{ijeong@kias.re.kr} 
	
	\author{Guolin Qin}
	\address{Academy of Mathematics and Systems Science, Chinese Academy of Sciences, Beijing 100190, China}
	\email{qinguolin18@mails.ucas.ac.cn} 
	
	\author{Kwan Woo$^*$}
	\address{Departement Mathematik und Informatik, Universit\"at Basel, CH-4051 Basel, Switzerland}
	\email{kwan.woo@unibas.ch}
	
	\begin{abstract} We construct von K\'arm\'an vortex streets as global kinetic energy maximizers for the 2D Euler flow in an infinite periodic channel: our variational framework imposes no a prior symmetry, support constraints, or sign conditions on the admissible class. More specifically, under $O(\varepsilon^{-2})$ constraint of the absolute vorticity maximum with other natural constraints, we show that the maximizers are traveling pairs of simply connected vortex patches of support radius $O(\varepsilon)$, converging to the von K\'arm\'an point-vortex configuration satisfying glide reflection symmetry as $\varepsilon \to 0$. We prove nonlinear stability of the von K\'arm\'an streets under a bound on the vertical support, using a compactness result valid for more general assumptions. Furthermore, we construct slightly tilted asymmetric  von K\'arm\'an vortex streets as energy maximizers under a tilted moment constraint.
	\end{abstract}
	
	\maketitle
	
	% ----------------------------------------------------------------

	%\tableofcontents \newpage 
	
	\section{Introduction}
	
	A von K\'arm\'an vortex street is a repeating pattern of counter-rotating vortices, frequently observed in flows behind blunt bodies.
    Once such a pattern is created, it is effectively described as a relative equilibrium of the inviscid two-dimensional fluid flow.
    More precisely, we consider von K\'arm\'an vortex streets as traveling waves of the $2\pi$-periodic infinite cylinder $\calC := \bbT \times \bbR$.
    In this domain, the vorticity equation for incompressible and inviscid flow takes the form 
	\begin{equation}\label{eq:Euler_cyl}
		\left\{
		\begin{aligned}
			&\partial_{t} \omega + u \cdot \nabla \omega = 0,\\
			&u = \nabla^\perp(-\Delta_{\calC})^{-1}\omega,
		\end{aligned}
		\right.
	\end{equation}
	where $\nabla^\perp = (\partial_{2}, -\partial_{1})$ and the \emph{Green's function} $G$ for $(-\Delta_{\calC})^{-1}$ is given explicitly by \begin{equation}\label{eq:Green}
		\begin{split}
			G(x, y) = G(x-y) = -\frac{1}{4\pi} \log \left( 2\cosh  \left  (x_2 - y_2 \right) - 2 \cos  \left( x_1 - y_1\right) \right),
		\end{split}
	\end{equation}
    so that \begin{equation}\label{eq:Biot-Savart-kernel}
		\begin{split}
			u(x)= \int_{\calC} K(x-y) \omega(y) \,dy, \quad K(x) := \nabla^{\perp}G(x) = \frac{1}{4\pi} \frac{ (-\sinh x_{2},  \sin x_{1} ) }{\cosh x_2 - \cos x_1}.
		\end{split}
	\end{equation}
    (The difference $x_1-y_1$ is represented by the unique element of $\bbT = [-\pi,\pi)$.%, and we write $|x-y| := {(|x_1-y_1|^2+|x_2-y_2|^2)^{1/2}}$.
	) For an initial datum $\omg_{0} \in (L^\infty \cap L^{1})(\calC)$, there is a unique global in time solution which belongs to $(L^\infty \cap L^{1})(\calC)$ for each $t \in \bbR$ (\cite{BD17}).  
	
\begin{figure}[h]
			\centering
			\includegraphics[scale=0.8]{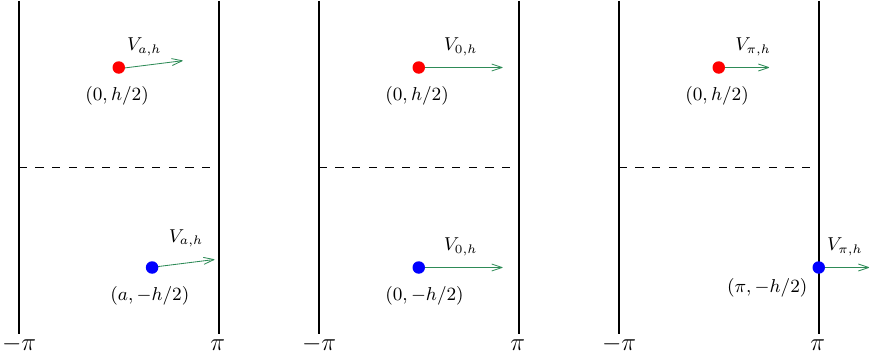}  
			\caption{The point vortex configuration; the right two figures represent symmetric ones. Asymmetric configurations have a nonzero vertical velocity component.} \label{fig:pv-Karman}
		\end{figure}

	The simplest way to model a von K\'arm\'an street is to employ \textit{point vortices} on $\calC$ (\cite{MP,Ga20,CQZZ23CVPDE}), as the following point vortex configuration of two Dirac deltas (see Figure \ref{fig:pv-Karman}): \begin{equation}\label{eq:pv-Karman}
		\begin{split}
			\omg_{a,h}^{pv} := \frac{1}{2}\dlt_{(0,h/2)} - \frac{1}{2}\dlt_{(a,-h/2)}, \quad \mbox{where} \quad a \in \bbT, h>0. 
		\end{split}
	\end{equation}
    It defines a uniformly translating point vortex configuration on $\calC$ with velocity \begin{equation}\label{eq:pv-Karman-vel}
		\begin{split}
			V_{a,h} := \frac{1}{8\pi} \frac{ (\sinh(h), \sin(a)) }{ \cosh(h) - \cos(a) }. 
		\end{split}
	\end{equation} We shall refer to the special cases $a = 0, \pi$ as the \textbf{symmetric} von K\'arm\'an vortex streets; they are precisely characterized by  \textit{vanishing} of the vertical velocity (\cite[Proposition 2.2]{Ga20}). This notion of symmetry carries over to general vorticity distributions: $a = 0$ corresponds to the vorticity having odd symmetry in $x_{2}$ and $a = \pi$ to \textit{glide reflection} symmetry - namely, invariance under the transform $T : (x_1,x_2) \mapsto (x_1+\pi,-x_2)$ in $\calC$. These symmetries are preserved by  the unique $L^\infty \cap L^1$ solutions of \eqref{eq:Euler_cyl}.

	In this work, we study the problem of finding \textbf{global} kinetic energy maximizers under natural constraints in $\calC$, and derive patch {desingularizations} of \eqref{eq:pv-Karman} as the maximizers.
    This calculus of variations approach not only gives existence of (a one-parameter family of) traveling wave solutions of \eqref{eq:Euler_cyl} approximating \eqref{eq:pv-Karman}, but also provides their dynamical stability under some conditions.
    To the best of our knowledge, there was no mathematical works on nonlinear stability of such desingularized von K\'arm\'an vortices.
    There are several difficulties in obtaining stability for them, and we shall discuss them after stating the main results.

	%In these two special cases, the  velocity is given by \begin{itemize}	\item $V_{0,h}=(\frac{1}{8\pi}\coth(\frac{h}{2}),0)$, and
		%\item $V_{\pi,h}=(\frac{1}{8\pi}\tanh(\frac{h}{2}),0)$,
		%\end{itemize} respectively.

		\subsection{Main results} To state our main results, let us introduce several functionals of the vorticity.
        To begin with, we consider vorticities in $\calC$ satisfying the \emph{zero total circulation}   $$\int_{\calC}\omega(x) dx= 0,$$ which is preserved by solutions of \eqref{eq:Euler_cyl}.
        This condition is automatic for $\omg = \nb \times u$ with $u\to 0$ as $|x_2|\to\infty$.
        Moreover, under this condition, the (kinetic) \textbf{energy} can be viewed as a functional of the vorticity:
		\[
		E[\omega] = \frac12\int_{\calC} |u(x)|^2 dx =  \frac{1}{2} \iint_{\calC \times \calC} G(x - y) \omega(x) \omega(y) dx dy = \frac{1}{2} \int_{\calC} \psi(x)\omega(x)dx,
		\]
        where $\psi[\omega](x) := \int_{\calC} G(x-y) \omega(y)dy$ is the \textbf{stream function} for $\omg$. 
		
		Next, we introduce the \textbf{impulse} $I(\omg)$ and the \textbf{absolute impulse} $\nrm{x_{2}\omg}_{L^1(\calC)}$: \begin{equation}\label{eq:impulse-and-absolute}
			\begin{split}
				I(\omg) = \int_{\calC} x_{2}\omg(x) dx , \qquad \nrm{x_{2}\omg}_{L^1(\calC)} =  \int_{\calC} |x_{2}\omg(x)| dx.  
			\end{split}
		\end{equation}
        While $I(\omg)$ is conserved in time (this requires the zero circulation condition and finiteness of absolute impulse), in general $\nrm{x_{2}\omg}_{L^1(\calC)}$ is not; they coincide only when $x_{2}\omg(x) \ge 0$ almost everywhere. 
		
		Recall that all the $L^{p}(\calC)$ norms ($1 \le p \le \infty$) of $\omg$ are conserved in time, for solutions of \eqref{eq:Euler_cyl}.
        In what follows, we shall use the shorthand $\nrm{\omg}_{p} := \nrm{\omg}_{L^{p}(\calC)}$. 
		
		Let us now introduce the following class of vorticities: 
		\begin{equation}\label{eq:adm-main}
			\begin{split}
				\calF^{\mathrm{cen}}_{\mu,\varepsilon} := \left\{ \omega \, :\, \|\omega\|_1 \le 1, \|\omega\|_{\infty} \le \varepsilon^{-2}, \|x_2\omega\|_1 \le \mu,\int_{\calC}  \omega(x)\,dx =0,  \int_{\calC} x_{2}|\omg| dx = 0 , \, \int_{\calC} x_{2}\omg dx \ge 0 \right\}.
			\end{split}
		\end{equation}
		For any $\varepsilon > 0$, it is not difficult to prove that the set of energy maximizers is always nonempty; see \S \ref{sec:prelim} below.
		Moreover, for sufficiently small $\varepsilon > 0$, imposing the last two requirements in \eqref{eq:adm-main} does not change the maximal energy; see Remark~\ref{rmk:cen_max} (the notation ``cen'' is for centered). 
		We denote \begin{equation*}
			\begin{split}
				\calS^{\mathrm{cen}}_{\mu,\varepsilon} := \left\{ \bar\omg \in \calF^{\mathrm{cen}}_{\mu,\varepsilon} \, : \, E[\bar\omg] = \sup_{ \omg \in \calF^{\mathrm{cen}}_{\mu,\varepsilon}} E[\omg] \right\} \ne \emptyset. 
			\end{split}
		\end{equation*}
		Our first main result describes the properties of each maximizer.
        Furthermore, we prove that, in the limit $\varepsilon\to 0^+$, any sequence of maximizers converges to the von K\'arm\'an vortex street with glide reflection symmetry. 
		\begin{theorem}[Global maximizer] \label{thm:max}
			For a fixed $\mu > 0$, the following holds for all sufficiently small $\varepsilon > 0$. Any   $\omg_{\varepsilon} \in \calS^{\mathrm{cen}}_{\mu,\varepsilon}$  takes the form \begin{equation*}
				\begin{split}
					\omg_{\varepsilon} = \varepsilon^{-2}\big( \mathbf{1}_{\Omg_+^\varepsilon} - \mathbf{1}_{\Omg_{-}^\varepsilon} \big)
				\end{split}
			\end{equation*}
            for some simply connected open sets $\Omg_{\pm}^{\varepsilon}$ disjoint from each other, satisfying $\mathrm{diam}(\Omg_{\pm}^{\varepsilon}) \le C(\mu) \varepsilon$, and
			\[
			\Omg_{+}^{\varepsilon} = \left\{  x_2 > 0 \, : \,  \psi[\omg_{\varepsilon}] - W_{\varepsilon}x_{2} - \gmm_{+,\varepsilon} > 0 \right\}, \quad \Omg_{-}^{\varepsilon} = \left\{  x_2 < 0 \, : \,  -\psi[\omg_{\varepsilon}] + W_{\varepsilon}x_{2} - \gmm_{-,\varepsilon} > 0 \right\}
			\]
			for some $\gmm_{\pm,\varepsilon}, W_{\varepsilon} > 0$. 
			In particular, $\omg_{\varepsilon}$ defines a traveling wave solution to \eqref{eq:Euler_cyl} with velocity $(W_\varepsilon,0)$.
            Furthermore, it saturates the impulse bound as follows:
            \begin{equation*}
				\begin{split}
					\int_{\calC} x_{2}\omg_{\varepsilon} dx = \int_{\calC} |x_{2}\omg_{\varepsilon}| dx = \mu. 
				\end{split}
			\end{equation*} 
			In the limit $\varepsilon \to 0^+$, we have: \begin{itemize}
				\item $W_{\varepsilon} \to \frac1{8\pi} \tanh(\mu)$,
				\item for $p_{+} := (0,  \mu)$ and $p_{-} := (\pi, - \mu)$, $\Omg^{\varepsilon}_{\pm}$ shrink to $p_{\pm}$ up to horizontal translations, and 
				\item for a proper choice of points $p^{\varepsilon}_{\pm} \in \Omg^{\varepsilon}_{\pm}$, we have the convergence of rescaled sets \begin{equation*}
					\begin{split}
						A_{\pm}^{\varepsilon} := \frac{\Omg^{\varepsilon}_{\pm}-p^{\varepsilon}_{\pm}}{\varepsilon} \longrightarrow B_{1/\sqrt{2\pi}}(0)
					\end{split}
				\end{equation*}
                in measure, where $B_{1/\sqrt{2\pi}}(0)$ is the ball of radius $1/\sqrt{2\pi}$ centered at the origin.
			\end{itemize}
		\end{theorem}

		The next main result gives compactness in $L^{1}$ (and hence in $L^{p}$ for all finite $p$ by interpolation) for the whole maximizing set, for $\varepsilon>0$ sufficiently small. 
		\begin{theorem}\label{thm:main-compact}
			Fix $\mu>0$ and take $\varepsilon>0$ sufficiently small.
            Let $\{ \omg_{n} \}_{n\ge0}$ be a maximizing sequence, which asymptotically belongs to $\calF^{\mathrm{cen}}_{\mu,\varepsilon}$.
            That is, \begin{equation*}
				\begin{split}
					&\int_{\calC}  \omega_n(x)\,dx =0, \quad \limsup_{n\to\infty}\|\omega_n\|_1 \le 1, \quad \limsup_{n\to\infty}\|\omega_n\|_{\infty} \le \varepsilon^{-2}, \quad \limsup_{n\to\infty}\|x_2\omega_n\|_1 \le \mu, \\
					&\int_{\calC} x_{2}|\omg_n| dx = 0 , \quad \int_{\calC} x_{2}\omg_n dx \ge 0, \qquad \mbox{and} \qquad \lim_{n\to\infty} E[\omg_{n}] = \sup_{ \omg \in \calF^{\mathrm{cen}}_{\mu,\varepsilon}} E[\omg].
				\end{split}
			\end{equation*}
            Then, there exist a subsequence $\{ \omega_{n_k}\}_k$ and $\bar\omg \in  \calS^{\mathrm{cen}}_{\mu,\varepsilon} $ such that $\nrm{\omg_{n_k} - \bar{\omg}}_{1} \to 0$ as $n\to\infty$. 
		\end{theorem}

		The above compactness result can be ``translated'' into nonlinear stability of the set $\calS_{\mu, \varepsilon}^{\rm cen}$ of maximizers, with one important caveat that the absolute impulse $\nrm{ x_2\omg }_{1}$ is not conserved.
        Still, when the vorticity is $L^{1}$ close to the energy maximizing configuration, this quantity can differ significantly from the conserved impulse only if the vorticity has some ``incorrectly'' signed parts very far from $x_{2} = 0$.
        Based on this intuition, one can obtain a conditional stability result, assuming that the support of the solution is not too widespread in the $x_{2}$ axis. 
		
		For convenience of notation, we introduce the translation (in both directions) orbit of the set of centered maximizers
		\begin{equation*}
			\begin{split}
				\calO^+_{\mu,\varepsilon} := \left\{ \tau_{a,b}\omega: \omega \in \calS^{\rm cen}_{\mu,\varepsilon}, \quad (a,b)\in\bbT\times\bbR \right\},
			\end{split}
		\end{equation*} and the following distance from a vorticity $\tht$ to the translation orbit: 
		\begin{equation*}
			\begin{split}
				d_{\calO^+}(\theta):= \inf_{\omega\in\calO^+_{\mu,\varepsilon}} \left\{ \| \tht - \omega\|_1 + |I(\tht) - \mu| \right\}. 
			\end{split}
		\end{equation*} 
		
		\begin{theorem}[Patch stability with support condition]\label{thm:main-stab} For a given $\mu > 0$, the following stability holds for all sufficiently small $\varepsilon>0$.
        Consider the following class of initial data: 
			\begin{equation}\label{eq:supp_patch_main}
				\theta(0) = \varepsilon^{-2} \big( \mathbf 1_{A_0^+}-\mathbf 1_{A_0^-} \big), \quad A_0^+ \cap A_0^- = \emptyset, \quad |A_0^{\pm}| = \frac{\varepsilon^2}{2}.
			\end{equation}
            Let $\theta(t)$ be the global solution to \eqref{eq:Euler_cyl} corresponding to $\tht(0)$.
            For every $\eta > 0$, there exist $\delta = \delta(\eta, \mu, \varepsilon) > 0$ and $c_0 = c_0(\mu) > 0$ such that the following holds:
            if  $d_{\calO^+}(\theta(0)) < \delta$ and
            \begin{equation}\label{eq:supp_stab_main}
            \sup_{t \in J} \sup_{x,y \in \supp\theta(t)} |x_2-y_2| \le c_0|\log\varepsilon|
			\end{equation}
            for an interval $0 \in J \subset \bbR$, then we have
            $$
            \sup_{t \in J} d_{\calO^+}(\theta(t)) < \eta.
            $$ 
		\end{theorem}

		\begin{remark}
			If the initial data $\tht(0)$ is supported close to the maximizer, it is expected that the support condition holds for a very long time, using the methods of Marchioro and Pulvirenti \cite[Chapter 4]{MP}.
            Furthermore, the support condition is just one way to obtain stability; we provide a more general condition as well as relevant discussions in \S \ref{sec:cpt_stab}. 
		\end{remark}
		
		Our next main result concerns tilted, or asymmetric, von K\'arm\'an vortex streets translating in the direction $(1,\kpp)$ for small $\kpp>0$, which corresponds to $a$ close to $\pi$. 
		
		\begin{theorem} \label{thm:main-tilt}
			For a fixed $\mu > 0$ and all sufficiently small $\kappa > 0$, let $a_\kappa \in (\pi/2, \pi)$ be the unique number satisfying
			\[
			\sin a_\kappa = \kappa \sinh (2\mu).
			\]
			For all sufficiently small $\varepsilon > 0$, there exists a traveling wave solution to \eqref{eq:Euler_cyl} with velocity $W_\varepsilon(1,\kappa)$, $W_\varepsilon > 0$, of the form
			\[
			\omega_\varepsilon = \varepsilon^{-2} \left( \mathbf 1_{\Omega_{+}^\varepsilon} - \mathbf 1_{\Omega_{-}^\varepsilon} \right)
			\]
			for some disjoint measurable sets $\Omega_{\pm}^\varepsilon \subset \calC$ satisfying $|\Omega_{\pm}^\varepsilon| = \varepsilon^2/2$.
			In the limit $\varepsilon \to 0^+$, we have:
			\begin{itemize}
				\item $W_\varepsilon \to \frac{1}{8\pi} \frac{\sinh (2\mu)}{\cosh (2\mu) - \cos a_\kappa}$,
				\item Up to translations, $\Omega_{\pm}^{\varepsilon}$ shrink to $p_{\pm}$, where $p_+ - p_- = (-a_\kappa,2\mu)$.
			\end{itemize}
			Moreover, $a_\kappa \to \pi$ as $\kappa \to 0$.
		\end{theorem}
		
		The existence is shown again by maximizing the kinetic energy in an admissible class, where the condition $\nrm{x_{2}\omg}_{1} \le \mu$ in \eqref{eq:adm-main} is replaced by an upper bound on a {tilted} moment \begin{equation*}
			\begin{split}
				\int_{[-\pi,\pi) \times \bbR} |x_{2} - \kpp {x}_1| |\omg(x)| dx. 
			\end{split}
		\end{equation*}
        Note that here we are ``unfolding'' the cylinder to identify it with the subset $[-\pi,\pi) \times \bbR$ of the plane.
        To realize the tilted and desingularized von K\'arm\'an vortex street as the \textbf{global} energy maximizer under this moment, we need smallness of $\kpp$.
        For solutions of \eqref{eq:Euler_cyl}, the quantity \begin{equation*}
			\begin{split}
				\int_{[-\pi,\pi) \times \bbR} x_{1}\omg(t,x) dx 
			\end{split}
		\end{equation*}
        is preserved as long as the support of the vorticity (interpreted as a subset of $[-\pi,\pi) \times \bbR$ in $\bbR^2$) does not intersect the ``cut'' $\{ x_{1} = \pi \}$.
        Based on this observation, it might be possible to deduce conditional nonlinear Lyapunov stability for traveling waves from Theorem \ref{thm:main-tilt}, but we do not claim such a result here. 
		
		\subsection{Mathematical literature on von K\'arm\'an street}
		
		The point vortices describing periodic vortex streets \eqref{eq:pv-Karman} were introduced by von K\'arm\'an \cite{Karman1911,Karman1912}.
        Many subsequent works employed point vortex model to study von K\'arm\'an vortex streets in various settings (\cite{Rosenhead1929,Rosenhead1929b,Stremler2003}). Linear instability at the point vortex level was considered in \cite{Jimenez1987,MacKay1987}. Numerical investigations for the vortex patch model for von K\'arm\'an streets were initiated in \cite{SS1,SS2,SS3}. 
        See also \cite{MSS1984} for a numerical study of the linear stability of vortex streets.
		
		While there is a vast literature on von K\'arm\'an vortex streets, there have been relatively few mathematically rigorous works.
        Keady \cite{Keady1985} studied asymptotics for symmetric vortex streets ($a = 0$) in a periodic strip of finite height using variational methods.
        Garcia \cite{Ga20} constructed vortex patch desingularization of \eqref{eq:pv-Karman} in the symmetric cases $a = 0, \pi$, following the bifurcation approach of \cite{HM2017}.
        %Considers generalized SQG equations (for which \eqref{eq:Euler_cyl} is a special case) as well as more general active scalars with radial multipliers. 
		Fan \cite{Fan2022} constructed desingularized solutions in the case $a = 0$.
        Cao--Qin--Zhan--Zou \cite{CQZZ23CVPDE} constructed $C^{1}$-smooth desingularization of \eqref{eq:pv-Karman} for all values of $a \in \bbT$, by solving a nonlinear elliptic equation of the form $\lap \psi_\varepsilon = f(\psi_\varepsilon)$ with prescribed monotone $f$.
        Their method is based on constructing a good approximation for small $\varepsilon > 0$ and applying the implicit function theorem.
        Moreover, there are works constructing hollow von K\'arm\'an vortex streets \cite{Traizet2015,OW} based on the implicit function theorem. 
		
		Somewhat surprisingly, there have been no previous works discussing nonlinear stability of desingularized von K\'arm\'an vortex streets.
        In particular, this seems to be the first work obtaining them as energy maximizers.
        We note that in the special case of $a = 0$, stability within odd-symmetric perturbations (equivalently, for domain $\bbT\times\bbR_+$) would follow from existing literature, but it is a very challenging problem to study stability without the odd symmetry. 
		
		%asymmetric case in the case $a = 0, \pi$, for the whole range $0<s\le 1$.
		%were widely investigated in the applied fluid literature, 

		\subsection{Ideas and difficulties of the proof}
		
		In this section, we discuss technical ideas that go into the proof of our main results, and mention several relevant works along the way. 
		
		Our starting point is the observation that von K\'arm\'an vortex streets are the energy maximizers in the periodic cylinder under natural constraints, at least when they are sufficiently concentrated.
        While Lord Kelvin declared back in 1880 \cite{Kelvin} that such energy maximizers are Lyapunov stable, in our concrete setting there are a few technical difficulties towards showing
        (i) von K\'arm\'an vortex streets are indeed energy maximizers and (ii) nonlinear Lyapunov stability.
        These difficulties come from the fact that these vortices have both positive and negative signs and that the fluid domain is unbounded.
        To the best of our knowledge, this is the first work dealing with these two difficulties at the same time, without imposing any symmetries of the vorticity. 
		
		The classical works of Turkington \cite{Tur1983,Tur1987} give stability of single-signed and concentrated vorticity in bounded two-dimensional domains, using that such configurations maximize the kinetic energy under $L^1$ and $L^\infty$ constraints.
        For bounded domains, such stability was soon extended to the case of mixed sign vorticities in \cite{MaPa}, based on the observation that energy is maximized when vorticities of opposite signs are as much separated as possible.
        Note that for two-dimensional compact manifolds without boundary, vorticity must have zero mean, as a consequence of Stokes' theorem.
        More systematic study of mixed-sign vorticities from calculus of variations point of view was initiated by Cao, Wang, and their collaborators \cite{CaoWangZuo2021,Wang2024b,WangZuo2023,CaoLiWang2025Sphere,CaoWangZhan2020}.
        In various domains, these works obtained steady vortex configurations concentrated near two points with opposite signs, and their nonlinear stability. 
		
		The case of unbounded domains presents significant technical challenges, already at the level of proving the existence of energy maximizers, even when the vorticity is assumed to be odd symmetric in $x_{2}$ and is nonnegative in the upper half-plane $\{ x_{2} \ge 0 \}$.
        (The historical breakthroughs were made in \cite{Bur1987,Bur1988,BNL2013}.)
        These two assumptions are \textit{essential}: the key conserved quantity controlling the vanishing of vorticity to infinity in $x_{2}$ is the impulse $\int x_{2}\omg(x)$, but this is coercive only under the sign condition of the vorticity.
        Furthermore, without forcing the odd symmetry, the vorticities of opposite signs will move away from each other in the $x_{1}$ direction, without changing the value of the impulse. 
		
		Therefore, our key technical contribution lies in the variational characterization of von K\'arm\'an vortex streets while facing two challenges at the same time: mixed-sign vorticity and having an unbounded domain.
        Let us clarify that it is possible to  variationally study von K\'arm\'an vortices while suppressing some, or even all, of these difficulties.
        For one thing, one may force the glide reflection symmetry on $\calC$ together with sign condition $\omg \ge 0$ on $\{x_{2} \ge 0\}$; our method in this paper immediately gives the existence of von K\'arm\'an vortices as the energy maximizer, but then the resulting stability statements hold only under these extra conditions, which are too restrictive.
        If one studies the same problem in the periodic strip $\bbT \times [-L,L]$ for a given $L > 0$, namely when $\bbR$ is replaced by a finite interval, none of these difficulties appear, and the analysis can be done in a parallel manner as in the aforementioned works of Cao, Wang, and others.
        Note however that in $\bbT \times [-L,L]$, one cannot study asymmetric (tilted) von K\'arm\'an vortex streets which do not travel in the horizontal direction. 
		
		\subsubsection{{Technical advances in the current work}}
        To obtain the desingularized von K\'arm\'an vortex patches as the energy maximizer without sign nor symmetry conditions, we need to use the absolute impulse $\int |x_2 \omg(x)|$ rather than the usual impulse $I(\omega)$.
        This allows us to control vorticity far away from $\{ x_2 = 0 \}$, and we conclude in the end that for the maximizer these quantities coincide up to sign.
        A similar idea was recently used in \cite{ACJQ} to obtain stability of the Lamb dipole without sign condition in $\bbR \times \bbR_+$.
        The price we pay for using the absolute impulse is that the energy maximizer is not automatically a traveling wave, unless we show additionally that the vorticity support avoids the line $\{ x_{2} = 0 \}$; see \S \ref{subsubsec:det} below for more details.

        Furthermore, unlike most of the previous works on variational study of mixed sign vorticities, we do not impose any a priori assumptions on the \textit{support of the vorticity} in the admissible class; see Definition~\ref{def:F-and-S}.
        Since one expects the energy maximum for small $\varepsilon > 0$ to approximate \eqref{eq:pv-Karman} (up to translations and reflection), it is convenient to assume from the beginning that in the admissible class, the vorticity is supported in small balls containing the point vortices.
        This approach gives the desingularized vortex streets only as a \textit{local} energy maximizer, which gives additional difficulties in deducing nonlinear stability; for this point, see the discussions in \cite{Burton2005,Wang2024,Wang2024b,CaoWang2019}.
        Moreover, this approach is not appropriate when the support of the maximizer is not highly localized. 
		
		Lastly, note that we are maximizing the energy directly, rather than a \textit{penalized} version of the energy which is done in many of the previous works.
        In the penalized method, one maximizes the quantity which is energy minus a suitable multiple ($\varepsilon$-dependent) of the $L^{1}$ of the vorticity, rather than placing the latter in the admissible class.
        %One may further subtract from the energy  other conserved quantities, including $L^p$ norm and the impulse, with well-designed coefficients. 
		While the penalization method has several technical advantages, its serious drawback appears when $\varepsilon$ is not so small: it could be difficult (or even false) to ensure that zero vorticity is not the penalized maximizer, which is always guaranteed in the direct method (cf. \cite{ACJSW}).
        Note that in \S\ref{sec:prelim} below when we prove the existence of energy maximizer, smallness of $\varepsilon$ is not needed at all.
        Based on this, we plan to address the existence and stability of \textit{large-scale} von K\'arm\'an vortex streets in the future: note that all the existing rigorous works including the current one only concern vortex streets with sufficiently small vortex core (support) radius.
        However, there is a lot of price to be paid to make the direct energy maximization work--compared to the penalization method, it is much more tedious to derive the Euler--Lagrange equations, study the properties of the maximizers, and to prove compactness for maximizing sequences; see the proof summary below.
        We believe that the relevant techniques developed here work in a much more general setting, and hope that it would justify the length of the current paper.  
		
		\subsubsection{{Proof details}}\label{subsubsec:det}
		
		We discuss a few technical details of our approach, focusing on core concentration and conditional stability.
		
		\begin{itemize}
			\item  \emph{Proof of concentration}: We successively extract information from the set of maximizer obtained without the last two conditions in \eqref{eq:adm-main}.
			The table below summarizes the chain of argument in \S\ref{sec:vari} and \S\ref{sec:unif}.

			\begin{center}
				\begin{minipage}{\linewidth}
					\small
					\renewcommand{\arraystretch}{1.2}
					\begin{tabularx}{\linewidth}{
							@{}>{\raggedright\arraybackslash\bfseries}p{0.27\linewidth}
							>{\raggedright\arraybackslash}X@{}
						}
						\toprule
						Existence
						& A maximizer $\omega$ with $\|x_2\omega\|_1 = \mu$. \\
						\addlinespace[0.3em]
						
						  Patch structure
						& $\omega = \varepsilon^{-2}(\mathbf1_{\Omega_+} - \mathbf1_{\Omega_-})$. \\
						\addlinespace[0.3em]
						
						Euler--Lagrange equation
						& $\Omega_\pm = \{\pm\psi[\omega] - W_\pm|x_2| - \gamma_\pm > 0\}$. \\
						\addlinespace[0.3em]
						
						Lagrange multipliers
						& $W_+ = W_- =: W \ge 0$ and $\gamma_+ + \gamma_- \ge 0$, identified by coupled variations. \\
						
						Pohozaev identity
						& $2W\mu + (\gamma_+ + \gamma_-)\|\omega\|_1/2 = E[\omega] + \|\partial_2\psi\|_2^2$. \\
						
						\midrule
						Weak mass saturation
						& $\|\omega_\pm\|_1 \to 1/2$ as $\varepsilon \to 0$. \\
						\addlinespace[0.3em]
						
						Weak confinement
						& $\displaystyle \int_{\{|x_2|\le R\}}\omega_\pm dx \ge c > 0$. \\
						\addlinespace[0.3em]
						
						Estimate of $W$
						& $W = O(1)$ as $\varepsilon \to 0$; subsequently identified with the traveling speed. \\
						\addlinespace[0.3em]
						
						Threshold estimates
						& $\gamma_\pm \sim |\log\varepsilon|$. \\
						
						\midrule
						Mass saturation
						& $\|\omega_\pm\|_1 = 1/2$, forced by $\gamma_+ + \gamma_- > 0$. \\
						\addlinespace[0.3em]
						
						Uniform confinement
						& $\Omega_\pm \subset \{|x_2| \le \bar R\}$;
						outside this strip, $\pm\psi[\omega]$ cannot reach thresholds. \\
						\addlinespace[0.3em]
						
						Core scale
						& $\operatorname{diam}(\Omega_\pm) \lesssim \varepsilon$. \\
						\bottomrule
					\end{tabularx}
					
					\par\smallskip
					{\footnotesize\emph{
							The same scheme extends perturbatively to the similar setting used to construct tilted streets; see \S\ref{sec:tilted_karman}.
						}\par}
				\end{minipage}
			\end{center}
			
			The key ingredient for the proof of concentration to show that $W=O(1)$ as $\varepsilon \to 0$.
            An important intermediate step is to derive the following Pohozaev-type identity from the Euler--Lagrange equation: 
			\[
			W\mu =   \int_{\calC} x_2\omega\partial_2\psi[\omega] dx +  \int_{\Omega_+} \varepsilon^{-2} \Psi_+ dx + \int_{\Omega_-} \varepsilon^{-2} \Psi_- dx, \quad\mbox{with}\quad  \Psi_\pm := (\pm\psi[\omega] - W|x_2| - \gamma_\pm)_+. 
			\]
			The first term can be bounded by total mass $\nrm{\omg}_{1}$ and $\mu$, after a 
			symmetrization which removes the singularity in the kernel.
			The remaining terms must be estimated while allowing the cores to intersect the line $\{x_2 = 0\}$.
			In particular, the weight $|x_2|$ may produce a measure defect, for example,
			\[
			-\Delta\Psi_\pm = \varepsilon^{-2} + 2W \delta_{\{x_2=0\}} \quad \textrm{in }\,\, \Omega_\pm,
			\]
			where $\delta_{\{x_2=0\}}$ denotes the Dirac delta on this line.
			We need to make the contribution from this line source small.
			For this, we consider an auxiliary function $\phi \in 
			H_0^1(\Omega_+)$ that satisfies $-\Delta\phi = \varepsilon^{-2}$ in $\Omega_+$.
			By testing $\Psi_+$ against $\phi$, we isolate the defect caused by $\delta_{\{x_2=0\}}$ as
			\[
			2W\int_{\bbT}\phi(x_1,0)dx_1.
			\]
			By Talenti's comparison principle,
			\[
			\|\phi\|_\infty \lesssim 1, \quad \|\nabla\phi\|_2^2 = \varepsilon^{-2}\int_{\Omega_+}\phi dx \lesssim 1,
			\]
			and then the area bound $|\Omega_+| \lesssim \varepsilon^2 $ makes $\|\phi\|_2$ small enough.
			The trace inequality transfers this smallness to the coefficient of $W$:
			\begin{align*}
				W\int_{\bbT}\phi(x_1,0)dx_1 \lesssim W\left( \|\phi\|_2 \|\partial_2\phi\|_2 \right)^{1/2} 
				\lesssim W|\Omega_+|^{1/4} \|\phi\|_\infty^{1/2}\|\partial_2\phi\|_2^{1/2}
				\lesssim  W\varepsilon^{1/2}.
			\end{align*}
			Thus
			\[
			\int_{\Omega_+} \Psi_+\omega_+ \, dx = \int_{\Omega_+} \varepsilon^{-2} \Psi_+ \, dx \lesssim 1 + W\varepsilon^{1/2},
			\]
			regardless of how $\Omega_+$ intersects $\{x_2=0\}$.
			The same argument for the negative patch yields $W = O(1)$.
			
			Once $W$ is bounded, the thresholds $\gamma_{\pm}$ must carry the logarithmic growth of the energy:
			\[
			|\log \varepsilon| \sim 2E[\omega] = \frac{\|\omega\|_1}{2}(\gamma_+ + \gamma_-) + \underbrace{W\mu+\int_{\calC}(\Psi_+\omega_++\Psi_-\omega_-)dx}_{O(1)}.
			\]
			For the potential to reach these thresholds throughout $\Omega_{\pm}$, most of core mass must lie within $O(\varepsilon)$ of each point, i.e. $\operatorname{diam}(\Omega_{\pm}) = O(\varepsilon)$.

			\item \emph{Proof of stability}: Compactness of maximizing sequences yields stability for a perturbed solution $v(t)$ provided that $\|x_2 v \|_1$ remains controlled.
			For the version under a support restriction, we compare two effects of ``incorrectly'' signed vorticity: its contribution to the excess $\|x_2v \|_1 - I(v)$ is weighted by the distance from the reference line $\{x_2=0\}$, whereas its energy penalty is of order $|\log \varepsilon|$ with respect to unit mass.
			The support bound connects these two quantities.
			
			Indeed, let $v(t) = v_+(t) - v_-(t)$ be a patch with the same distribution as $\omega \in \calS_{\mu, \varepsilon}^{\mathrm{cen}}$, sufficiently close in $L^1$ and impulse.
			Let $R \ge 1$ bound its vertical support diameter.
			Then
			\[
			\|x_2v\|_1 - I(v) = 2\int_{\{x_2<0\}}|x_2|v_+ dx  + 2\int_{\{x_2>0\}}|x_2|v_- dx \lesssim R\cdot m,
			\]
			where $m = \int_{\{x_2<0\}}v_+ dx + \int_{\{x_2>0\}}v_- dx$ is the incorrectly signed mass.
			
			For $R \le c_0|\log \varepsilon|$ with small $c_0$, the estimates $\gamma \sim |\log\varep|$ and $W = O(1)$ give a gap of order $|\log \varepsilon|$ on the incorrectly signed regions:
			\begin{equation*}
				\begin{cases}
					\psi[\omega]-Wx_2 \le \gamma_+-c|\log\varepsilon| & \textrm{on }\,\, \supp v_+\cap\{x_2<0\}, \\
					\psi[\omega]-Wx_2 \ge -\gamma_-+c|\log\varepsilon| & \textrm{on }\,\, \supp v_-\cap\{x_2>0\}.
				\end{cases}
			\end{equation*}
			Rearranging this mass into the corresponding cores gains energy $\gtrsim |\log\varepsilon|m$.
			Although the impulse increases, the rearranged patch has energy at most $E[\omega] + C|I(v)-\mu| + O(Rm)$, and hence
			\[
			m|\log\varepsilon| \lesssim E[\omega] - E[v] + C|I(v)-\mu|.
			\]
			Together with $\|x_2v\|_1 - I(v)\lesssim Rm$, we derive
			\[
			0 \le \|x_2v\|_1 - I(v) \lesssim \frac{R}{|\log\varepsilon|} \big( E[\omega]-E[v]+C|I(v)-\mu| \big).
			\]
			Conservation of energy and impulse therefore provides the control of $\|x_2v\|_1$ needed to apply compactness.
		\end{itemize}
		
		\subsubsection{Potential extensions}\label{subsubsec:circ} Let us list a few extensions of the current work which seem to be within reach of methods developed here.
        As we mentioned earlier, the case of $\bbT \times [-L,L]$ can be studied in parallel.
        Moreover, in the case of $\calC$, one can study von K\'arm\'an vortex streets with {arbitrary total circulation}: zero circulation condition is an obstruction, rather than a simplification in all of our proofs.

		\subsection{Organization of the paper}
		The rest of the paper is organized as follows. \begin{itemize}
			\item In \S \ref{sec:prelim}, we collect several basic estimates in the cylinder. 
			\item In \S \ref{sec:vari}, we introduce the admissible classes and derive the Euler--Lagrange equation satisfied by maximizers.
            Preliminary estimates for the Lagrange multipliers, including a Pohozaev-type identity are obtained in this section as well. 
			\item More detailed quantitative estimates for the multipliers and the maximizers are then obtained in \S \ref{sec:unif}. 
			\item The asymptotic information of the maximizers as $\varepsilon\to 0$ is given in \S \ref{sec:karman}, which reveals that the maximizers converge to the point vortices \eqref{eq:pv-Karman}. 
			\item We prove the compactness of maximizing sequences and stability of the maximizing set in \S \ref{sec:cpt_stab}.
            This concludes in particular Theorems \ref{thm:max}, \ref{thm:main-compact}, and \ref{thm:main-stab}. 
			\item Lastly in \S \ref{sec:tilted_karman}, we study asymmetric vortex streets, proving Theorem \ref{thm:tilt_small} which covers Theorem \ref{thm:main-tilt}. 
		\end{itemize}
		
		\subsection*{Acknowledgments} IJ was supported by the NRF grant from the Korea government No. 2022R1C1C1011051, RS-2024-00406821, and by a KIAS Individual Grant.
        GQ was supported by National Key R\&D Program of China (Grant 2025YFA1018400) and NNSF of China (Grant 12471190).
        KW was supported by SNF Project 212573 FLUTURA – Fluids, Turbulence, Advection.
		
		\subsection*{AI Declarations} We declare that no AI tools were used in any part of this work.

		\section{Preliminaries}\label{sec:prelim}
		
		\subsection{Basic estimates}
		
		\subsubsection{Green's function}
		
		\begin{lemma}\label{lem_green_local}
			Let $r_0 \in ( 0, 1) $ be sufficiently small.
			Then there exists a constant $C = C(r_0) > 0$ such that
			\begin{equation}\label{green_local}
				\left|G(x-y) + \frac{1}{2\pi} \log\left| x-y \right| \right| \le C
			\end{equation}
			for any $x, y \in \calC$ such that $0 < |x-y| \le r_0$.
		\end{lemma}
		
		\begin{proof}
			By translation invariance, we may take $y=0$ and then
			\[
			G(x, 0) = -\frac{1}{4\pi}\log \left( 2\cosh  x_2   - 2\cos   x_1  \right).
			\]
			Since $2\cosh  x_2   - 2\cos   x_1 = |x|^2 + O(|x|^4)$ for $|x| \le r_0 < 1/2$, \eqref{green_local} follows.
		\end{proof}
		
		\begin{lemma}\label{lem_green_decom}
			For any $x, y \in \calC$ such that $ x \neq y$, 
			\begin{equation}\label{green_lin}
				G(x-y) =  H(x-y) -\frac{1}{4\pi}|x_2 - y_2|,
			\end{equation}
			where 
			\begin{equation}\label{green_sing}
				H(x-y) := - \frac{1}{4\pi} \log \left( 1 - 2e^{-|x_2 - y_2|}\cos \left(  x_1 - y_1  \right) + e^{-2 |x_2-y_2|} \right).
			\end{equation}
			Moreover, the singular term $H$ satisfies 
			\begin{equation}\label{green_sing_2}
				-\frac{1}{2\pi} \log 2 \le H(x-y) \le -\frac{1}{2\pi} \log \left(1 - e^{- |x_2 - y_2|} \right).
			\end{equation}
			In particular, for some absolute constant $C > 0$,
			\begin{equation}\label{green_sing_3}
				H(x-y) \le \frac{1}{2\pi} \log_+ \frac{1}{|x-y|} + C,
			\end{equation}    
			where $\log_+s := \max\{ \log s, 0 \}$ for $s > 0$.    
		\end{lemma}
		
		\begin{proof}
			From the identity
			\[
			2\cosh (x_2 - y_2) - 2\cos(x_1 - y_1)   =  e^{ |x_2 - y_2|} \left( 1 - 2e^{- |x_2 - y_2|}\cos (x_1 - y_1) + e^{-2 |x_2 - y_2|} \right),
			\]
			\eqref{green_lin} is directly follows.
			Since $(1 + e^{-|x_2 - y_2|})^2 \ge (1 - e^{- |x_2 - y_2|})^2$, we get \eqref{green_sing_2} from the expression \eqref{green_sing}.
			To prove \eqref{green_sing_3}, we first note that, for $t = |x_2 - y_2| \in [0, 1]$ and $s = x_1 - y_1 \in [-\pi, \pi)$,
			\[
			1-e^{-t}\ge t/2, \quad  1 - \cos s\ge \frac{2s^2}{\pi^2}.
			\]
			Thus, if $|x - y| \le 1$,
			\[
			1 - 2e^{-|x_2 - y_2|}\cos \left(  x_1 - y_1  \right) + e^{-2 |x_2-y_2|} = (1 - e^{-t})^2 + 2e^{-t}(1 - \cos s) \gtrsim   |x-y|^2,
			\]
			whereas the left-hand side is bounded away from zero when $|x-y| > 1$.
			Substituting in \eqref{green_sing} yields \eqref{green_sing_3}.
		\end{proof}
		
		\subsubsection{Kinetic energy}
		
		Although the kinetic energy is naturally associated with zero-circulation vorticities on $\calC$, we shall use the same notation even if $\int_{\calC}\omega \neq 0$:
		\[
		E[\omg]:=\frac{1}{2} \iint_{\calC\times\calC}G(x-y)\omega(x)\omega(y)\,dxdy
		\]
		whenever it is well-defined.
		We define the interaction energy between two vorticities by \begin{equation}\label{eq:E-int}
			\begin{split}
				E_{\mathrm{int}}[\omega, \eta] :=  \iint_{\calC \times \calC} G(x-y) \omega(x) \eta(y) dx dy
			\end{split}
		\end{equation} so that $
		E[\omega]=\frac12 E_{\mathrm{int}}[\omega,\omega]$
		and $
		E[\omega + \eta] = E[\omega] + E[\eta] + E_{\mathrm{int}}[\omega,\eta].$

		\begin{lemma}\label{lem:en_pos}
			Let $A \subset \calC$ be a measurable set with finite measure, $\lambda > 0$, and consider the patch $\omega = \lambda \mathbf{1}_A$. Assume $\|x_2 \omega\|_{L^1(\calC)} < \infty$.
			Then there exists a constant $C > 0$ such that
			\begin{equation} \label{eq:en_pos}
				-C \left( \|\omega\|_1 \|x_2\omega\|_1 + \|\omega\|_1^2 \right) \le E[\omega] \le \frac{1}{8\pi} \|\omega\|_1^2 \log_+ \left( \frac{\lambda}{\|\omega\|_1} \right) + C \|\omega\|_1^2.
			\end{equation}
		\end{lemma}
		
		\begin{proof}
			By \eqref{green_lin} in Lemma~\ref{lem_green_decom}, we split the energy as
			\[
			E[\omega] = \frac{1}{2} \int_{\calC} \int_{\calC} H(x-y) \omega(x) \omega(y) \,dx dy - \frac{1}{8\pi} \int_{\calC} \int_{\calC} |x_2 - y_2| \omega(x) \omega(y) \,dx dy.
			\]
			Since $\omega \ge 0$, the second term is non-positive.
			Using \eqref{green_sing_3} and choosing a suitable cut of $\calC$, we have
			\[
			E[\omega] \le \frac{1}{4\pi} \int_{\bbR^2} \int_{\bbR^2} \log_+ \frac{1}{|x-y|} (\omega \mathbf{1}_{\calC}) (x) (\omega \mathbf{1}_{\calC})(y) \,dx dy + C \|\omega\|_1^2.
			\]
			By the Riesz rearrangement inequality,
			\[
			\iint_{\bbR^2 \times \bbR^2} \log_+\frac1{|x-y|}\omega(x)\omega(y)\,dxdy
			\le
			\iint_{\bbR^2 \times \bbR^2} \log_+\frac1{|x-y|}\omega^*(x)\omega^*(y)\,dxdy,
			\]
			where $\omega^* = \lambda \mathbf{1}_{D_R}$ is the symmetric decreasing rearrangement of $\omega = \lambda \mathbf{1}_A$ with radius $R = \sqrt{\|\omega\|_1 / (\pi \lambda)}$.
			We estimate the right-hand side depending on the size of $R$.
			
			If $R \ge 1$, then $\lambda / \|\omega\|_1 \le 1/\pi < 1$.
			Since $\log_+ \frac{1}{|x-y|}$ is supported on $|x-y| \le 1$, we see
			\begin{align*}
				\iint_{D_R \times D_R}  \log_+ \frac{1}{|x-y|} \omega^*(x) \omega^*(y) \,dx dy &\le \lambda^2 \int_{D_R} \left( \int_{B_1(x)} \log \frac{1}{|x-y|} \,dy \right) dx \\
				&= \lambda^2 |D_R| \int_{B_1(0)} -\log|z| \,dz \lesssim \|\omega\|_1^2,
			\end{align*}
			where we used $\lambda |D_R| = \|\omega\|_1$ and $\lambda \le \pi^{-1}\|\omega\|_1$ in the last line.
			Thus, the upper bound \eqref{eq:en_pos} holds when $R \ge 1$.
			
			For $R < 1$, we have $|x-y|  < 2$ for all $x, y \in D_R$.
			In this range, the inequality $\log_+ \frac{1}{|x-y|} \le \log \frac{1}{|x-y|} + \log 2$ holds.
			Since 
			\[
			\iint_{D_R \times D_R} -\log|x-y| \,dx dy = |D_R|^2 \left( -\log R + \frac{1}{4} \right),
			\]
			we obtain
			\begin{align*}
				&\iint_{D_R \times D_R} \log_+ \frac{1}{|x-y|} \omega^*(x) \omega^*(y) \,dx dy \le \iint_{D_R \times D_R} \left(-\log|x-y| + \log 2\right) \lambda^2 \,dx dy \\
				&\quad = \lambda^2 |D_R|^2 \left( -\log R + \frac{1}{4} + \log 2 \right)  = \frac{1}{2} \|\omega\|_1^2 \log \left( \frac{\lambda}{\|\omega\|_1} \right) + C \|\omega\|_1^2.
			\end{align*}
			Dividing by $4\pi$ yields the upper bound.
			For the lower bound, we use \eqref{green_lin} and \eqref{green_sing_2} in Lemma~\ref{lem_green_decom}.
			Since $H(x-y) \ge -\frac{\log 2}{2\pi}$ and $|x_2 - y_2| \le |x_2| + |y_2|$, we obtain
			\begin{align*}
				E[\omega] \ge -\frac{\log 2}{4\pi} \|\omega\|_1^2 - \frac{1}{8\pi} \int_{\calC} \int_{\calC} (|x_2| + |y_2|) \omega(x) \omega(y) \,dx dy  \ge -C \|\omega\|_1^2 - \frac{1}{4\pi} \|\omega\|_1 \|x_2 \omega\|_1.
			\end{align*}
			This completes the proof.
		\end{proof}
		
		\begin{lemma}\label{lem:en_cyl}
			Let $\omega, \eta \in L^1\cap L^{\infty}(\calC)$ satisfy $\|x_2\omega\|_{L^1(\calC)}, \|x_2\eta\|_{L^1(\calC)} < \infty$.
			For any sufficiently small $r>0$, we have
			\begin{equation}
				\label{en_bi}
				\begin{aligned}
					\left|E_{\mathrm{int}}[\omega, \eta] \right| \lesssim &\|x_2\omega\|_1\|\eta\|_1 + \|\omega\|_1\|x_2\eta\|_1 + \|\omega\|_{\infty} \|\eta\|_1 r^2 (1 + |\log r|) + \left(1 + |\log r| \right)\|\omega\|_1 \|\eta\|_1.
				\end{aligned}
			\end{equation}
			In particular, 
			\begin{equation}
				\label{energy}
				|E[\omega]| \lesssim  \|\omega\|_1 \|x_2\omega\|_1 + \|\omega\|_1^2 \left(1 + \log_+\frac{\|\omega\|_{\infty}}{\|\omega\|_1} \right) ,
			\end{equation}
			where $\log_+s := \max\{\log s, 0\}$ for $s > 0$.
		\end{lemma}
		
		\begin{proof}
			By using \eqref{green_lin}, we decompose
			\begin{align*}
				E_{\mathrm{int}}[\omega, \eta] &=  \int_{\calC} \int_{\calC} \left( -\frac{1}{4\pi} |x_2 - y_2| \right) \omega(x) \eta(y) \,dxdy +  \int_{\calC} \int_{\calC} H(x-y) \omega(x)\eta(y) \,dxdy  =: E_{\mathrm{lin}}[\omega, \eta] + E_{H}[\omega, \eta].
			\end{align*}
			It is clear that
			\begin{equation}
				\label{eq_tail}
				|E_{\mathrm{lin}}[\omega, \eta]| \lesssim \iint \left( |x_2 | + | y_2| \right) |\omega(x)| |\eta(y)| \,dxdy = \|x_2\omega\|_1\|\eta\|_1 + \|\omega\|_1\|x_2\eta\|_1.
			\end{equation}
			For $E_H[\omega, \eta]$, we further decompose it in the following way: for $r>0$,
			\[
			E_H[\omega, \eta] = \iint_{|x-y| \le r} (\cdots) + \iint_{|x-y| > r} (\cdots) =: I_{\mathrm{near}} + I_{\mathrm{far}}.
			\]
			We fix $ 0 < r_0 \ll 1$ and let $r < r_0$.
			Then by Lemma~\ref{lem_green_local} and \eqref{green_lin} of Lemma~\ref{lem_green_decom}, it holds that
			\begin{equation}
				\label{H_sing}
				|H(x-y)| \lesssim_{r_0} 1 + |\log|x-y|| \quad \textrm{ if }\,\, |x-y| \le r,
			\end{equation}
			and then, we see
			\begin{align*}
				|I_{\mathrm{near}}| & \lesssim  \int_{\calC} |\eta(y)| \left( \int_{|x-y| \le r} \left( 1 + \left |\log \left| x - y \right| \right| \right) |\omega(x)| dx \right) dy \\
				& \le  \|\eta\|_1 \|\omega\|_{\infty} \left( \int_{|z| \le r} 1 + |\log|z|| dz \right)  \simeq  \|\eta\|_1 \|\omega\|_{\infty} r^{2}(1 + |\log r|).
			\end{align*} 
			On the other hand, by \eqref{green_sing} and \eqref{green_sing_2}, we see that $|I_{\mathrm{far}}| \lesssim \left(1 + |\log r| \right)\|\omega\|_1 \|\eta\|_1$, and therefore,
			\[
			|E_H[\omega, \eta]| \le |I_{\mathrm{near}}| + |I_{\mathrm{far}}| \lesssim \|\omega\|_{\infty} \|\eta\|_1 r^2 (1 + |\log r|) + \left(1 + |\log r| \right)\|\omega\|_1 \|\eta\|_1.
			\]
			This with \eqref{eq_tail} proves \eqref{en_bi}.
			By taking $\omega = \eta$ in \eqref{en_bi} and  using $2E[\omega] = E_{\mathrm{int}}[\omega, \omega]$, we obtain \eqref{energy} after choosing $r = r_0 \wedge \sqrt{\| \omega \|_1 / \| \omega \|_{\infty} }$.
		\end{proof}

		\begin{remark}\label{rmk:H_unif_1}
			One verifies that $H(x-\cdot)\in L^1(\calC)$ uniformly in $x\in \calC$.
			Indeed, $|H(x-y)| \lesssim 1 + \log_+(1/|x-y|)$ near the singularity, and \eqref{green_sing} gives $|H(x-y)| \lesssim e^{-|x_2-y_2|}$ away from singularity.
		\end{remark}
		
		\begin{lemma}\label{rem_lin}
			For $\omega$ satisfying $\int_{\calC}\omega = 0$ and $\omega, x_2\omega \in L^1(\calC)$, we have
			\begin{equation}
				\label{eq_en_lin}
				E_{\mathrm{lin}}[\omega]:= -\frac{1}{8\pi} \int_{\calC}\int_{\calC} |x_2 - y_2| \omega(x) \omega(y) dxdy = \frac{1}{4\pi} \int_{\bbR} |A[\omega](\tau)|^2 d\tau,
			\end{equation}
			where
			\begin{equation}
				\label{A}
				A[\omega](\tau) =  \int_{x_2 \ge \tau } \omega(x)dx = \int_\tau^{\infty} \int_{\bbT} \omega(x_1, x_2) \,dx_1dx_2.
			\end{equation}    
		\end{lemma}
		
		\begin{proof}
			By direct computation, we see that $-4\pi E_{\mathrm{lin}}[\omega]$ is equal to
			\begin{align*}
				&\iint_{x_2 > y_2} (x_2 - y_2) \omega(x) \omega(y) \,dxdy =  \iint \omega(x) \omega(y) \left( \int_{\bbR} \mathbf{1}_{y_2 < \tau}\mathbf{1}_{x_2 > \tau} d\tau \right) dxdy \\
				&\quad =  \int_{\bbR} \left( \iint \omega(x)\mathbf{1}_{x_2 > \tau}\, \omega(y)  \mathbf{1}_{y_2 < \tau} dxdy \right)d\tau
				=  \int_{\bbR} \left( \int_{x_2 > \tau} \omega(x)   dx \right) \left( \int_{y_2 < \tau}  \omega(y)  dy \right)  d\tau.
			\end{align*}
			Since $\int_{\calC} \omega =0$, we have $\int_{y_2 < \tau} \omega(y)  = -\int_{y_2 > \tau}\omega(y)  = -A[\omega](\tau)$ and then \eqref{eq_en_lin} follows.
		\end{proof}

		\subsubsection{Stream function}
		
		\begin{lemma}\label{lem_str_bdd}
			Let $\omega \in L^1 \cap L^{\infty}(\calC)$ satisfy $x_2 \omega \in L^1(\calC)$ and $\int_{\calC}\omega = 0$.
			Then, the stream function $\psi[\omega]$ satisfies 
			\begin{equation}\label{str_pt_op}
				\|\psi\|_{\infty} \lesssim \|\omega\|_{1} \left(1 + \log_+\frac{\|\omega\|_{\infty}}{\|\omega\|_{1}} \right) + \|x_2\omega\|_{1}.
			\end{equation}
		\end{lemma}
		
		\begin{proof}
			Using \eqref{green_lin}, we decompose $\psi$ into two parts:
			\[
			\psi(x) = \int_{\calC}H(x-y) \omega(y) dy - \frac{1}{4\pi} \int_{\calC} |x_2 - y_2|\omega(y)dy =: \psi_{H}(x) + \psi_0(x).
			\]
			Since $\int_{\calC}\omega = 0$, 
			\[
			4\pi |\psi_0(x)| = \left| \int_{\calC} \left( |x_2 - y_2| - |x_2| \right)\omega(y)dy \right| \le \int_{\calC} |y_2| |\omega(y)|dy.
			\]
			To estimate $\psi_{H}$, fix $ r \in (0, 1/2)$ and further split
			\[
			\psi_{H}(x) = \int_{|x-y| < r} H(x-y) \omega(y) \,dy + \int_{|x-y| \ge r} H(x-y) \omega(y) \,dy.
			\]
			Using \eqref{green_local} and \eqref{green_lin},
			\[
			\int_{|x-y| < r} |H(x-y)| |\omega(y)| dy \lesssim \|\omega\|_{\infty} \int_{|y| < r} 1 + |\log|y|| dy \lesssim \|\omega\|_{\infty} r^2(1 + |\log r|).
			\]
			On the region $\{|x-y|\ge r\}$, Lemma~\ref{lem_green_decom} gives $|H(x-y)|\lesssim 1+|\log r|$, and hence
			\[
			\int_{|x-y|\ge r}|H(x-y)||\omega(y)|\,dy \lesssim (1+|\log r|)\|\omega\|_1.
			\]
			Combining the estimates for $\psi_0$ and $\psi_H$, we obtain
			\[
			\|\psi\|_{\infty} \lesssim   r^2(1 + |\log r|) \|\omega\|_{\infty}  +  (1 + |\log r|)\|\omega\|_1 + \|x_2\omega\|_1.
			\]
			Then, \eqref{str_pt_op} follows by choosing $r= \frac{1}{2} \wedge \sqrt{\|\omega\|_1/\|\omega\|_{\infty}}$.
		\end{proof}
		
		\begin{lemma}\label{lem_L2}
			Let $\omega \in L^1\cap L^{\infty}(\calC)$ satisfy $x_2\omega \in L^1(\calC)$ and $\int_{\calC}\omega = 0$.
			Then $\nabla \psi \in L^2(\calC)$ and
			\[
			0\le \frac{1}{2}\int_{\calC} |\nabla \psi|^2 dx = \frac{1}{2}\int_{\calC} \psi \omega dx = E[\omega].
			\]
		\end{lemma}
		
		\begin{proof}
			By Lemma~\ref{lem_str_bdd}, we have $\psi \in L^{\infty}(\calC)$.
			Let $\chi \in C_c^{\infty}(\bbR)$ satisfy $\chi \equiv 1$ on $[-1, 1]$, and define $\chi_R(x_2) = \chi(x_2/R)$.
			Testing $-\Delta \psi = \omega$ against $\chi_R^2 \psi$ and integrating by parts yields
			\begin{equation}
				\label{ell_local}
				\int_{\calC} \chi_R^2 |\nabla \psi|^2 dx = \int_{\calC} \chi_R^2 \, \psi \, \omega \,dx - 2 \int_{\calC} \chi_R \psi \, \nabla\psi \cdot \nabla \chi_R\,dx.
			\end{equation}
			By Young's inequality with the property of $\chi_R$, we see
			\[
			\int_{\calC} \chi_R^2 |\nabla \psi|^2 dx \lesssim \int_{\calC} \chi_R^2 \, \psi \, \omega \,dx + \frac{1}{R} \le \|\psi\|_{\infty}\|\omega\|_1 + \frac{1}{R}.
			\]
			This proves the lemma.
		\end{proof}

		\subsection{Existence of maximizer}\label{subsec:max-prelim}
		Let $M_0, M_1, M_2 > 0$ be constants and define
		\[
		\calA := \left\{ \omega \in L^1\cap L^{\infty}(\calC): \|\omega\|_{1} \le M_0, \quad \|\omega\|_{\infty} \le M_1, \quad \|x_2\omega\|_1 \le M_2, \quad \int_{\calC}\omega(x)\, dx = 0 \right\}.
		\]
		In the rest of this section, we prove the existence of an energy maximizer in the class $\calA$. 
		
		\begin{proposition}\label{prop:existence}
			There exists $\omega^* \in \calA$ such that $0 < E[\omega^*] = \sup_{\omega \in \calA} E[\omega] < \infty.$
		\end{proposition}

		\begin{remark}\label{rem:ex}
			The vertical absolute impulse bound $\|x_2\omega\|_1 \le \mu$ is essential: to see this, 
			fix a set $M = [0, 1] \times [0, 1/4] \subset \calC$ and define $\omega_L := 2\,\mathbf 1_{M} - \mathbf 1_{M+(0,L)} - \mathbf 1_{M+(0,-L)}$ for some $L > 1.$ 
			We see that $\|\omega_L\|_{1} + \|\omega_L\|_{\infty} \lesssim 1$, $\int \omega_L  = \int x_2 \omega_L  = 0$, but $E_{\mathrm{lin}}[\omega_L] \gtrsim L \to \infty$ as $L \to \infty$, while $E[\cdot] = E_{\mathrm{lin}}[\cdot] + E_H[\cdot]$ and $|E_H[\omega_L]| \lesssim \|\omega_L\|_{L^1 \cap L^{\infty}} \lesssim 1$.
			Note that $\|x_2\omega_L\|_1 \sim L$.
		\end{remark}

		\begin{lemma}
			\label{lem_conv_lin}
			Let $\omega^* \in \calA$ and let $\{ \omega_n \}_n \subset \calA$ be a sequence such that $\omega_n \rightharpoonup^* \omega^*$ in $L^{\infty}(\calC)$.
			Then, $E_{\mathrm{lin}}[\omega_n] \to E_{\mathrm{lin}}[\omega^*]$ as $n \to \infty$, where $E_{\mathrm{lin}}[\cdot]$ is defined by \eqref{eq_en_lin}.
		\end{lemma}
		
		\begin{proof}
			By \eqref{eq_en_lin} in Lemma~\ref{rem_lin}, it suffices to prove $A[\omega_n] \to A[\omega^*]$ in $L^2(\bbR)$, where $A[\eta](\tau) = \int_{x_2 \ge \tau} \eta  $.
			We first note that, for each $\tau \in \bbR$ and $R > |\tau|$,
			\begin{align*}
				A[\omega_n](\tau)  &= \int_{\calC} \omega_n(x) \mathbf{1}_{R > x_2 > \tau} dx + \int_{\calC} \omega_n(x) \mathbf{1}_{x_2 > R} dx,
			\end{align*}
			where
			\[
			\left| \int_{\calC} \omega_n(x) \mathbf{1}_{x_2 > R} dx \right| \le \frac{M_2}{R} \quad \textrm{and} \quad \lim_{n \to \infty} \int_{\calC} \omega_n(x) \mathbf{1}_{R > x_2 > \tau} dx = \int_{\calC} \omega^*(x) \mathbf{1}_{R > x_2 > \tau} dx.
			\]
			By taking the limit as $n \to \infty$ followed by $R \to \infty$, we have $A[\omega_n](\tau) \to A[\omega^*](\tau)$ for each $\tau \in \bbR$.
			Moreover,  it is straightforward to check
			\[
			|A[\omega_n](\tau)| \le M_0\mathbf{1}_{|\tau| < 1} +  |\tau|^{-1}M_2 \mathbf{1}_{|\tau| \ge 1} \in L^2(\bbR),
			\]
			where we use the relation $A[\omega_n](\tau) = - \int_{x_2 \le \tau} \omega_n\,dx$ for the case $\tau < -1$.
			By dominated convergence theorem, this completes the proof.
		\end{proof}
		
		\begin{lemma}\label{lem:weak_energy_conv}
			Let $\{\omega_n\}_n \subset \calA$.
			Then there exist a subsequence, still denoted by $\omega_n$, and a function $\omega^* \in \calA$ such that $\omega_n \rightharpoonup^* \omega^*$ in $L^{\infty}(\calC)$ and $E[\omega_n]\to E[\omega^*]$ as $n \to \infty$.
		\end{lemma}
		
		\begin{proof}
			Since $\sup_n \|\omega_n\|_{\infty} \le M_1 < \infty$, there exist a subsequence, still denoted by $\omega_n$, and a function $\omega^* \in L^{\infty}(\calC)$ such that $\omega_n \rightharpoonup^* \omega^*$ in $L^{\infty}(\calC)$.
			By standard properties of weak-$^*$ convergence with the help of $\|x_2\omega_n\|_1 \le M_2$, one verifies $\omega^* \in \calA$.
			Set
			\[
			E[\omega] = E_{\mathrm{lin}}[\omega] + E_H[\omega], \quad
			E_H[\omega] := \frac12\iint H(x-y)\omega(x)\omega(y)\,dxdy = \frac{1}{2}E_H[\omega, \omega],
			\]
			where $H$ is defined by \eqref{green_sing} in Lemma~\ref{lem_green_decom}.
			Then by Lemma~\ref{lem_conv_lin}, it suffices to prove
			\begin{equation}
				\label{eq_H_conv}
				\lim_{n \to \infty} \iint H(x-y) \,\omega_n(x) \,\omega_n(y)\,dxdy = \iint H(x-y) \,\omega^*(x) \,\omega^*(y)\,dxdy.
			\end{equation}
			To see this, we set $\psi_H[\omega] := H \ast \omega$ and decompose
			\begin{equation}
				\label{eq_H_conv_2}
				\begin{aligned}
					|E_H[\omega^*] - E_H[\omega_n] | \le &\left| \int_{\calC}\psi_H[\omega^*] (\omega^* - \omega_n) dx \right| + \left| \int_{\calC} (\psi_{H}[\omega_n] - \psi_{H}[\omega^*])  \omega_n dx \right|   =: I_n + II_n.
				\end{aligned}
			\end{equation}
			By Remark~\ref{rmk:H_unif_1}, we have $H(x - \cdot) \in L^1(\calC)$ and this implies $\left\| \psi_H[\omega^*] \right\|_{L^1(\calC)} \le M_0\|H\|_1 < \infty$ and thus $I_n \to 0$ by the weak-$^*$ convergence.
			For $II_n$, we further split the integral:
			\begin{align*}
				II_n \le \int_{|x_2| > R} |\psi_{H}[\omega_n] - \psi_{H}[\omega^*]| |\omega_n| \,dx + \int_{|x_2| \le R} |\psi_{H}[\omega_n] - \psi_{H}[\omega^*]| |\omega_n| \,dx =: II_n^1 + II_n^2.
			\end{align*}
			Since $\|\psi_{H}[\omega_n]\|_{\infty}, \|\psi_{H}[\omega^*]\|_{\infty} \le M_1 \|H\|_1$, it holds that
			\[
			II_n^1 \le \|\psi_H[\omega_n] - \psi_H[\omega^*]\|_{\infty}  \int_{ |x_2| > R} \frac{|x_2|}{R}|\omega_n| \,dx \le 2M_1\|H\|_1\frac{M_2}{R}.
			\]
			For $II_n^2$, note that $\psi_{H}[\omega_n](x) \to \psi_{H}[\omega^*](x)$ for each $x \in \calC$ and
			\begin{align*}
				\left| \psi_{H}[\omega_n](x+h) - \psi_{H}[\omega_n](x) \right| = &\left| \int_{\calC} \left( H(x+h-y) - H(x-y) \right) \omega_n(y) \,dy \right| \\
				\le & M_1 \int_{\calC} \left| H(z + h) - H(z) \right|  \,dz \to 0 \quad \textrm{ as }\, h \to 0.
			\end{align*}
			Therefore, by the Arzela--Ascoli theorem, $\psi_{H}[\omega_n](x) \to \psi_{H}[\omega^*](x)$ uniformly on $\{x : |x_2| \le R\}$.
			Consequently, $ II_n^2 \le  M_0 \left\|\psi_{H}[\omega_n] - \psi_{H}[\omega^*] \right\|_{L^{\infty}(\{|x_2| \le R\})} \to 0$ as $n \to \infty$.
			Letting $n \to \infty$ and then $R \to \infty$ in \eqref{eq_H_conv_2}, we obtain \eqref{eq_H_conv}.
		\end{proof}

		\begin{proof}[\textbf{Proof of Proposition \ref{prop:existence}}]
			By Lemma~\ref{lem:en_cyl}, the supremum is finite.
			Moreover, choosing a nontrivial vorticity in $\calA$, Lemma~\ref{lem_L2} gives a strictly positive energy.
			Hence $0 < \sup_{\omega\in\calA}E[\omega]=:E^*<\infty$.
			Let $\{\omega_n\}_n \subset  \calA$ be a maximizing sequence such that $\lim_{n \to \infty} E[\omega_n] = E^*$.
			Then the result directly follows from Lemma~\ref{lem:weak_energy_conv}.
		\end{proof}

		\section{Variational problem for desingularized von K\'arm\'an vortex}\label{sec:vari}
		
		In this section, we introduce the small parameter $\varepsilon > 0$ and introduce the admissible class $\calF_{\mu, \varepsilon}$ as well as the absolute impulse saturating maximizing set $\calS_{\mu, \varepsilon}$, where $\mu > 0$ is fixed in the remainder of the paper. 
		
		\begin{definition}[Admissible class and Maximizing set]\label{def:F-and-S}
			For $\varepsilon > 0$, we define the admissible class
			\begin{align*}
				\calF_{\mu, \varepsilon} &:= \left\{ \omega \in L^1\cap L^{\infty}(\calC): \|\omega\|_1 \le 1, \|\omega\|_{\infty} \le \varepsilon^{-2}, \|x_2\omega\|_1 \le \mu,\int_{\calC}  \omega(x)\,dx =0 \right\} 
			\end{align*}
			and the $\mu$-saturated set of maximizers as
			\begin{align*}
				\calS_{\mu, \varepsilon} &:= \left\{ \omega   \in \calF_{\mu, \varepsilon} \, : \, E[\omega] = \sup_{\theta \in \calF_{\mu, \varepsilon}}E[\theta], \,  \|x_2\omega\|_1 = \mu \right\}.
			\end{align*}
			Observe that $\calS_{\mu, \varepsilon}$ is \textbf{nonempty}: applying Proposition~\ref{prop:existence} with $M_0 = 1$, $M_1 = \varepsilon^{-2}$, and $M_2 = \mu$, 
			\begin{align*}
				\tld{\calS}_{\mu, \varepsilon} &:= \{ \omega \in \calF_{\mu, \varepsilon}: E[\omega] = \sup_{\theta \in \calF_{\mu, \varepsilon}}E[\theta] \} \neq \emptyset.
			\end{align*} Next, note that $E$, $\|\cdot\|_1$, $\|\cdot\|_\infty$, and the zero-circulation constraint are invariant under translations.
			Therefore, if $\omega \in \tld{\calS}_{\mu, \varepsilon}$ satisfies $\|x_2\omega\|_1 < \mu$, then a vertical translation of $\omg$ can saturate the constraint $\mu$:  $\omega(x_1, x_2 - c) \in \tld{\calS}_{\mu, \varepsilon}$ and $\|x_2 \omega(x_1, x_2 - c)\|_1 = \mu$ for some $c \in \bbR$.
		\end{definition}
		
		The main result of this section is the \textbf{Euler--Lagrange equation} for elements in $\calS_{\mu, \varepsilon}$: \begin{proposition}[Euler--Lagrange equation]\label{prop:EL}
			There exists $\varepsilon_0 = \varepsilon_0(\mu)$ such that the following holds:
			For any $\varepsilon \in (0, \varepsilon_0)$ and for any $\omega \in \calS_{\mu, \varepsilon}$, there exist unique constants $\gamma_+ = \gamma_+(\omega), \gamma_- = \gamma_-( \omega)$ with $\gamma_+ + \gamma_- \ge 0$ and $W = W(\omega) \ge 0$ such that 
			\begin{equation}\label{eq:EL}
				\omega = \varepsilon^{-2} \left(\mathbf{1}_{ \{\psi[\omega] - W|x_2| - \gamma_+ > 0 \} } - \mathbf{1}_{ \{ -\psi[\omega] - W|x_2| - \gamma_- > 0 \} } \right).
			\end{equation}
		\end{proposition}
		
		This result is proved after deriving a sequence of lemmas in \S \ref{subsec:EL}. Then, in \S \ref{subsec:Po} we derive an identity satisfied by the Lagrange multipliers.

		\subsection{Deriving the Euler--Lagrange equation}\label{subsec:EL}
		We now derive the Euler--Lagrange equation for $\calS_{\mu, \varepsilon}$.

		\begin{lemma}[Energy lower bound]\label{lem:en_low}
			There exist $\varepsilon_0 = \varepsilon_0(\mu)$ and $C = C(\mu) > 0$ such that the following holds:
			For any $\varepsilon \in (0, \varepsilon_0)$ and any maximizer $\omega \in \calS_{\mu, \varepsilon}$, we have
			\begin{equation} \label{eq:en_low}
				E[\omega] \ge \frac{1}{8\pi} |\log \varepsilon| - C.
			\end{equation}
		\end{lemma}
		
		\begin{proof}
			Since $\omega$ is a maximizer, $E[\omega] \ge E[\rho]$ for any competitor $\rho \in \calF_{\mu, \varepsilon}$.
			We consider $\rho = \rho_+ - \rho_-$, where $\rho_{\pm}$ as uniform disk patches of density $\varepsilon^{-2}$ and mass $1/2$.
			The radius of each disk is $R = \varepsilon / \sqrt{2\pi} < \pi \wedge \frac{\mu}{2}$.
			By placing the centers of $\rho_+$ and $\rho_-$ at $(0, \mu)$ and $(\pi, -\mu)$ respectively, we see that $\|x_2\rho\|_1 = \mu$ which ensures $\rho \in \calF_{\mu, \varepsilon}$.
			By Lemma~\ref{lem_green_decom}, a direct calculation shows that
			\[
			E[\rho_\pm] \ge \frac{1}{16\pi} |\log \varepsilon| - C, \quad |E_{\mathrm{int}}[\rho_+, \rho_-]| \le C(\mu),
			\]
			where, for the interaction energy, we use the fact that $\dist ( \supp \rho_+, \supp \rho_-) \ge 2\mu - 2R \ge \mu$ for sufficiently small $\varepsilon$.
			This proves the lemma.
		\end{proof}
		
		\begin{lemma}\label{lem:patch}
			Every maximizer $\omega \in \calS_{\mu, \varepsilon}$ is of the patch type, that is, any $\omega \in \calS_{\mu, \varep}$ satisfies
			\begin{equation}
				\label{eq:max_patch}
				\omega = \varepsilon^{-2} \left(\mathbf{1}_{\Omega_+} - \mathbf{1}_{\Omega_-} \right)
			\end{equation}
			for some measurable sets $\Omega_{\pm} \subset \calC$ such that $|\Omega_+|= |\Omega_-| \le \varepsilon^2/2$. 
		\end{lemma}
		
		\begin{proof}
			Let $\omega \in \calS_{\mu, \varepsilon}$.
			It suffices to treat the positive part $\omega_+ := \max\{ \omega, 0\}$.
			
			Towards contradiction, suppose that $|\{0 < \omega < \varepsilon^{-2}\}| > 0$; the positive part is not saturated.
			Then there exists $\delta_0 > 0$ such that the set $K = K_{\delta_0} := \{ \delta_0 < \omega < \varep^{-2} - \delta_0\}$ satisfying $|K_{\delta_0}| > 0$.
			We take bounded sets $A, B \subset K$ with positive measure such that
			\[
			\frac{1}{|A|}\int_A |x_2|\,dx \neq \frac{1}{|B|}\int_B |x_2|\,dx.
			\]
			Indeed, if this were not possible, $|x_2|$ would be constant almost everywhere on $K$, which contradicts $|K| > 0$.
			Note that
			\[
			\det \begin{pmatrix}
				|A| & |B| \\
				\int_A |x_2| \,dx & \int_B |x_2| \,dx
			\end{pmatrix} \neq 0.
			\]
			%Thus, by taking appropriate linear combinations $\mathbf{1}_A$ and $\mathbf{1}_B$, we define functions $h_1, h_2 \in L_c^{\infty}(\calC)$ satisfying $\supp h_i \subset K$, $i = 1, 2$, and
			Thus by solving linear systems, we find auxiliary functions $h_i = a_i \mathbf{1}_A + b_i \mathbf{1}_B \in L_c^{\infty}(\calC)$, $i = 1, 2$, satisfying $\supp h_i \subset K$ and
			\begin{equation*}
				\begin{cases}
					\int_{\calC} h_1 \,dx = 1, \quad \int_{\calC} |x_2|h_1 \,dx = 0, \\[4pt]
					\int_{\calC} h_2 \,dx = 0, \quad \int_{\calC} |x_2| h_2 \,dx = 1.
				\end{cases}
			\end{equation*}
			Let $\delta \in (0, \delta_0)$.
			Fix an arbitrary function $h \in  L_c^{\infty}(\calC)$ supported in $\{ \omega \ge 0 \}$ such that
			\begin{equation}
				\label{eq:sign_h}
				\begin{cases}
					h \ge 0 \quad \textrm{ on }\,\,\{ 0 \le \omega \le \delta \},\\
					h \le 0 \quad \textrm{ on }\,\, \{ \varepsilon^{-2} - \delta \le \omega \le \varepsilon^{-2} \}.
				\end{cases}
			\end{equation}
			We set the perturbation as
			\[
			\eta = h - \left( \int_{\calC} h \, dx \right)h_1 - \left( \int_{\calC} |x_2|h \,dx \right)h_2.
			\]
			By construction, it is clear that $\supp \eta \subset \{ \omega \ge 0 \}$ and $\int_{\calC} \eta \,dx = \int_{\calC} |x_2 |\eta \,dx = 0$.\\
			
			\noindent \textbf{Claim:} \textit{$\omega + \sigma \eta \in \calF_{\mu, \varepsilon}$ for sufficiently small $\sigma > 0$.}\\
			
			\noindent To see this, first note that
			\[
			\omega + \sigma \eta = \begin{cases}
				\omega_+ + \sigma \eta  &\textrm{ on }\,\, \{ \omega \ge 0 \},\\
				-\omega_-  &\textrm{ on }\,\, \{\omega < 0\},
			\end{cases}
			\]
			and that
			\[
			\begin{cases}
				h_1 = h_2 = 0 \ge h \quad \textrm{ on }\,\,  \{ \varepsilon^{-2} - \delta \le \omega \le \varepsilon^{-2} \},\\[4pt]
				h_1 = h_2 = 0 \le h\quad \textrm{ on }\,\,  \{ 0 \le \omega \le \delta \}.
			\end{cases}
			\]
			These properties ensure that for $\sigma > 0$ sufficiently small,
			\[
			\| \omega + \sigma \eta\|_{\infty} \le \varep^{-2}, \quad \| \omega + \sigma \eta\|_{1} = \|\omega\|_1 \le 1, \quad \| x_2(\omega + \sigma \eta)\|_{1} = \|x_2\omega\|_1 \le \mu.
			\]
			Since $\int (\omega + \sigma \eta)  = 0$, the claim is verified.\\
			
			Since $\omega$ is a maximizer, we have $E[\omega] \ge E[\omega + \sigma \eta]$ for small $\sigma > 0$ by the claim.
			This yields
			\begin{equation}\label{eq:en_var_1}
				0 \ge \int_{\calC} \psi[\omega](x) \eta(x) \,dx = \int_{\{\omega \ge 0\}} \left( \psi[\omega](x) - W|x_2| - \gamma \right) h(x) \,dx,
			\end{equation}
			where we define
			\[
			W = \int_{\calC} \psi[\omega](x) h_2(x) \,dx, \quad \gamma = \int_{\calC} \psi[\omega](x) h_1(x) \,dx.
			\]
			Since $h$ is arbitrary with the sign property \eqref{eq:sign_h} and $\delta$ is also arbitrary small, from \eqref{eq:en_var_1} we deduce
			\[
			\psi[\omega](x) - W|x_2| - \gamma \,\,
			\begin{cases} \le 0    &\textrm{ on }\,\, \{\omega  = 0\}, \\
				= 0 &\textrm{ on }\,\, \{0 < \omega < \varepsilon^{-2} \}, \\
				\ge 0 &\textrm{ on } \,\, \{ \omega = \varepsilon^{-2}\}.
			\end{cases}
			\]
			In particular, $\psi[\omega](x) = W|x_2| + \gamma $ almost everywhere on the set $\{ 0 < \omega < \varepsilon^{-2} \}$.
			%Since $\omega \in L^\infty(\calC)$, elliptic regularity implies $\psi[\omega] \in W_{loc}^{2,p}(\calC)$ for any $1 < p < \infty$.
			%By Stampacchia's theorem, the Laplacian of $\psi[\omega] - W|x_2| - \gamma$ must vanish almost everywhere on its level set.
			Noting that $\Delta(|x_2|) = 0$ away from the measure-zero set $\{x_2 = 0\}$, we obtain 
			\[
			\omega = -\Delta \psi = 0 \quad \textrm{ a.e. on }\,\, \{ 0 < \omega < \varepsilon^{-2} \},
			\]
			which contradicts $|\{ 0 < \omega < \varepsilon^{-2} \}| > 0$.
			Applying an identical argument to the negative part $\omega_- := \max \{ - \omega, 0 \}$, we obtain \eqref{eq:max_patch} for some measurable sets $\Omega_{\pm} \subset \calC$.
			Finally, the constraints $\int \omega  = 0$ and $\|\omega\|_1 \le 1$ yield $|\Omega_+| = |\Omega_-| \le \varepsilon^{2}/2$.
		\end{proof}
		
		\begin{lemma}\label{lem:no_strip}
			There exists $\varepsilon_0 = \varepsilon_0(\mu) > 0$ such that for any
			$\varepsilon\in(0,\varepsilon_0)$ and any maximizer
			\[
			\omega=\varepsilon^{-2}(\mathbf 1_{\Omega_+}-\mathbf 1_{\Omega_-})\in\calS_{\mu,\varepsilon},
			\]
			neither $\Omega_+$ nor $\Omega_-$ can be of the form $\bbT \times (a,b)$ for some $a<b$, up to a null set.
		\end{lemma}

		\begin{proof}
			We argue by contradiction.
			Assume, for instance, that $\Omega_+=\bbT\times(a,b)$.
			Set $M:=\|\omega_+\|_1=\|\omega_-\|_1 \le 1/2$.
			Then $|\Omega_+|=2\pi(b-a)=M\varepsilon^2$.
			
			We first estimate the self-energy of the strip patch $\omega_+$.
			By Lemma~\ref{lem_green_decom},
			\[
			E[\omega_+] = \frac{1}{2} \iint H(x-y)\omega_+(x)\omega_+(y)\,dxdy -  \frac{1}{8\pi}\iint |x_2-y_2|\omega_+(x)\omega_+(y)\,dxdy.
			\]
			Since $\omega_+\ge0$, the last term is non-positive.
			Therefore, using \eqref{green_sing_3},
			\[
			E[\omega_+] \le C M^2 + C\varepsilon^{-4} \iint_{\Omega_+\times\Omega_+} \log_+\frac1{|x-y|}\,dxdy.
			\]
			For each fixed $x \in \Omega_+$,
			\[
			\int_{\Omega_+}\log_+\frac1{|x-y|}\,dy \le \int_{|t|\le b-a}\int_{-\pi}^{\pi} \log_+\frac1{\sqrt{s^2+t^2}}\,dsdt \le C(b-a).
			\]
			Then, since $|\Omega_+|=2\pi(b-a) = M\varepsilon^2$ ,
			\[
			\iint_{\Omega_+\times\Omega_+} \log_+\frac1{|x-y|}\,dxdy \lesssim \varepsilon^4,
			\]
			and it follows that $E[\omega_+]\le C M^2\le C$.
			
			On the other hand, by Lemma~\ref{lem:en_pos} applied to $\omega_-$,
			\[
			E[\omega_-]
			\le
			\frac{1}{8\pi}M^2
			\log_+\left(\frac{\varepsilon^{-2}}{M}\right)
			+
			C M^2
			\le
			\frac{M^2}{4\pi}|\log\varepsilon|+C
			\le
			\frac{1}{16\pi}|\log\varepsilon|+C.
			\]
			For the interaction term, using Lemma~\ref{lem_green_decom}, we have $-G(x-y) \lesssim 1 + |x_2| + |y_2|$.
			Therefore,
			\[
			-E_{\mathrm{int}}[\omega_+,\omega_-] \lesssim \iint (1+|x_2|+|y_2|) \omega_+(x)\omega_-(y)\,dxdy \le C(\mu).
			\]
			
			Combining the above, we obtain
			\[
			E[\omega] \le \frac{1}{16\pi}|\log\varepsilon|+C(\mu).
			\]
			This contradicts Lemma~\ref{lem:en_low}, which gives $E[\omega]\ge \frac{1}{8\pi} |\log\varepsilon| - C$.
			The case where $\Omega_-$ is a strip is identical.
			This completes the proof.
		\end{proof}
		
		\begin{lemma}[Auxiliary functions]\label{lem:aux_basis}
			Let $\Omega\subset\calC$ be a measurable set.
			Assume that there exist two points $y, z \in \calC$ with $|z_2| > |y_2|$ such that every neighborhood of each of $y$ and $z$ intersects both $\Omega$ and $\Omega^c$ in sets of positive measure.
			Then, there exist a sequence $\rho_n \to 0$ and corresponding functions $h_{n,1}^{\pm},  h_{n,2}^{\pm}\in L_c^\infty(\calC)$ such that
			\[
			\supp h_{n, i}^{\pm} \subset B_{\rho_n}(y)\cup B_{\rho_n}(z), \,\, i = 1, 2,
			\]
			satisfying
			\begin{equation}\label{eq:aux_basis_c1}
				\int_{\calC} h_{n,1}^{\pm}\,dx=1,
				\quad
				\int_{\calC}|x_2|h_{n,1}^{\pm}\,dx=0,
			\end{equation}
			and
			\begin{equation}\label{eq:aux_basis_c2}
				\int_{\calC} h_{n,2}^{\pm}\,dx = 0,
				\quad
				\int_{\calC}|x_2|h_{n,2}^{\pm}\,dx = 1.
			\end{equation}
			Moreover, for all sufficiently large $n$, these functions satisfy
			\begin{equation}\label{eq:aux_basis_sign}
				\begin{cases}
					h_{n,i}^{+}\ge0 \text{ on } \Omega,\quad h_{n,i}^{+}\le0 \text{ on } \Omega^c,\\[4pt]
					h_{n,i}^{-}\le0 \text{ on } \Omega,\quad h_{n,i}^{-}\ge0 \text{ on } \Omega^c,
				\end{cases}
				\quad i=1,2,
			\end{equation}
			and, for any continuous function $\phi$,
			\begin{equation}\label{eq:aux_basis_limit_1}
				\lim_{n\to\infty}\int_{\calC}\phi h_{n,1}^{+}\,dx
				=
				\lim_{n\to\infty}\int_{\calC}\phi h_{n,1}^{-}\,dx
				=
				\frac{|z_2|\phi(y)-|y_2|\phi(z)}{|z_2|-|y_2|},
			\end{equation}
			\begin{equation}\label{eq:aux_basis_limit_2}
				\lim_{n\to\infty}\int_{\calC}\phi h_{n,2}^{+}\,dx
				=
				\lim_{n\to\infty}\int_{\calC}\phi h_{n,2}^{-}\,dx
				=
				\frac{\phi(z)-\phi(y)}{|z_2|-|y_2|}.
			\end{equation}
		\end{lemma}
		
		\begin{proof}
			We follow the argument in \cite[Proposition 5.3]{Choi2024}.
			Let $\rho_n \searrow 0$.
			By the assumption on $y$ and $z$, after passing to a subsequence if necessary, we may choose measurable sets
			\[
			Y_n^+\subset \Omega\cap B_{\rho_n}(y),\quad Y_n^-\subset \Omega^c\cap B_{\rho_n}(y),
			\]
			and
			\[
			Z_n^+\subset \Omega\cap B_{\rho_n}(z),\quad Z_n^-\subset \Omega^c\cap B_{\rho_n}(z),
			\]
			all of positive measure.
			Define
			\[
			f_n^\pm := \frac{1}{|Y_n^\pm|} \mathbf 1_{Y_n^\pm}, \quad g_n^\pm := \frac{1}{|Z_n^\pm|}\mathbf 1_{Z_n^\pm}.
			\]
			Then
			\[
			\int_{\calC}|x_2|f_n^\pm\,dx \to |y_2|, \quad \int_{\calC}|x_2|g_n^\pm\,dx \to |z_2|.
			\]
			Set
			\[
			\alpha_n^\pm := \int_{\calC}|x_2|f_n^\pm\,dx, \quad \beta_n^\pm := \int_{\calC}|x_2|g_n^\pm\,dx.
			\]
			%Since $|z_2|>|y_2|$, for all sufficiently large $n$ all denominators below are positive.
			Define
			\[
			h_{n,1}^+ := \frac{\beta_n^- f_n^+ - \alpha_n^+g_n^-}{\beta_n^- - \alpha_n^+}, \quad h_{n,2}^+ :=  \frac{ - f_n^- + g_n^+}{\beta_n^+ - \alpha_n^-},
			\qquad
			\mbox{and}
			\qquad 
			h_{n,1}^- :=  \frac{\beta_n^+f_n^- -\alpha_n^- g_n^+}{\beta_n^+ - \alpha_n^-}, \quad h_{n,2}^- := \frac{-f_n^+ + g_n^-}{\beta_n^- - \alpha_n^+},
			\]
			where the denominators satisfy $\beta_n^- - \alpha_n^+ > 0$ and $\beta_n^+ - \alpha_n^- > 0$ for all sufficiently large $n$, since $\alpha_n^\pm\to |y_2|$, $\beta_n^\pm\to |z_2|$ and $|z_2|>|y_2|$.
			With this choice, a direct computation proves \eqref{eq:aux_basis_c1} and \eqref{eq:aux_basis_c2}.
			Here, for example, $h_{n, 1}^+$  is determined by requiring
			\[
			\int_{\calC} h_{n,1}^+\,dx = 1, \quad \int_{\calC}|x_2|h_{n,1}^+\,dx = 0,
			\]
			equivalently, by solving the linear system $\mathbf{a} - \mathbf{b} = 1$ and $\mathbf{a} \alpha_n^+ - \mathbf{b} \beta_n^- = 0$ for $( \mathbf{a}, \mathbf{b}) \in \bbR^2$.
			The other functions are obtained in the same way. Since $\beta_n^- - \alpha_n^+, \beta_n^+ - \alpha_n^- > 0$ and  $\supp f_n^+, \supp g_n^+ \subset \Omega$, $\supp f_n^-, \supp g_n^- \subset \Omega^c$, the sign property \eqref{eq:aux_basis_sign} follows.
			
			It remains to verify the limits.
			Since $\phi$ is continuous and the supports of $f_n^\pm$ and $g_n^\pm$ shrink to $y$ and $z$, respectively, we have $\int \phi f_n^\pm \to \phi(y)$ and  $\int \phi g_n^\pm \to \phi(z)$.
			Together with the explicit choice above, we find \eqref{eq:aux_basis_limit_1} and \eqref{eq:aux_basis_limit_2}.
		\end{proof}

		\begin{lemma}\label{lem:local_EL}
			Let $\varepsilon_0=\varepsilon_0(\mu)>0$ be as in Lemma~\ref{lem:no_strip}.
			For any $\varepsilon \in (0,\varepsilon_0)$ and any $\omega\in\calS_{\mu,\varepsilon}$, there exist unique constants $W_{\pm},\gamma_{\pm}\in\bbR$ such that the sets $\Omega_{\pm}$ in \eqref{eq:max_patch} satisfy
			\begin{equation}\label{eq:local_EL_patches}
				\begin{aligned}
					\Omega_{\pm} = \{ x \in \calC : \pm\psi[\omega](x) - W_{\pm}|x_2| - \gamma_{\pm} > 0 \}.
				\end{aligned}
			\end{equation}
		\end{lemma}
		
		\begin{proof}
			\textbf{Step 1: Auxiliary functions for $\omega_+$}.
			By the Lebesgue density theorem and Lemma~\ref{lem:no_strip}, we choose two points $y, z \in \calC$ with $|z_2| > |y_2|$ such that every neighborhood of each of $y$ and $z$ intersects both $\Omega_+$ and $\Omega_+^c$ with positive measure.
			
			Applying Lemma~\ref{lem:aux_basis} with $(\Omega, \phi) = (\Omega_+, \psi[\omega])$, we obtain functions $h_{n,1}^{\pm},\ h_{n,2}^{\pm}\in L_c^\infty(\calC)$ such that
			\[
			\int_{\calC} h_{n,1}^{\pm}\,dx=1,
			\quad
			\int_{\calC}|x_2|h_{n,1}^{\pm}\,dx=0, \qquad\mbox{and}
			\qquad 
			\int_{\calC} h_{n,2}^{\pm}\,dx=0,
			\quad
			\int_{\calC}|x_2|h_{n,2}^{\pm}\,dx=1.
			\]
			Moreover,
			\begin{equation}\label{eq:bas_sign_pos}
				\begin{cases}
					h_{n,i}^+ \ge 0 \text{ on } \Omega_+, \quad h_{n,i}^+ \le 0 \text{ on } \Omega_+^c, \\[4pt]
					h_{n,i}^- \le 0 \text{ on } \Omega_+, \quad h_{n,i}^- \ge 0 \text{ on } \Omega_+^c,
				\end{cases}
				\quad i=1,2,
			\end{equation}
			and
			\begin{equation}\label{eq:aux_lim_positive}
				\begin{aligned}
					\lim_{n\to\infty}\int_{\calC}\psi[\omega] h_{n,1}^+\,dx
					&=
					\lim_{n\to\infty}\int_{\calC}\psi[\omega] h_{n,1}^-\,dx
					=:\gamma_+,\\
					\lim_{n\to\infty}\int_{\calC}\psi[\omega] h_{n,2}^+\,dx
					&=
					\lim_{n\to\infty}\int_{\calC}\psi[\omega] h_{n,2}^-\,dx
					=:W_+.
				\end{aligned}
			\end{equation}
			
			\medskip
			
			\noindent \textbf{Step 2: Admissible perturbation on $\Omega_+$}.
			Let $h \in L_c^\infty(\calC)$ be any function satisfying
			\begin{equation}\label{eq:h_sign_positive}
				h \le 0 \quad \text{on } \Omega_+, \quad h \ge 0 \quad \text{on } \Omega_+^c.
			\end{equation}
			To ensure the perturbation preserves this sign property, we select the basis functions by the following:
			\[
			h_{n,1} := \begin{cases} 
				h_{n,1}^+ & \text{if } \int_{\calC} h \,dx \ge 0, \\ 
				h_{n,1}^- & \text{if } \int_{\calC} h \,dx < 0, 
			\end{cases} \quad
			h_{n,2} := \begin{cases} 
				h_{n,2}^+ & \text{if } \int_{\calC} |x_2|h \,dx \ge 0, \\
				h_{n,2}^- & \text{if } \int_{\calC} |x_2|h \,dx < 0. 
			\end{cases}
			\]
			We then define the perturbation as
			\[
			\eta_n := h - \left(\int_{\calC} h \,dx\right) h_{n,1} - \left(\int_{\calC} |x_2|h \,dx\right) h_{n,2}.
			\]
			With this choice, we have $\eta_n \le 0$ on $\Omega_+$ and $\eta_n \ge 0$ on $\Omega_+^c$.
			By construction, $\int_{\calC} \eta_n \,dx = 0$ and $\int_{\calC} |x_2| \eta_n \,dx = 0$.
			
			We fix $n$ and claim that, for any sufficiently small $\sigma > 0$, $\omega + \sigma \eta_n \in \calF_{\mu, \varepsilon}$.
			The $L^\infty$ bound and $\int  (\omega + \sigma\eta_n)  = 0$ are immediate.
			For the $L^1$ constraint, from Lemma~\ref{lem:patch} we see that
			\begin{align*}
				\|\omega + \sigma \eta_n\|_1 &= \int_{\Omega_+} (\varepsilon^{-2} + \sigma \eta_n) \,dx + \int_{\Omega_-} (\varepsilon^{-2} - \sigma \eta_n) \,dx + \int_{ (\Omega_+ \cup \Omega_-)^c} \sigma \eta_n \,dx \\
				&= \|\omega\|_1 + \sigma \left( \int_{\Omega_+} \eta_n \,dx - \int_{\Omega_-} \eta_n \,dx + \int_{(\Omega_+ \cup \Omega_-)^c} \eta_n \,dx \right).
			\end{align*}
			Using $\int \eta_n  = 0$, we find $\|\omega + \sigma \eta_n\|_1 = \|\omega\|_1 - 2\sigma \int_{\Omega_-} \eta_n  \le \|\omega\|_1 \le 1$ since $\eta_n \ge 0$ on $\Omega_- \subset \Omega_+^c$.
			Similarly, $\|x_2(\omega + \sigma \eta_n)\|_1 = \|x_2\omega\|_1 - 2\sigma \int_{\Omega_-} |x_2|\eta_n  \le \mu$.
			Thus $\omega + \sigma \eta_n \in \calF_{\mu, \varep}$.
			
			\medskip
			
			\noindent \textbf{Step 3: Variational inequality}.
			Since $\omega$ is a maximizer, we have $\frac{d}{d\sigma} E[\omega + \sigma \eta_n] \big|_{\sigma \to  0+} \le 0$ and then $\int \psi[\omega] \eta_n  \le 0$.
			By taking $n \to \infty$, we obtain
			\[
			\int_{\calC} ( \psi[\omega] - W_+|x_2| - \gamma_+ ) h \,dx \le 0,
			\]
			where $W_+$ and $\gamma_+$ are defined in \eqref{eq:aux_lim_positive}.
			Note that $W_+$ and $\gamma_+$ are independent of the choice of $h$.
			Since $h$ is arbitrary under the sign condition \eqref{eq:h_sign_positive}, it follows that
			\[
			\begin{cases}
				\psi[\omega](x) - W_+|x_2| - \gamma_+ \ge 0 \quad \text{a.e. on } \Omega_+, \\
				\psi[\omega](x) - W_+|x_2| - \gamma_+ \le 0 \quad \text{a.e. on } \Omega_+^c.
			\end{cases}
			\]
			Let $\Psi_+(x) := \psi[\omega](x) - W_+|x_2| - \gamma_+$.
			The above inequalities imply $\{\Psi_+ > 0\} \subset \Omega_+ \subset \{\Psi_+ \ge 0\}$ modulo null sets. 
			Suppose for contradiction that the level set $E := \Omega_+ \cap \{\Psi_+ = 0\}$ has strictly positive measure.
			By elliptic regularity, $\Psi_+ \in W^{2,p}(K)$ for any $p \in (1, \infty)$ and compact subset $K \subset \calC \setminus \{x_2 = 0\}$, and $\Delta \Psi_+ = 0$ a.e. on $E$. 
			Thus, we see
			\[
			-\Delta \Psi_+(x) = -\Delta \psi[\omega](x) + W_+ \Delta(|x_2|) = \omega(x) = \varepsilon^{-2} > 0 \quad \text{a.e. on } E \setminus \{x_2 = 0\}.
			\]
			Since $|\{x_2 = 0\}| = 0$, this implies $-\Delta \Psi_+ > 0$ a.e. on $E$, which directly contradicts $\Delta \Psi_+ = 0$ a.e. on $E$.
			Thus, $|E| = 0$ and we conclude that, modulo null sets,
			\[
			\Omega_+ = \{ \psi[\omega](x) - W_+|x_2| - \gamma_+ > 0 \}.
			\]
			
			\medskip
			
			\noindent \textbf{Step 4: Negative patch}.
			We apply the same argument with $(\Omega,\phi)=(\Omega_-,-\psi[\omega])$ in Lemma~\ref{lem:aux_basis}.
			More precisely, Lemma~\ref{lem:aux_basis} gives auxiliary functions $k_{n,1}^{\pm},\ k_{n,2}^{\pm}\in L_c^\infty(\calC)$ such that
			\[
			\int_{\calC} k_{n,1}^{\pm}\,dx = 1, \quad \int_{\calC}|x_2|k_{n,1}^{\pm}\,dx =  0,
			\qquad
			\mbox{and}
			\qquad 
			\int_{\calC} k_{n,2}^{\pm}\,dx =  0, \quad \int_{\calC}|x_2|k_{n,2}^{\pm}\,dx =  1.
			\]
			Moreover,
			\begin{equation}\label{eq:basis_sign_negative}
				\begin{cases}
					k_{n,i}^+ \ge 0 \text{ on } \Omega_-, \quad k_{n,i}^+ \le 0 \text{ on } \Omega_-^c, \\
					k_{n,i}^- \le 0 \text{ on } \Omega_-, \quad k_{n,i}^- \ge 0 \text{ on } \Omega_-^c,
				\end{cases}
				\quad i=1,2,
			\end{equation}
			and
			\begin{equation}\label{eq:aux_lim_negative}
				\begin{aligned}
					\lim_{n\to\infty}\int_{\calC}(-\psi[\omega]) k_{n,1}^+\,dx
					&=
					\lim_{n\to\infty}\int_{\calC}(-\psi[\omega]) k_{n,1}^-\,dx
					=:\gamma_-,\\
					\lim_{n\to\infty}\int_{\calC}(-\psi[\omega]) k_{n,2}^+\,dx
					&=
					\lim_{n\to\infty}\int_{\calC}(-\psi[\omega]) k_{n,2}^-\,dx
					=:W_-.
				\end{aligned}
			\end{equation}
			Repeating \textbf{Step 2} and \textbf{Step 3} with the sign condition
			\[
			h \ge 0 \,\, \text{ on }\,\, \Omega_-, \quad h \le 0  \,\, \text{ on }\,\,\Omega_-^c,
			\]
			we deduce
			\[
			\Omega_- = \{-\psi[\omega](x)-W_-|x_2|-\gamma_->0\}.
			\]
			
			\medskip
			
			\noindent \textbf{Step 5: Uniqueness of multipliers.}
			Finally, we show that the constants $W_{\pm}$ and $\gamma_{\pm}$ are uniquely determined by $\omega$.
			We first consider the positive patch.
			After redefining $\Omega_+$ on a null set, we may assume that
			\[
			\Omega_+ = \{x \in \calC : \psi[\omega](x)-W_+|x_2|-\gamma_+ > 0\}.
			\]
			Since $\psi[\omega]$ is continuous, every point $x\in\partial\Omega_+$ satisfies $\psi[\omega](x)-W_+|x_2|-\gamma_+=0$.
			By Lemma~\ref{lem:no_strip}, we can choose two boundary points $y,z\in\partial\Omega_+$ such that $|y_2|\neq |z_2|$.
			Hence $W_+$ and $\gamma_+$ are uniquely determined by the linear system
			\[
			\begin{pmatrix}
				|y_2| & 1 \\
				|z_2|  & 1
			\end{pmatrix}
			\begin{pmatrix}
				W_+ \\
				\gamma_+
			\end{pmatrix} = \begin{pmatrix}
				\psi[\omega](y) \\
				\psi[\omega](z)
			\end{pmatrix}.
			\]
			The same argument shows that $W_-$ and $\gamma_-$ are uniquely determined as well.
		\end{proof}
		
		\begin{lemma}\label{lem:impulse_multipl}
			Let $\varepsilon_0=\varepsilon_0(\mu)>0$ be as in Lemma~\ref{lem:no_strip}.
			For any $\varepsilon\in(0,\varepsilon_0)$ and $\omega\in\calS_{\mu,\varepsilon}$, the multipliers $W_{\pm}$ obtained in Lemma~\ref{lem:local_EL} satisfy $W_+ = W_- =: W \ge 0$.
		\end{lemma}
		
		\begin{proof}
			We first prove that $W_+=W_-$.
			
			\medskip
			
			\noindent \textbf{Auxiliary setup}.
			Let $h_{n,2}^{\pm}$ be the sequence of auxiliary functions for $\Omega_+$ in \textbf{Step 1} (the positive patch part) of the proof of Lemma~\ref{lem:local_EL}.
			Define $\eta_{\mathrm{up},n}^+ := h_{n,2}^-, \quad \eta_{\mathrm{down},n}^+ := -h_{n,2}^+.$ Then both satisfy
			\begin{equation}
				\label{eq:cross_sign}
				\eta_{\mathrm{up},n}^+,\,\, \eta_{\mathrm{down},n}^+ \le 0 \,\, \text{ on }\,\,\Omega_+,
				\quad
				\eta_{\mathrm{up},n}^+,\,\, \eta_{\mathrm{down},n}^+ \ge 0 \,\,  \text{ on } \,\, \Omega_+^c,
			\end{equation}
			and
			\[
			\int_{\calC}\eta_{\mathrm{up},n}^+\,dx = \int_{\calC}\eta_{\mathrm{down},n}^+\,dx = 0, \quad \int_{\calC}|x_2|\eta_{\mathrm{up},n}^+\,dx =  -\int_{\calC}|x_2|\eta_{\mathrm{down},n}^+\,dx = 1,
			\]
			that is, `up' for the vertical moment $1$ and `down' for $-1$.
			Furthermore,
			\[
			\lim_{n\to\infty}\int_{\calC}\psi[\omega]\eta_{\mathrm{up},n}^+\,dx=W_+, \quad
			\lim_{n\to\infty}\int_{\calC}\psi[\omega]\eta_{\mathrm{down},n}^+\,dx=-W_+.
			\]
			We also use the auxiliary functions $k_{n, 2}^{\pm}$ in \textbf{Step 4} (the negative patch part) of the proof of Lemma~\ref{lem:local_EL}.
			Define $\eta_{\mathrm{up},n}^- := k_{n,2}^+, \quad \eta_{\mathrm{down},n}^- := -k_{n,2}^-.$ 
			Then both satisfy
			\[
			\eta_{\mathrm{up},n}^-, \eta_{\mathrm{down},n}^- \ge 0 \,\, \text{ on } \,\, \Omega_-,
			\quad
			\eta_{\mathrm{up},n}^-, \eta_{\mathrm{down},n}^- \le 0 \,\, \text{ on } \,\, \Omega_-^c,
			\]
			and
			\[
			\int_{\calC}\eta_{\mathrm{up},n}^-\,dx = \int_{\calC}\eta_{\mathrm{down},n}^-\,dx = 0, \quad \int_{\calC}|x_2|\eta_{\mathrm{up},n}^-\,dx = - \int_{\calC}|x_2|\eta_{\mathrm{down},n}^-\,dx = 1.
			\]
			Also,
			\[
			\lim_{n\to\infty}\int_{\calC}\psi[\omega]\eta_{\mathrm{up},n}^-\,dx = -W_-,
			\quad
			\lim_{n\to\infty}\int_{\calC}\psi[\omega]\eta_{\mathrm{down},n}^-\,dx = W_-.
			\]
			
			\medskip
			
			\noindent \textbf{Cross perturbations}.
			We first show $W_+ = W_-$.
			Set
			\[
			\eta_{A,n}:=\eta_{\mathrm{up},n}^+ + \eta_{\mathrm{up},n}^-.
			\]
			Then $\int \eta_{A,n}  = \int \omega + \sigma\eta_{A,n}  = 0$.
			We claim that $\omega+\sigma\eta_{A,n}\in\calF_{\mu,\varepsilon}$ for sufficiently small $\sigma>0$.
			Indeed, both $\eta_{\mathrm{up},n}^+$ and $\eta_{\mathrm{up},n}^-$ are nonpositive on $\Omega_+$, and nonnegative on $\Omega_-$.
			Hence, for sufficiently small $\sigma>0$, we have $\| \omega + \sigma \eta_{A, n}\|_{\infty} \le \varepsilon^{-2}$.
			Also, on $(\Omega_+\cup\Omega_-)^c$, we have $\eta_{\mathrm{up},n}^+ \ge 0$, $\eta_{\mathrm{up},n}^- \le 0$ and then,
			\[
			|\eta_{\mathrm{up},n}^+ + \eta_{\mathrm{up},n}^-| \le \eta_{\mathrm{up},n}^+ - \eta_{\mathrm{up},n}^- \quad \textrm{ on } \,\, (\Omega_+\cup\Omega_-)^c.
			\]
			Thus
			\begin{align*}
				\|\omega+\sigma\eta_{A,n}\|_1  & = \|\omega\|_1 + \sigma\bigg[
				\int_{\Omega_+}(\eta_{\mathrm{up},n}^+ + \eta_{\mathrm{up},n}^-)\,dx
				-\int_{\Omega_-}(\eta_{\mathrm{up},n}^+ + \eta_{\mathrm{up},n}^-)\,dx +\int_{(\Omega_+\cup\Omega_-)^c}
				|\eta_{\mathrm{up},n}^+ + \eta_{\mathrm{up},n}^-|\,dx
				\bigg] \\
				& \le \|\omega\|_1 +
				\sigma\bigg[ \int_{\Omega_+}(\eta_{\mathrm{up},n}^+ + \eta_{\mathrm{up},n}^-)\,dx
				-\int_{\Omega_-}(\eta_{\mathrm{up},n}^+ + \eta_{\mathrm{up},n}^-)\,dx 
				+\int_{(\Omega_+\cup\Omega_-)^c}
				(\eta_{\mathrm{up},n}^+ - \eta_{\mathrm{up},n}^-)\,dx
				\bigg] \\
				& =  \|\omega\|_1 + 2\sigma\left[ -\int_{\Omega_-}\eta_{\mathrm{up},n}^+\,dx + \int_{\Omega_+}\eta_{\mathrm{up},n}^-\,dx \right]
				\le \|\omega\|_1 \le 1,
			\end{align*}
			where we used $\int \eta_{\mathrm{up},n}^{\pm}  = 0$ with the sign property \eqref{eq:cross_sign}.
			Similarly,
			\begin{align*}
				\|x_2(\omega+\sigma\eta_{A,n})\|_1
				&\le
				\|x_2\omega\|_1 + \sigma\bigg[ \int_{\Omega_+}|x_2|(\eta_{\mathrm{up},n}^+ + \eta_{\mathrm{up},n}^-)\,dx -\int_{\Omega_-}|x_2|(\eta_{\mathrm{up},n}^+ + \eta_{\mathrm{up},n}^-)\,dx \\ &\hspace{4cm} +\int_{(\Omega_+\cup\Omega_-)^c} |x_2|(\eta_{\mathrm{up},n}^+ - \eta_{\mathrm{up},n}^-)\,dx \bigg] \\
				&=
				\|x_2\omega\|_1 + 2\sigma\left[-\int_{\Omega_-}|x_2| \eta_{\mathrm{up},n}^+\,dx + \int_{\Omega_+}|x_2| \eta_{\mathrm{up},n}^- \,dx \right]
				\le \mu.
			\end{align*}
			Here we used $\int |x_2|\eta_{\mathrm{up},n}^{\pm} = 1$, again with \eqref{eq:cross_sign}.
			Therefore $\omega+\sigma\eta_{A,n}\in\calF_{\mu,\varepsilon}$.
			By the maximality of $\omega$,
			\[
			0\ge
			\lim_{n\to\infty}\int_{\calC}\psi[\omega]\eta_{A,n}\,dx
			=
			W_+-W_-.
			\]
			Hence $W_+\le W_-$.
			Conversely, set
			\[
			\eta_{B,n}:=\eta_{\mathrm{down},n}^+ + \eta_{\mathrm{down},n}^-.
			\]
			The same argument gives $\omega+\sigma\eta_{B,n}\in\calF_{\mu,\varepsilon}$ for sufficiently small $\sigma>0$, and then,
			\[
			0 \ge \lim_{n\to\infty}\int_{\calC}\psi[\omega]\eta_{B,n}\,dx = -W_+ + W_-, \quad \textrm{ i.e.} \,\, W_-\le W_+.
			\]
			Therefore, $W_+=W_-=:W$.
			
			\medskip 
			
			\noindent \textbf{Non-negativity: $W \ge 0$}.
			Suppose for contradiction that $W<0$.
			By Lemma~\ref{lem:local_EL},
			\[
			\Omega_+ = \{x\in\calC : \psi[\omega](x) - W|x_2| - \gamma_+ > 0 \}.
			\]
			Since $\psi[\omega]\in L^\infty(\calC)$ (Lemma~\ref{lem_str_bdd}) and $-W |x_2| \to + \infty$ as $|x_2| \to \infty$, there exists $R > 0$ such that $\psi[\omega](x) - W|x_2| - \gamma_+ > 0$ for $|x_2| > R$.
			This contradicts $|\Omega_+| \le \varepsilon^{2}/2$.
			Therefore $W \ge 0$.
		\end{proof}
		
		Finally, we establish $\gamma_+ + \gamma_- \ge 0$ by a similar argument to Lemma~\ref{lem:impulse_multipl}.
		
		\begin{lemma}\label{lem:mass_multip}
			Let  $\varepsilon_0 = \varepsilon_0(\mu) > 0$ be as in Lemma~\ref{lem:no_strip}.
			For any $\varepsilon\in(0,\varepsilon_0)$ and any $\omega\in\calS_{\mu,\varepsilon}$, the unique multipliers $\gamma_{\pm}$ obtained in Lemma~\ref{lem:local_EL} satisfy $\gamma_+ + \gamma_- \ge 0$.
		\end{lemma}
		
		\begin{proof}
			Let $h_{n,1}^+$ and $k_{n,1}^+$ be the functions constructed in the proof of Lemma~\ref{lem:local_EL}.
			They satisfy $h_{n,1}^+ \ge 0$ on $\Omega_+$, $h_{n,1}^+ \le 0$ on $\Omega_+^c$, and $k_{n,1}^+ \ge 0$ on $\Omega_-$, $k_{n,1}^+ \le 0$ on $\Omega_-^c$.
			Their integrals are
			\begin{align*}
				\int_{\calC} h_{n,1}^+ \,dx = 1, \quad \int_{\calC} |x_2|h_{n,1}^+ \,dx = 0,  \quad \textrm{and} \quad \int_{\calC} k_{n,1}^+ \,dx = 1, \quad \int_{\calC} |x_2|k_{n,1}^+ \,dx = 0.
			\end{align*}
			We also have
			\[
			\lim_{n \to \infty} \int_{\calC} \psi[\omega] h_{n,1}^+ \,dx = \gamma_+, \quad \lim_{n \to \infty} \int_{\calC} \psi[\omega] k_{n,1}^+ \,dx = -\gamma_-.
			\]
			We define
			\[
			\eta_{C, n} := -h_{n,1}^+ + k_{n,1}^+.
			\]
			By construction, $\int \eta_{C, n}  = 0$ and $\int |x_2|\eta_{C, n}  = 0$.
			Since $h_{n,1}^+ \le 0$ and $k_{n,1}^+ \le 0$ on $(\Omega_+ \cup \Omega_-)^c$, we have $|\eta_{C, n}| \le -h_{n,1}^+ - k_{n,1}^+$ on $(\Omega_+ \cup \Omega_-)^c$. 
			For small $\sigma > 0$, we compute
			\begin{align*}
				\|\omega + \sigma \eta_{C, n}\|_1 &\le \|\omega\|_1 + \sigma \left[ \int_{\Omega_+} (-h_{n,1}^+ + k_{n,1}^+) - \int_{\Omega_-} (-h_{n,1}^+ + k_{n,1}^+) + \int_{(\Omega_+ \cup \Omega_-)^c} (-h_{n,1}^+ - k_{n,1}^+) \right] \\
				&= \|\omega\|_1 + \sigma \left[ \left( \int_{\calC} (-h_{n,1}^+) + 2\int_{\Omega_-} h_{n,1}^+ \right) + \left( -\int_{\calC} k_{n,1}^+ + 2\int_{\Omega_+} k_{n,1}^+ \right) \right].
			\end{align*}
			Since $\int h_{n,1}^+  = 1$ with $h_{n,1}^+ \le 0$ on $\Omega_- \subset \Omega_+^c$, and since $\int k_{n,1}^+  = 1$ with $k_{n,1}^+ \le 0$ on $\Omega_+ \subset \Omega_-^c$, we obtain
			\[
			\|\omega + \sigma \eta_{C, n}\|_1 \le \|\omega\|_1 - 2\sigma \le 1.
			\]
			Similarly, $\|x_2(\omega + \sigma\eta_{C, n})\|_1$ is bounded above by
			\begin{align*}
				\|x_2\omega\|_1 
				+ \sigma \left[ \left( \int_{\calC} |x_2|(-h_{n,1}^+) + 2\int_{\Omega_-} |x_2|h_{n,1}^+ \right) + \left( -\int_{\calC} |x_2|k_{n,1}^+ + 2\int_{\Omega_+} |x_2|k_{n,1}^+ \right) \right].
			\end{align*}
			Since $\int |x_2|h_{n,1}^+  = 0$ and $\int |x_2|k_{n,1}^+  = 0$, we have $\|x_2(\omega + \sigma\eta_{C, n})\|_1 \le \|x_2\omega\|_1 \le \mu$.
			Therefore, $\omega + \sigma\eta_{C, n} \in \calF_{\mu, \varepsilon}$.
			From the maximality of $\omega$,
			\[
			0 \ge \lim_{n \to \infty} \int_{\calC} \psi[\omega] \eta_{C, n} \,dx = \lim_{n \to \infty} \int_{\calC} \psi[\omega] (-h_{n,1}^+) \,dx + \lim_{n \to \infty} \int_{\calC} \psi[\omega] k_{n,1}^+ \,dx = -\gamma_+ - \gamma_-.
			\]
			This implies $\gamma_+ + \gamma_- \ge 0$.
		\end{proof}

		\begin{proof}[\textbf{Proof of Proposition \ref{prop:EL}}]
			This follows from combining Lemmas~\ref{lem:patch}, \ref{lem:local_EL}, \ref{lem:impulse_multipl}, and \ref{lem:mass_multip}.
		\end{proof}

		\subsection{Pohozaev identity}\label{subsec:Po}
		\begin{lemma}\label{lem:pohozaev_identity}
			Let $\varepsilon_0 = \varepsilon_0(\mu) > 0$ be as in Proposition~\ref{prop:EL}.
			Then, for any $\varepsilon \in (0,\varepsilon_0)$ and any $\omega\in\calS_{\mu,\varepsilon}$ satisfying \eqref{eq:EL}, we have
			\begin{equation}\label{eq:poho}
				2W\mu + \frac{(\gamma_+ + \gamma_-)\|\omega\|_1}{2}= E[\omega] + \|\partial_2\psi\|_2^2.
			\end{equation}
		\end{lemma}
		
		\begin{proof}
			Recall that $\|x_2\omega\|_1=\mu$ and that $M:=\|\omega_+\|_1=\|\omega_-\|_1=\|\omega\|_1/2\le 1/2$.
			By Lemmas~\ref{lem_str_bdd} and~\ref{lem_L2}, $\psi \in L^\infty(\calC)$ and $\nabla \psi \in L^2(\calC)$.
			For the given $\omega \in \calS_{\mu, \varepsilon}$ satisfying \eqref{eq:EL}, $W$ and $\gamma_\pm$ are finite constants.
			Let us define
			\begin{equation*}
				\Psi_+(x) := \max(\psi - W|x_2| - \gamma_+, 0), \quad \Psi_-(x) := \max(-\psi - W|x_2| - \gamma_-, 0).
			\end{equation*}
			Notice that if $\Psi_+(x) > 0$ and $\Psi_-(x) > 0$ for the same $x \in \calC$, it would imply $-2W|x_2| - (\gamma_+ + \gamma_-) > 0$, contradicting $W \ge 0$ and $\gamma_+ + \gamma_- \ge 0$.
			Thus, $\varepsilon^{-2}(\Psi_+ + \Psi_-) = \Psi_+ \omega_+ + \Psi_- \omega_-$.
			
			Since $\omega \in L^1(\calC)$, $\psi \in L^\infty(\calC)$, and $x_2\omega \in L^1(\calC)$, we have $\Psi_\pm \omega_\pm \in L^1(\calC)$.
			Therefore, the function
			\begin{equation*}
				K(x_2) := \int_{\mathbb{T}} |\nabla \psi(x)|^2 + \Psi_+(x)\omega_+(x) + \Psi_-(x)\omega_-(x) \, dx_1
			\end{equation*}
			belongs to $L^1(\bbR)$.
			%By considering $\tilde{K}(r) := K(r) + K(-r) \in L^1(\bbR_+)$, the integrability implies $\liminf_{r \to \infty} r \tilde{K}(r) = 0$. 
			Thus, there exists a sequence $R_n \to \infty$ such that
			\begin{equation} \label{eq:bdry_vanish}
				\lim_{n \to \infty} R_n K(R_n) = 0 \quad \text{and} \quad \lim_{n \to \infty} R_n K(-R_n) = 0.
			\end{equation}
			We multiply the equation $-\Delta \psi = \omega$ by $x_2 \partial_2 \psi$ and integrate over $\calC_n := \mathbb{T} \times [-R_n, R_n]$.
			For the left-hand side, integration by parts yields
			\begin{align} \label{eq:poho_lhs}
				\int_{\calC_n} (-\Delta \psi)(x_2 \partial_2 \psi) \, dx &= \frac{1}{2} \int_{\calC_n} |\partial_2 \psi|^2 - |\partial_1 \psi|^2 dx + B_n,
			\end{align}
			where the boundary term $B_n = \left[ \frac{x_2}{2} \int_{\mathbb{T}} |\partial_1\psi|^2 - |\partial_2 \psi|^2 dx_1  \right]_{x_2 = -R_n}^{x_2 = R_n}$ satisfies $|B_n| \lesssim R_n K(R_n) + R_n K(-R_n) \to 0$ as $n \to \infty$ due to \eqref{eq:bdry_vanish}.
			
			For the right-hand side, taking the (weak) derivatives of $\Psi_\pm$ yields $\partial_2 \Psi_\pm = (\pm\partial_2 \psi - W\text{sgn}(x_2)) \mathbf{1}_{\{ \Psi_\pm > 0 \}}$.
			Multiplying by $\varepsilon^{-2} x_2$ and using $x_2 \text{sgn}(x_2) = |x_2|$ as well as $\varepsilon^{-2}\mathbf{1}_{\{ \Psi_\pm > 0 \}} = \omega_\pm$, we get
			\begin{align*}
				\pm \omega_{\pm} (x_2 \partial_2 \psi) &= \varepsilon^{-2} x_2 \partial_2 \Psi_{\pm} + W|x_2|\omega_{\pm}. %\\				(-\omega_-) (x_2 \partial_2 \psi) &= \varepsilon^{-2} x_2 \partial_2 \Psi_- + W|x_2|\omega_-.
			\end{align*}
			Adding these two identities and using $\omega = \omega_+ - \omega_-$, $|\omega| = \omega_+ + \omega_-$, we obtain
			\begin{equation*}
				\omega (x_2 \partial_2 \psi) = x_2 \partial_2 (\varepsilon^{-2}\Psi_+ + \varepsilon^{-2}\Psi_-) + W|x_2||\omega|.
			\end{equation*}
			Integrating this over $\calC_n$ and applying integration by parts to the first term, we have
			\begin{equation} \label{eq:poho_rhs}
				\int_{\calC_n} \omega (x_2 \partial_2 \psi) \, dx = - \int_{\calC_n} (\Psi_+ \omega_+ + \Psi_- \omega_-) \, dx + W \int_{\calC_n} |x_2 \omega| \, dx + C_n,
			\end{equation}
			where we used $\varepsilon^{-2}\Psi_\pm = \Psi_\pm \omega_\pm$.
			The boundary term
			\[
			C_n = \left[ x_2 \int_{\mathbb{T}} (\Psi_+ \omega_+ + \Psi_- \omega_-) \, dx_1 \right]_{x_2=-R_n}^{x_2=R_n},
			\]
			satisfies $C_n \le R_n K(R_n) + R_n K(-R_n) \to 0$ as $n \to \infty$.
			Letting $n \to \infty$, \eqref{eq:poho_rhs} becomes
			\begin{equation} \label{eq:poho_rhs_2}
				\int_{\calC} \omega(x_2 \partial_2 \psi) \, dx = - \int_{\calC} (\Psi_+ \omega_+ + \Psi_- \omega_-) \, dx + W \|x_2\omega\|_1.
			\end{equation}
			Here, we observe that
			\begin{align*}
				\int_{\calC} (\Psi_+ \omega_+ + \Psi_- \omega_-) \,dx &= \int_{\calC} (\psi - W|x_2| - \gamma_+)\omega_+ \,dx + \int_{\calC} (-\psi - W|x_2| - \gamma_-)\omega_- \,dx \\
				&= \int_{\calC} \psi \omega \,dx - W \int_{\calC} |x_2||\omega| \,dx - \gamma_+\|\omega_+\|_1 - \gamma_-\|\omega_-\|_1 \\
				&= 2E[\omega] - W\|x_2\omega\|_1 - M(\gamma_+ + \gamma_-).
			\end{align*}
			Substituting this back into \eqref{eq:poho_rhs_2} and using \eqref{eq:poho_lhs}, we deduce
			\begin{equation*}
				\frac{1}{2} \int_{\calC} |\partial_2 \psi|^2 - |\partial_1 \psi|^2 dx  =- 2E[\omega] + 2W\mu + M(\gamma_+ + \gamma_-).
			\end{equation*}
			This proves \eqref{eq:poho}.
		\end{proof}
		
		\section{Quantitative estimates}\label{sec:unif}
		
		In this section, we fix $\mu > 0$ and obtain quantitative estimates for $\omega \in \calS_{\mu,\varepsilon}$ which hold for all sufficiently small $\varepsilon$. From the previous section, we know that any $\omega \in \calS_{\mu,\varepsilon}$ satisfies, for uniquely defined $W, \gmm_{\pm}$, \begin{equation*}
			\begin{split}
				\omg = \varepsilon^{-2}\left( \mathbf{1}_{\Omg_+} - \mathbf{1}_{\Omg_-} \right), \qquad \Omg_{\pm} := \left\{ x \in \calC \, : \, \pm \psi[\omg](x) - W|x_{2}| - \gmm_{\pm} > 0 \right\}. 
			\end{split}
		\end{equation*} The main results of this section are summarized in the following proposition. 
		
		\begin{proposition}\label{prop:unif}
			There exist $\varepsilon_0 = \varepsilon_0(\mu) > 0$ such that the following holds for any $\varepsilon \in (0,\varepsilon_0)$ and $\omega \in \calS_{\mu,\varepsilon}$. \begin{enumerate}
				\item[(i)] Lagrange multiplier estimates: for some $C = C(\mu)>0$,   \begin{equation*}
					\begin{split}
						0 \le W \le C, \qquad \left| \gmm_{\pm} - \frac1{4\pi} |\log \varepsilon  | \right| \le C. 
					\end{split}
				\end{equation*}
				\item[(ii)] Estimates on $\Omg_{\pm}$: for some $C, \bar{R}, c > 0$ depending only on $\mu$,  \begin{itemize}
					\item Mass saturation: $|\Omg_{\pm}| = \varepsilon^2/2$. 
					\item Uniform compact support: $\supp(\Omg_{\pm}) \subset \left\{ |x_{2}| \le \bar{R} \right\}$. 
					\item Diameter bound: $\mathrm{diam}(\Omg_{\pm}) \le C \varepsilon$. 
					\item Separation: $\mathrm{dist}(\Omg_+, \Omg_-) \ge c$.
				\end{itemize}
			\end{enumerate}
		\end{proposition}
		
		\begin{proof}
			The item (i) follows from Lemmas \ref{lem:W} and \ref{lem:gamma}. The four estimates in the item (ii) are proved in Lemmas \ref{lem:mass_sat}, \ref{lem:supp}, \ref{lem:core}, and \ref{lem:sep}, respectively.
		\end{proof}

		\subsection{Weak concentration} We prove some preliminary lower bounds on the $L^{1}$ norm.
		
		\begin{lemma}\label{lem:mass_asymp}
			There exist $\varepsilon_0 = \varepsilon_0(\mu) > 0$ and $C = C(\mu)>0$ such that for any $\varepsilon \in (0,\varepsilon_0)$ and any $\omega \in \calS_{\mu,\varepsilon}$, 
			\[
			\|\omega_+\|_1^2 = \|\omega_-\|_1^2 \ge \frac14-\frac{C}{|\log\varepsilon|}.
			\]
			In particular,
			\[
			\lim_{\varepsilon \to 0}\inf_{\omega\in\calS_{\mu,\varepsilon}} \|\omega_+\|_1 = \lim_{\varepsilon \to 0}\inf_{\omega\in\calS_{\mu,\varepsilon}} \|\omega_-\|_1 = \frac{1}{2}.
			\]
		\end{lemma}
		
		\begin{proof}
			Let $M:=\|\omega_+\|_1=\|\omega_-\|_1\le 1/2$.
			Recall that $E[\omega] = E[\omega_+] + E[\omega_-] - E_{\mathrm{int}}[\omega_+,\omega_-]$.
			By Lemma~\ref{lem:patch}, we know that $\omega_\pm=\varepsilon^{-2}\mathbf 1_{\Omega_\pm}$.
			Applying Lemma~\ref{lem:en_pos} to each of $\omega_\pm$, we obtain
			\[
			E[\omega_+] +  E[\omega_-] \le \frac{M^2}{2\pi} |\log\varepsilon| + C.
			\]
			By \eqref{green_lin} and \eqref{green_sing_2} in Lemma~\ref{lem_green_decom}, we have $-G(x-y)\lesssim 1+|x_2|+|y_2|$ on $\calC\times\calC$.
			Therefore,
			\begin{align*}
				-E_{\mathrm{int}}[\omega_+,\omega_-] &\le C\iint_{\calC\times\calC} \big( 1 + |x_2| + |y_2| \big) \omega_+(x)\omega_-(y)\,dxdy \le C,
			\end{align*}
			since $M \le 1/2$ and $\|x_2\omega\|_1 = \mu$.
			Consequently,
			\[
			E[\omega] \le \frac{M^2}{2\pi}|\log\varepsilon| + C.
			\]
			On the other hand, Lemma~\ref{lem:en_low} gives $E[\omega]\ge \frac{1}{8\pi}|\log\varepsilon| - C$.
			Combining the two estimates, we obtain
			\[
			\frac{1}{8\pi} |\log\varepsilon| - C \le \frac{M^2}{2\pi}|\log\varepsilon| + C \quad  \Longrightarrow \quad M^2 \ge \frac{1}{4} - \frac{C(\mu)}{|\log \varepsilon|}.
			\]
			Since $M \le 1/2$, this implies $M \to 1/2$ as $\varepsilon \to 0$ along any sequence of maximizers.
		\end{proof}

		\begin{lemma} \label{lem:weak_bulk}
			There exist $\varepsilon_0 = \varepsilon_0(\mu) > 0$ and $R = R(\mu) > 0$ such that the following holds:
            for any $\varepsilon \in (0, \varepsilon_0)$ and $\omega =  \omega_+ - \omega_- \in \calS_{\mu,\varepsilon}$, we have
			\begin{equation}\label{eq:bulk_mass}
				\int_{\Omega_+\cap\{|x_2|\le R\}}\omega_+\,dx > \frac{1}{4} - \frac{C}{|\log\varepsilon|}, \quad \int_{\Omega_-\cap\{|x_2|\le R\}}\omega_-\,dx > \frac{1}{4} - \frac{C}{|\log\varepsilon|}
			\end{equation}
            for some $C = C(\mu) > 0$.
		\end{lemma}
		
		\begin{proof}
			It suffices to prove the estimate for $\omega_+$.
			By Lemma~\ref{lem:mass_asymp}, we have
			\[
			M:=\|\omega_+\|_1=\|\omega_-\|_1 \ge \frac{1}{2}-\frac{C}{|\log\varepsilon|}.
			\]
			For any $R>0$,
			\[
			\int_{\Omega_+\cap\{|x_2|>R\}}\omega_+\,dx \le \frac{1}{R} \int_{\Omega_+}|x_2|\omega_+\,dx \le \frac{\mu}{R}.
			\]
			Hence
			\[
			\int_{\Omega_+\cap\{|x_2|\le R\}}\omega_+\,dx \ge M - \frac{\mu}{R} \ge \frac{1}{2} - \frac{\mu}{R}-\frac{C}{|\log\varepsilon|}.
			\]
			Choose $R>4\mu$ so that $1/2 - \mu/R > 1/4$.
			By enlarging $C = C(\mu)$ if necessary, we obtain the first inequality in \eqref{eq:bulk_mass}.
			The proof for the second inequality is identical.
		\end{proof}

		\subsection{Uniform Lagrange multipliers}

		\begin{lemma}[$W = O(1)$]\label{lem:W}
			There exist $\varepsilon_0 = \varepsilon_0(\mu)$ and $C = C(\mu) > 0$ such that the following holds:
			for any $\varepsilon \in  (0, \varepsilon_0)$ and any $\omega \in \calS_{\mu, \varepsilon}$, its Lagrange multiplier $W = W(\omega) \ge 0$ in \eqref{eq:EL} satisfies
			\begin{equation}\label{eq:W_unif}
				0\le W \le C.
			\end{equation}
		\end{lemma}
		
		\begin{proof}
			If $W = 0$, there is nothing to prove.
			Let $W > 0$.
			By \eqref{eq:poho_rhs_2} in the proof of Lemma~\ref{lem:pohozaev_identity}, we have
			\begin{align*}
				0 < W\mu = W\|x_2\omega\|_1 &= \int_{\calC} x_2\omega \partial_2 \psi \,dx + \int_{\calC} \Psi_+ \omega_+ \, dx + \int_{\calC} \Psi_- \omega_- \, dx  =: I_1 + I_2 + I_3,
			\end{align*}
			where $\Psi_\pm(x) = \max\{\pm\psi - W|x_2| - \gamma_\pm, 0\}$.
			
			\medskip
			
			\noindent \textit{Estimate of $I_1$}: Since $\psi(x) = \int G(x, y) \omega(y) dy$, we have
			\[
			I_1 = \iint_{\calC \times \calC} x_2 \partial_{x_2} G(x,y) \omega(x) \omega(y) \, dx dy = \iint_{\calC \times \calC} y_2 \partial_{y_2} G(y,x) \omega(x) \omega(y) \, dx dy.
			\]
			Due to the symmetry $G(x,y) = G(y,x)$ and the translation invariance $x_2$, we have $\partial_{y_2} G(y,x) = \partial_{y_2} G(x,y) = -\partial_{x_2} G(x,y)$.
			This yields
			\begin{equation} \label{eq:I_1_sym}
				I_1 = \frac{1}{2} \iint_{\calC \times \calC} (x_2 - y_2) \partial_{x_2} G(x,y) \omega(x) \omega(y) \, dx dy.
			\end{equation}
			Observe that, near $x = y$, $G(x,y) \simeq -\frac{1}{2\pi} \log |x-y|$, and $G(x,y) \simeq -|x_2 - y_2|$ when $|x_2 - y_2| \gg 1$ (see Lemmas~\ref{lem_green_local} and~\ref{lem_green_decom}).
			Thus we have
			\[
			|(x_2 - y_2) \partial_{x_2} G(x,y)| \lesssim  1 + |x_2 - y_2|  \lesssim  1 + |x_2| + |y_2|,
			\]
			which, together with \eqref{eq:I_1_sym}, implies
			\begin{align*}
				|I_1| \lesssim  \|\omega\|_{1}^2 +  \|\omega\|_{1} \|x_2\omega\|_1 \lesssim 1 + \mu.
			\end{align*}
			
			\bigskip
			
			\noindent \textit{Estimate of $I_2$ and $I_3$}:
			We only consider $I_2$ since $I_3$ can be dealt in the same way.
			Let $\Omega_{+} = \{\omega_+ > 0\} $.
			Note that $\omega_+ = \varepsilon^{-2} \mathbf{1}_{\Omega_+}$, and $|\Omega_+| = \varepsilon^2 \|\omega_+\|_1 \le \varepsilon^2/2$.
			On $\Omega_+=\{\Psi_+>0\}$, we have $\Psi_+=\psi-W|x_2|-\gamma_+$, and hence, in the sense of distribution,
			$$
			-\Delta\Psi_+=\varepsilon^{-2}+2W\delta_{\{x_2=0\}} \quad \textrm{in }\,\, \Omega_+.
			$$
			where $\delta$ is the Dirac delta.
			Let $u \in H^1_0(\Omega_+)$ be the weak solution to $-\Delta u = \varepsilon^{-2}$ in $\Omega_+$.
			Since $\Psi_+$ is also in $H^1_0(\Omega_+)$, we have
			\begin{align*}
				I_2 = &\int_{\Omega_+} \Psi_+ \varepsilon^{-2} \, dx  = \int_{\Omega_+} \nabla \Psi_+  \cdot \nabla  u \, dx\\ 
				= &\int_{\Omega_+} u\, \varepsilon^{-2} \, dx + 2W \int_{\Omega_+ \cap \{x_2 = 0\}} u(x_1, 0) \, dx_1  =: J_1 + J_2.
			\end{align*}
			By Talenti's rearrangement principle (see \cite{Tal1976} and Remark~\ref{rmk:talenti}), for any $p \in [1, \infty]$ we have
			\[
			\| u \|_{L^p(\Omega_+)} \le \| v \|_{L^p(\Omega^*)}, \quad \textrm{where} \quad v(r) = \frac{\varepsilon^{-2}}{4}\left( \frac{|\Omega_+|}{\pi} - r^2 \right).
			\]
			Here, $v$ is the radially symmetric solution to $-\Delta v = \varepsilon^{-2}$ on the disk $\Omega^*$ such that $|\Omega^*| = |\Omega_+|$.
			Therefore,
			\[
			J_1 = \int_{\Omega_+} u \varepsilon^{-2} \, dx \le \int_{\Omega^*} v \varepsilon^{-2} \, dx = \frac{\varepsilon^{-4}}{8\pi}|\Omega_+|^2  \le \frac{1}{32\pi}.
			\]
			Moreover,
			\begin{equation}\label{eq:u_inf_Tal}
				\|u\|_{L^\infty(\Omega_+)} \le \|v\|_{L^\infty(\Omega^*)} = v(0) = \frac{\varepsilon^{-2}}{4\pi}|\Omega_+| \le \frac{1}{8\pi}.
			\end{equation}
			
			For $J_2$, we extend $u$ by zero to the entire $\calC$ and use the same notation for the extension.
			Note that the trace $u(x_1, 0)$ is well-defined and satisfies
			\[
			u(x_1, 0)^2 = -2 \int_0^\infty u(x_1, x_2) \partial_{x_2} u(x_1, x_2) \, dx_2.
			\]
			Integrating over $\Omega_+ \cap \{x_2 = 0\} \subset \bbT$, we obtain
			\begin{align*}
				\int_{\Omega_+ \cap \{x_2 = 0\}} u(x_1, 0)^2 \, dx_1 \le 2 \int_{\bbT} \int_0^\infty |u| |\nabla u| \, dx_2 dx_1  \le 2 \|u\|_{L^2(\Omega_+)} \|\nabla u\|_{L^2(\Omega_+)}.
			\end{align*}
			Observe that
			\begin{align*}
				\|\nabla u\|_{L^2(\Omega_+)}^2 &= \int_{\Omega_+} u \varepsilon^{-2} \, dx = J_1 \le \frac{1}{32\pi}, \\
				\|u\|_{L^2(\Omega_+)}^2 &\le \|u\|_{L^\infty(\Omega_+)}^2 |\Omega_+| \le \left( \frac{1}{8\pi} \right)^2 \frac{\varepsilon^2}{2} = \frac{\varepsilon^2}{128\pi^2},
			\end{align*}
			where we used \eqref{eq:u_inf_Tal} in the second inequality.
			Combining these yields $\int_{\Omega_+ \cap \{x_2 = 0\}} u(x_1, 0)^2 \, dx_1 \lesssim  \varepsilon$.
			By H\"older's inequality on $\Omega_+ \cap \{x_2 = 0\} \subset \bbT$, we get
			\[
			\int_{\Omega_+ \cap \{x_2 = 0\}} u(x_1, 0) \, dx_1 \le |\Omega_+ \cap \{x_2 = 0\}|^{1/2} \left( \int_{\Omega_+ \cap \{x_2 = 0\}} u(x_1, 0)^2 \, dx_1 \right)^{1/2} \lesssim  \varepsilon^{1/2}.
			\]
			Consequently, $I_2 = J_1 + J_2 \le \frac{1}{32\pi} +  C W \varepsilon^{1/2}$. 
			Applying the identical procedure for $I_3$, we find $I_3 \le \frac{1}{32\pi} +  C W \varepsilon^{1/2}$.
			
			\medskip
			
			\noindent Summing the estimates gives
			\[
			W\mu = W \|x_2 \omega\|_1 = I_1 + I_2 + I_3 \le C(1+\mu) + \frac{1}{16\pi} +  C W \varepsilon^{1/2},
			\]
			that is, $ W \left( \mu-  C \varepsilon^{1/2} \right) \le C(1+\mu)$.
			Then, the desired inequality \eqref{eq:W_unif} holds provided that $  C\varepsilon^{1/2} \le \mu/2$.
		\end{proof}
		
		\begin{remark}\label{rmk:talenti}
			The use of Talenti's principle on $\calC = \bbT \times \bbR$ is justified by the isoperimetric inequality:
			\[
			\operatorname{Per}(E) \ge \min \left\{ \sqrt{4\pi |E|}, 4\pi \right\},
			\]
			where $\operatorname{Per}(E)$ denotes the perimeter of $E$.
			In particular, if $|E| \le 4\pi$, then $\operatorname{Per}(E)^2 \ge 4 \pi |E|$.
			Talenti's argument \cite{Tal1976} depends only on this inequality.
			Since $|\Omega_{\pm}| \le \varepsilon^2/2 \le 4\pi$ in the proof of Lemma~\ref{lem:W}, Talenti's original argument applies without change.
		\end{remark}
		
		\begin{lemma}[$\gamma_{\pm} = O(|\log \varepsilon|)$]\label{lem:gamma}
			There exist $\varepsilon_0 = \varepsilon_0(\mu) > 0$ and $C = C(\mu) > 0$ such that the following holds:
			For any $\varepsilon \in (0, \varepsilon_0)$ and any $\omega \in \calS_{\mu, \varepsilon}$, its Lagrange multipliers $\gamma_{\pm} =  \gamma_{\pm}(\omega)$ in \eqref{eq:EL} satisfy
			\begin{equation}\label{eq:gam}
				\frac{1}{4\pi}|\log \varepsilon| - C \le \gamma_{\pm} \le \frac{1}{4\pi}|\log \varepsilon| + C.
			\end{equation}
		\end{lemma}
		
		\begin{proof}
			Let $M:=\|\omega_+\|_1=\|\omega_-\|_1 \le 1/2$.
			By Lemma~\ref{lem:mass_asymp}, for $\varepsilon > 0$ sufficiently small, we may assume that $M\ge 1/4$.
			By \eqref{eq:en_low} of Lemma~\ref{lem:en_low}, we have
			$
			E[\omega]\ge \frac{1}{8\pi}|\log\varepsilon|-C.
			$
			
			\medskip
			
			\noindent \textbf{Lower bound for $\gamma_++\gamma_-$.}
			Let $\Psi_{\pm}= \max\{\pm\psi-W|x_2|-\gamma_\pm, 0\}$.
			Using \eqref{eq:EL}, we have
			\begin{align*}
				2E[\omega]
				&=\int_{\Omega_+}\psi\omega_+\,dx-\int_{\Omega_-}\psi\omega_-\,dx \\
				&=M(\gamma_++\gamma_-)+W\mu+\int_{\calC}(\Psi_+\omega_+ + \Psi_-\omega_-)\,dx.
			\end{align*}
			By Lemma~\ref{lem:W}, $W\le C(\mu)$. Moreover, as in the proof of Lemma~\ref{lem:W}, the rearrangement argument gives
			$$
			\int_{\calC}(\Psi_+\omega_+ + \Psi_-\omega_-)\,dx\le C(\mu).
			$$
			Hence,
			$$
			M(\gamma_+ + \gamma_-) = 2E[\omega]-W\mu-\int_{\calC}(\Psi_+\omega_+ + \Psi_-\omega_-)\,dx \ge \frac{1}{4\pi}|\log\varepsilon| - C.
			$$
			Since $M \le 1/2 $, this implies
			\begin{equation}\label{eq:gam_sum_low}
				\gamma_+ + \gamma_- \ge \frac{1}{2\pi}|\log\varepsilon| - C.
			\end{equation}
			
			\medskip
			
			\noindent \textbf{Upper bound for $\gamma_{\pm}$.}
			We estimate $\gamma_-$; the proof for $\gamma_+$ is identical. For $x\in\Omega_-$, \eqref{eq:EL} gives $\gamma_-\le -\psi(x)-W|x_2|$. Multiplying by $\omega_-$ and integrating over $\Omega_-$, we obtain
			\begin{equation}\label{eq:gam_avg}
				M\gamma_-\le -\int_{\Omega_-}\psi\omega_-\,dx
				=\int_{\Omega_-}\psi[\omega_-]\omega_-\,dx-\int_{\Omega_-}\psi[\omega_+]\omega_-\,dx,
			\end{equation}
			where we dropped the non-positive term $-W\int_{\Omega_-}|x_2|\omega_- dx$.
			For the self-interaction term, Lemma~\ref{lem:en_pos} gives
			$$
			\int_{\Omega_-}\psi[\omega_-]\omega_-\,dx=2E[\omega_-]\le \frac{M^2}{2\pi}|\log\varepsilon|+C.
			$$
			For the interaction term, using \eqref{green_lin} and \eqref{green_sing_2} of Lemma~\ref{lem_green_decom}, we have $G(x,y)\gtrsim -(1+|x_2|+|y_2|)$.
			Therefore
			$$
			\int_{\Omega_-}\psi[\omega_+]\omega_-\,dx\ge -C(1+\mu).
			$$
			Substituting these estimates into \eqref{eq:gam_avg}, and using $1/4 \le M\le 1/2$, we get
			$$
			\gamma_-\le \frac{M}{2\pi}|\log\varepsilon|+C(\mu)\le \frac{1}{4\pi}|\log\varepsilon|+C.
			$$
			The same argument gives
			\begin{equation}\label{eq:gam_up_2}
				\gamma_+\le \frac{1}{4\pi}|\log\varepsilon|+C.
			\end{equation}
			
			\medskip
			
			\noindent Finally, combining \eqref{eq:gam_sum_low} and \eqref{eq:gam_up_2}, we obtain
			$$
			\gamma_+\ge \left(\frac{1}{2\pi}|\log\varepsilon|-C\right)-\left(\frac{1}{4\pi}|\log\varepsilon|+C\right)
			=\frac{1}{4\pi}|\log\varepsilon|-C,
			$$
			and the same lower bound holds for $\gamma_-$. This proves \eqref{eq:gam}.
		\end{proof}
		
		\begin{lemma}[Mass saturation]\label{lem:mass_sat}
			There exists $\varepsilon_0 = \varepsilon_0(\mu)>0$ such that for any $\varepsilon\in(0,\varepsilon_0)$ and any $\omega\in\calS_{\mu,\varepsilon}$, we have $\|\omega_+\|_1=\|\omega_-\|_1 = 1/2$.
		\end{lemma}
		
		\begin{proof}
			Let $M:=\|\omega_+\|_1=\|\omega_-\|_1 \le 1/2$.
			Suppose for contradiction that $M < 1/2$.
			
			Let $h_{n,1}^-$ and $k_{n,1}^-$ be the auxiliary functions for $\Omega_+$ and $\Omega_-$ used in the proof of Lemma~\ref{lem:local_EL}, respectively.
			Define $\eta_n:=h_{n,1}^- - k_{n,1}^-$.
			Then $\int \eta_n  = \int|x_2|\eta_n =0$.
			Moreover, $\eta_n \le 0$ on $\Omega_+$ and $\eta_n \ge 0$ on $\Omega_-$.
			Observe that, since $M<1/2$, the $L^1$ constraint is unsaturated.
			Hence, by the same admissibility argument as in the proof of Lemma~\ref{lem:mass_multip}, we have $\omega + \sigma\eta_n \in \calF_{\mu,\varepsilon}$ for all sufficiently small $\sigma > 0$.
			By maximality,
			$$
			0\ge \int_{\calC}\psi[\omega]\eta_n\,dx.
			$$
			Taking $n\to\infty$ and using \eqref{eq:aux_lim_positive} and \eqref{eq:aux_lim_negative}, we obtain $0\ge \gamma_+ + \gamma_-$.
			This contradicts \eqref{eq:gam} of Lemma~\ref{lem:gamma} for all sufficiently small $\varepsilon$.
			Therefore $M=1/2$.
		\end{proof}
		
		\subsection{Localization of patches}
		
		\begin{lemma}[Uniform compact support] \label{lem:supp}
			There exist $\varepsilon_0 = \varepsilon_0(\mu) > 0$ and $\bar{R} = \bar{R}(\mu) > 0$ such that the following holds:
			for any $\varepsilon \in (0, \varepsilon_0)$ and any $\omega \in \calS_{\mu, \varepsilon}$, we have
			\begin{equation}\label{eq:supp}
				\supp \omega \subset \{ x \in \calC :  |x_2| \le \bar{R} \}.
			\end{equation}
		\end{lemma}
		
		\begin{proof}
			We only treat $\Omega_+ = \supp \omega_+$ since the proof for $\Omega_-$ is symmetric.
			Let $R > 0$ be the constant from Lemma~\ref{lem:weak_bulk}.
			We claim that $\Omega_+ \subset \{|x_2| \le R + 1\}$ holds for all sufficiently small $\varepsilon$.
			Suppose for contradiction that there exists $x \in \Omega_+$ such that $|x_2| > R + 1$.
			By \eqref{eq:EL} with $W \ge 0$ and Lemma~\ref{lem:gamma}, 
			\begin{equation} \label{eq:el_psi_low}
				%\psi(x) > \gamma_+\ge \frac{\|\omega\|_1}{4\pi}|\log \varepsilon| - C.
				\psi(x) > \gamma_+\ge \frac{1}{4\pi}|\log \varepsilon| - C.
			\end{equation}
			Next, we establish an upper bound for $\psi(x)$. 
			Using Lemma~\ref{lem_green_decom}, we split $\psi$ into two terms:
			\begin{align*}
				\psi(x) = -\frac{1}{4\pi} \int_{\calC} |x_2 - y_2|\omega(y) \,dy + \int_{\calC} H(x-y)\omega(y) \,dy.
			\end{align*}
			For the first term, we use $-|x_2 - y_2| \le -|x_2| + |y_2|$ and $|x_2 - y_2| \le |x_2| + |y_2|$.
			Since $\|\omega_+\|_1 = \|\omega_-\|_1$, we find
			\begin{align*}
				-\frac{1}{4\pi}\int_{\calC} |x_2 - y_2|\omega(y) \,dy &= -\frac{1}{4\pi}\int_{\Omega_+} |x_2 - y_2|\omega_+(y) \,dy + \frac{1}{4\pi}\int_{\Omega_-} |x_2 - y_2|\omega_-(y) \,dy \\
				&\le \frac{1}{4\pi}\int_{\Omega_+} \big( -|x_2| + |y_2| \big)\omega_+(y) \,dy + \frac{1}{4\pi}\int_{\Omega_-} \big( |x_2| + |y_2| \big)\omega_-(y) \,dy \\
				&=  \frac{1}{4\pi}\int_{\calC} |y_2||\omega(y)| \,dy = \frac{\mu}{4\pi}.
			\end{align*}
			For the second term,
			\[
			\int_{\calC} H(x-y)\omega(y) \,dy = \int_{\Omega_+} H(x-y)\omega_+(y) \,dy - \int_{\Omega_-} H(x-y)\omega_-(y) \,dy,
			\]
			the lower bound $H(x-y) \ge - \frac{\log 2}{2\pi}$ from \eqref{green_sing_2} ensures
			\[
			-\int_{\Omega_-} H(x-y)\omega_-(y) \,dy \le \frac{\log 2}{2\pi} \|\omega_-\|_1 \le C.
			\]
			For the contribution from $\Omega_+$, we split the domain into $A := \Omega_+ \cap \{|y_2| \le R\}$ and $B := \Omega_+ \setminus A$. 
			For $y \in A$, we have $|x_2 - y_2| \ge |x_2| - |y_2| > (R+1) - R = 1$.
			By the upper bound in \eqref{green_sing_2},
			\[
			\int_A H(x-y)\omega_+ \,dy  \le C\|\omega_+\|_1 \le C.
			\]
			For the tail $B$, the bound \eqref{green_sing_3} and the rearrangement of $\omega_+ \mathbf{1}_{B} = \varepsilon^{-2} \mathbf{1}_{\Omega_+\cap B}$ yields
			\[
			\int_B H(x-y)\omega_+ \,dy \le \frac{\|\omega_+\|_{L^1(B)}}{2\pi}|\log \varepsilon| + C.
			\]
			By Lemma~\ref{lem:weak_bulk}, we have $\|\omega_+\|_{L^1(A)} > 1/4 -C|\log\varepsilon|^{-1}$.
			Since $\|\omega_+\|_1\le 1/2$, it follows that $\|\omega_+\|_{L^1(B)} < 1/4 + C|\log\varepsilon|^{-1}$.
			This implies
			\[
			\int_B H(x-y)\omega_+ \,dy \le \frac{1}{8\pi}|\log \varepsilon| + C.
			\]
			Collecting all the estimates for $\psi(x)$ and combining them with \eqref{eq:el_psi_low}, we obtain
			\[
			\frac{1}{4\pi}|\log \varepsilon| - C <\psi(x) \le  \frac{1}{8\pi}|\log \varepsilon| + C,
			\]
			that is, $|\log \varepsilon| \le C'$, which is a contradiction for all sufficiently small $\varepsilon > 0$.
			Therefore, $\Omega_+ \subset \{|x_2| \le R + 1\}$. Taking $\bar{R} = R + 1$ concludes the proof.
		\end{proof}
		
		\begin{lemma}[Core size]\label{lem:core}
			Let  $\varepsilon_0 = \varepsilon_0 (\mu) > 0$ be sufficiently small so that Lemmas~\ref{lem:W},~\ref{lem:gamma}, and~\ref{lem:supp} apply.
			For any $\varepsilon \in (0, \varepsilon_0)$ and any $\omega =  \varepsilon^{-2}(\mathbf{1}_{\Omega_+} - \mathbf{1}_{\Omega_-}) \in \calS_{\mu, \varepsilon}$, we have $\operatorname{diam} (\Omega_{\pm}) \le C \varepsilon$
			where $C = C(\mu) > 0$.
		\end{lemma}
		
		\begin{proof}
			We present the proof only for $\Omega_+$, since the argument for $\Omega_-$ is symmetric.
			By Lemma~\ref{lem:supp}, there exists $R = R(\mu) > 0$ such that $\supp \omega \subset \bbT \times [-R, R]$.
			%Let $M := \|\omega_+\|_1 = \|\omega_-\|_1$ so that $\|\omega\|_1 = 2M$.
			For any $x \in \Omega_+$, \eqref{eq:EL} and Lemma~\ref{lem:gamma} imply
			\[
			\psi[\omega](x) - W|x_2| > \gamma_+ \ge \frac{1}{4\pi}|\log \varepsilon| - C.
			\]
			Since $|x_2| \le R$ and $W = O(1)$ by Lemma~\ref{lem:W}, the term $W|x_2|$ is uniformly bounded.
			Thus, recalling $\psi[\omega] = \psi[\omega_+] - \psi[\omega_-]$, we have
			\[
			\psi[\omega_+](x) \ge \psi[\omega_-](x) + \frac{1}{4\pi} |\log \varepsilon| - C'.
			\]
			By Lemma~\ref{lem_green_decom}, we have the  lower bound 
			\[
			G(x, y) \ge -C_0 - C_1 \big( |x_2| + |y_2| \big) \ge -C(R) \quad \textrm{ for } \,\, x, y \in \supp \omega.
			\]
			Therefore, we have $\psi[\omega_-](x) \ge -C$.
			This implies that for any $x \in \Omega_+$,
			\begin{equation}\label{eq:psi_lower_2}
				\psi[\omega_+](x) \ge \frac{1}{4\pi} |\log \varepsilon| - C,
			\end{equation}
			where $C > 0$ is independent of $\varepsilon$.
			
			\medskip
			
			\noindent We now establish the upper bound for $d := \operatorname{diam}( \Omega_+)$.
			By Lemma~\ref{lem:supp}, $d < \infty$.
			We take two points $x, y \in \overline{\Omega}_+$ such that $|x - y| = d$.
			By \eqref{eq:psi_lower_2} and the continuity of $\psi[\omega]$, we have
			\begin{equation}\label{eq:psi_lower_3}
				\psi[\omega_+](x) + \psi[\omega_+](y) \ge \frac{1}{2\pi} |\log \varepsilon| - 2C.
			\end{equation}
			On the other hand, since $\Omega_+$ is bounded, Lemma~\ref{lem_green_local} ensures that $G(x, z) \le -\frac{1}{2\pi}\log|x-z| + C$ for any $x, z \in \Omega_+$.
			Since $\omega_+(z) = \varepsilon^{-2}$ on $\Omega_+$, we obtain
			\begin{align*}
				\psi[\omega_+](x) + \psi[\omega_+](y) &\le \underbrace{\frac{\varepsilon^{-2}}{2\pi} \int_{\Omega_+} \left( \log \frac{1}{|x-z|} + \log \frac{1}{|y-z|} \right)\, dz}_{=: I(\varepsilon)} + C.
			\end{align*}
			To estimate $I(\varepsilon) = I$, we split the domain $\Omega_+$ by
			\[
			\Omega_+ = \{z \in \Omega_+ : |x-z| \le |y-z| \} \, \cup \, \{z \in \Omega_+ : |x-z| > |y-z| \} =: \Omega_1 \cup \Omega_2.
			\]
			Let $I_1$ and $I_2$ be the integrals over $\Omega_1$ and $\Omega_2$, respectively. 
			Since $d = |x - y| \le |x - z| + |y - z|$, for any $z \in \Omega_1$ we have $|y - z| \ge d/2$. 
			Thus, 
			\[
			I_1 = \int_{\Omega_1} \log \frac{1}{|x-z|} \, dz + \int_{\Omega_1} \log \frac{1}{|y-z|} \, dz \le \int_{\Omega_1} \log \frac{1}{|x-z|} \, dz + |\Omega_1| \log \frac{2}{d}.
			\]
			Note that the remaining integral is maximized when $\Omega_1$ is the disk $B_{R_1}(x)$ with area $\pi R_1^2 = |\Omega_1|$, that is,
			\[
			\int_{\Omega_1} \log \frac{1}{|x-z|} \, dz \le \int_{B_{R_1}(x)} \log \frac{1}{|x-z|} \, dz = \frac{|\Omega_1|}{2} + \frac{|\Omega_1|}{2} \log \frac{\pi}{|\Omega_1|}.
			\]
			Applying the symmetric argument for $I_2$ over $\Omega_2$ (where $|x-z| \ge d/2$), we obtain
			\begin{align*}
				I = \frac{\varepsilon^{-2}}{2\pi}\left( I_1 + I_2 \right) &\le \frac{\varepsilon^{-2}}{2\pi} \left\{ |\Omega_+| \log \frac{2}{d} + \frac{|\Omega_+|}{2} + \frac{1}{2} \left( |\Omega_1| \log \frac{\pi}{|\Omega_1|} + |\Omega_2| \log \frac{\pi}{|\Omega_2|} \right) \right\}.
			\end{align*}
			Note that the function $t \mapsto t \log(\pi/t)$ is concave for $t > 0$.
			Since $|\Omega_1| + |\Omega_2| = |\Omega_+| = \varepsilon^2/2$ by Lemma~\ref{lem:mass_sat}, Jensen's inequality yields
			\[
			|\Omega_1| \log \frac{\pi}{|\Omega_1|} + |\Omega_2| \log \frac{\pi}{|\Omega_2|} \le  |\Omega_+|  \log \frac{2\pi}{|\Omega_+|} = \frac{\varepsilon^2}{2} \log \frac{4\pi}{\varepsilon^2 }.
			\]
			Therefore, $I \le  \frac{1}{4\pi} \left( |\log \varepsilon| + \log \frac{2}{d} \right) + C$.
			Combining this with \eqref{eq:psi_lower_3}, we deduce
			\begin{align*}
				2 |\log \varepsilon| - C &\le  |\log \varepsilon| +  \log \frac{2}{d} + C \quad
				\Longrightarrow \quad \log \frac{1}{\varepsilon} \le \log \frac{2}{d} + C.
			\end{align*}
			Therefore, we conclude that $\operatorname{diam}(\Omega_+) = d \le C' \varepsilon$. 
		\end{proof}

		\begin{lemma}[$O(1)$-Separation]\label{lem:sep}
			There exist $\varepsilon_0=\varepsilon_0(\mu)>0$ and $c=c(\mu)>0$ such that for any $\varepsilon \in (0,\varepsilon_0)$ and any  $\omega = \varepsilon^{-2}(\mathbf 1_{\Omega_+}-\mathbf 1_{\Omega_-})\in\calS_{\mu,\varepsilon}$, we have $\dist(\Omega_+,\Omega_-) \ge c$.
		\end{lemma}
		
		\begin{proof}
			Let $d_\varepsilon:=\dist(\Omega_+,\Omega_-)$, and fix $r_0 \in (0, 1)$.
			Without loss of generality, we may assume $d_\varepsilon < r_0/2$.
			By Lemma~\ref{lem:mass_sat}, $\|\omega_+\|_1=\|\omega_-\|_1 = 1/2$.
			Recall that $E[\omega] = E[\omega_+] + E[\omega_-] - E_{\mathrm{int}}[\omega_+,\omega_-]$.
			By Lemmas~\ref{lem:en_pos} and~\ref{lem:en_low},
			$$
			E[\omega_\pm]\le \frac{1}{16\pi}|\log\varepsilon|+C, \quad E[\omega]\ge \frac{1}{8\pi}|\log\varepsilon| - C,
			$$
			and then, we have $E_{\mathrm{int}}[\omega_+,\omega_-]  \le C$. On the other hand, $\operatorname{diam}(\Omega_\pm)\le C_0\varepsilon$ by Lemma~\ref{lem:core}.
			Thus
			$$
			|x-y|\le d_\varepsilon + C'(\mu) \varepsilon \quad \textrm{for any }\,\, x \in \Omega_+, \,\, y \in \Omega_-.
			$$
			We may assume $C'(\mu) \varepsilon < r_0/2$, and hence $|x-y|<r_0$.
			Then by Lemma~\ref{lem_green_local},
			$$
			G(x,y)\ge -\frac{1}{2\pi} \log(d_\varepsilon + C'(\mu) \varepsilon) - C_1.
			$$
			Integrating over $\Omega_+\times\Omega_-$, we get
			$$
			E_{\mathrm{int}}[\omega_+,\omega_-]
			=\int_{\Omega_-}\int_{\Omega_+}G(x,y)\omega_+(x)\omega_-(y)\,dxdy \ge -\frac{1}{8\pi}\log(d_\varepsilon + C'(\mu) \varepsilon) - C_2.
			$$
			Since  $d_\varepsilon + C'(\mu) \varepsilon \le r_0 < 1$, this implies $E_{\mathrm{int}}[\omega_+,\omega_-]\ge -C \log(d_\varepsilon + C'(\mu) \varepsilon) - C_2$.
			Combining this with the upper bound $E_{\mathrm{int}}[\omega_+,\omega_-]\le C(\mu)$, we obtain
			$
			-\log(d_\varepsilon + C'(\mu)\varepsilon)\le C''(\mu).
			$
			Hence $d_\varepsilon + C'(\mu)\varepsilon\ge c_0$ for some $c_0=c_0(\mu)>0$.
			Taking $\varepsilon_0$ smaller so that $C'(\mu) \varepsilon<c_0/2$, we conclude $d_\varepsilon\ge c_0/2$.
			This proves the lemma with $c:=c_0/2$.
		\end{proof}

		\section{Emergence of von K\'arm\'an vortex street}\label{sec:karman}   
		
		{In this section, we complete the proof of Theorem \ref{thm:max}. First, in \S \ref{subsec:asym}, we prove that the support of maximizing patches converge to a pair of points in $\calC$. This pair of points corresponds to the glide reflection symmetric von K\'arm\'an vortex pair, as shown in \S \ref{subsec:local}. Then, we derive the asymptotic translation speed and simple connectedness of the vortex supports in \S \ref{subsec:speed} and \S \ref{subsec:simp-conn}, respectively. 
		}
		
		\subsection{Asymptotics as $\varepsilon \to 0$}\label{subsec:asym}
		
		\begin{lemma}\label{lem:asy_pt}
			Let $\varepsilon_j \to 0$ as $j \to \infty$ and $\omega_j = \varepsilon_j^{-2}(\mathbf{1}_{\Omega_{j, +}} - \mathbf{1}_{\Omega_{j, -}}) \in \calS_{\mu,\varepsilon_j}$.
			Then, there exist points $p_{\pm} \in \calC$ such that, up to a subsequence,
			\[
			\lim_{j \to \infty} \varepsilon_j^{-2} \int_{\Omega_{j, \pm}} \varphi(x) dx = \frac{1}{2} \varphi(p_{\pm}) \quad \textrm{ for every } \,\, \varphi \in C_c(\calC).
			\]
			Moreover, for arbitrary points $p_{j,\pm}\in \Omega_{j,\pm}$, we have $p_{j,\pm} \to p_\pm$.
		\end{lemma}
		
		\begin{proof}
			We only prove the lemma for $\varepsilon_j^{-2}\mathbf{1}_{\Omega_{j, +}} =: \omega_{j,+}$ since the proof for $\varepsilon_j^{-2}\mathbf{1}_{\Omega_{j, -}}$ is identical.
			By Lemmas~\ref{lem:supp} and~\ref{lem:core}, there exists a compact set $K\subset \mathcal C$, independent of $j$, such that $\Omega_{j,+} \subset K$.
			Also, by Lemma~\ref{lem:core}, $\operatorname{diam}(\Omega_{j,+})\lesssim_{\mu} \varepsilon_j$.
			Choose a sequence $\{q_{j, +}\}_{j}$ such that $ q_{j, +}\in \Omega_{j, +}$.
			Since $K$ is compact, by taking a subsequence still denoted by $q_{j, +}$, we may assume that $q_{j, +} \to p_{+}$ as $j \to \infty$ for some $p_+ \in K \subset \calC$.
			Moreover, for any $x \in \Omega_{j, +}$, we see
			\[
			|x - p_+| \le |x - q_{j,+}| + | q_{j,+} - p_+| \le \operatorname{diam}(\Omega_{j,+}) + | q_{j,+} - p_+| \to 0.
			\]
			
			We now prove the measure convergence.
			From the above properties, we see
			\begin{align*}
				\left|
				\int_{\calC} \left( \varphi(x) - \varphi(q_{j,+})\right) \omega_{j,+}(x)\,dx \right| \le \|\omega_{j,+}\|_{1} \sup_{x\in\Omega_{j,+}} |\varphi(x)-\varphi(q_{j,+})| \to 0, \quad \varphi \in C_c(\calC).
			\end{align*}
			By Lemma~\ref{lem:mass_sat}, for sufficiently large $j$, we have $\int \omega_{j,+}  = 1/2$, and then
			\[
			\lim_{j \to \infty} \int_{\mathcal C}\varphi(x)\omega_{j,+}(x)\,dx = \lim_{j \to \infty}  \frac{1}{2} \varphi(q_{j,+}) = \frac{1}{2} \varphi(p_+).
			\]
			This completes the proof.
		\end{proof}
		
		In what follows, whenever $p \in \Omega \subset \calC$ and $\operatorname{diam}(\Omega) < \pi$, the rescaled set $(\Omega - p)/\varepsilon$ is understood via the (local) lift of $\Omega$ and $p$ to $\bbR^2$.

		\begin{lemma}[Renormalized energy]
			\label{lem:ren_en_decom_2}
			Let $\varepsilon_j \to 0$ as $j \to \infty$, and $\omega_j = \varepsilon_j^{-2} \left( \mathbf 1_{\Omega_{j,+}} - \mathbf 1_{\Omega_{j,-}} \right) \in\calS_{\mu,\varepsilon_j}$.
			Let $p_{j,\pm}\in\Omega_{j,\pm}$ be arbitrary points and let $p_\pm\in\calC$ be the limit points, i.e. $p_{j,\pm}\to p_\pm$, ensured by Lemma~\ref{lem:asy_pt}.
			For all sufficiently large $j$, set
			%choose local lifts of $\Omega_{j, \pm}$ to $\bbR^2$
			\begin{equation}\label{eq:res_patch}
				A_{j,\pm} := \frac{\Omega_{j,\pm} - p_{j,\pm} }{\varepsilon_j}.
			\end{equation}
			Then we have $|A_{j,\pm}|= 1/2$, $\operatorname{diam}(A_{j,\pm})\le C(\mu)$, and
			\begin{equation}
				\label{eq:renorm}
				E[\omega_j] = \frac{1}{8\pi}|\log\varepsilon_j| + \mathcal{E}(A_{j,+}) + \calE (A_{j,-}) - \frac{1}{4}G(p_+-p_-) + o(1),
			\end{equation}
			where
			\begin{equation}
				\label{eq:en_log}
				\calE(A) := -\frac{1}{4\pi} \iint_{A\times A}\log|x-y|\,dxdy.
			\end{equation}
			Furthermore, $\calE(A_{j,\pm}) \le \calE(B_*) $, where $B_* := B_{1/\sqrt{2\pi}}(0)$.
		\end{lemma}

		\begin{proof}
			By Lemmas~\ref{lem:mass_sat} and~\ref{lem:core}, we see
			\begin{equation}
				\label{eq:renorm_2}
				|A_{j,\pm}| = \frac{1}{2},\quad \operatorname{diam}(A_{j,\pm}) \le C(\mu).
			\end{equation}
			We now show \eqref{eq:renorm}.
			For the self-energy term, we define
			\begin{equation}\label{eq:0806_01}
				R(z) := G(z) + \frac{1}{2\pi}\log|z| = -\frac{1}{4\pi} \log\left( \frac{2\cosh z_2-2\cos z_1}{|z|^2} \right), \quad z = (z_1, z_2) \neq 0.
			\end{equation}
			Note that the function $R$ continuously extends to $z = 0$ by setting $R(0) = 0$.
			By changing variables $x \mapsto p_{j,+} + \varepsilon_j x$ and  $y \mapsto p_{j,+} + \varepsilon_j y$, we see that
			\begin{align*}
				E[\omega_{j,+}] = & \frac{1}{2} \iint_{A_{j,+}\times A_{j,+}} G\left(\varepsilon_j(x - y)\right)\,dxdy \\
				= & \frac{|A_{j,+}|^2}{4\pi} |\log\varepsilon_j| + \calE(A_{j, +})
				+ \frac{1}{2} \iint_{A_{j,+}\times A_{j,+}} R\left(\varepsilon_j(x - y)\right)\,dxdy.
			\end{align*}
			By \eqref{eq:renorm_2} and the continuity of $R$, we have $\sup_{x,y\in A_{j,+}} \left| R(\varepsilon_j(x - y))\right| = o(1)$ as $j \to \infty$.
			By the same argument for $A_{j, -}$, we obtain
			\begin{equation}\label{eq:renorm_self}
				E[\omega_{j, \pm}] = \frac{1}{16\pi}|\log\varepsilon_j| + \calE(A_{j,\pm}) + o(1).
			\end{equation}
			
			We next consider the interaction term.
			By Lemmas~\ref{lem:asy_pt} and~\ref{lem:sep}, $p_{j, \pm} \to p_{\pm}$ for some $p_{\pm}$ such that $p_+\ne p_-$.
			Moreover, $x - y \to p_+ - p_-$ uniformly for $(x, y) \in \Omega_{j,+} \times \Omega_{j,-}$.
			Since $G(z)$ is smooth away from the origin and $\|\omega_{j,\pm}\|_1  = 1/2$, it follows that
			%\[
			%\sup_{x\in\Omega_{j,+}, y\in\Omega_{j,-}} |G(x-y)-G(p_+-p_-)| \to 0.
			%\]
			%Using $\|\omega_{j,+}\|_1 = \|\omega_{j,-}\|_1 = 1/2$, we obtain
			\[
			E_{\mathrm{int}}[\omega_{j,+},\omega_{j,-}] = \frac{1}{4} G(p_+-p_-) + o(1).
			\]
			This with \eqref{eq:renorm_self} proves \eqref{eq:renorm}.
			The last assertion follows from the Riesz rearrangement inequality for the kernel $\log_+\frac{1}{|z|}$ with the fact that $-\log |z| \le \log_+ \frac{1}{|z|}$ and $\operatorname{diam} B_* = \sqrt{2/\pi} < 1$.
		\end{proof}
		
		\begin{remark}\label{rem:ren_en_gen}
			The maximality of $\omega_j \in \calS_{\mu, \varepsilon_j}$ is not used in the proof of Lemma~\ref{lem:ren_en_decom_2}.
			In fact, the same conclusion holds for any sequence $\rho_j = \varepsilon_j^{-2} \left( \mathbf 1_{\calD_{j,+}}-\mathbf 1_{\calD_{j,-}} \right)$ satisfying
			\[
			|\calD_{j,\pm}|=\frac{\varepsilon_j^2}{2}, \quad \operatorname{diam}(\calD_{j,\pm})\le C\varepsilon_j,
			\]
			and corresponding points $p_{j,\pm} \in \calD_{j,\pm}$ such that $p_{j,\pm} \to p_\pm$ and $p_+ \neq p_-$.
		\end{remark}

		\begin{lemma}[Point-vortex maximization] \label{lem:max_pv}
			For $\mu > 0$, define
			\[
			\calP_\mu := \left\{ (p_+,p_-) \in \calC \times \calC :  \frac{|p_{+,2}|+|p_{-,2}|}{2} = \mu \right\},
			\]
			where $p_{\pm} = (p_{\pm, 1}, p_{\pm, 2})$.
			Then
			\[
			\max_{(p_+,p_-) \in \calP_\mu} \{-G(p_+-p_-)\} = -G(\pi,2\mu), \quad (-G(0) := -\infty).
			\]
			In particular, $(p_+,p_-) \in \calP_\mu$ is the maximum point if and only if $p_{+,1} - p_{-,1} = \pi$ in $\bbT$ and  $|p_{+,2}-p_{-,2}|=2\mu$.
		\end{lemma}
		%Equivalently, every maximizing pair has the form
		%\[
		%p_+ = (a, t + s\mu), \quad p_-=(a + \pi, t - s\mu),
		%\]
		%for some $a \in \bbT$,  $t\in[-\mu,\mu]$, and $s \in \{+1, -1 \}$.
		
		\begin{proof}
			For $d := p_+-p_- = (d_1,d_2)$, we have $|d_2| \le |p_{+,2}|+|p_{-,2}| =2\mu$.
			Note that $-4\pi G(d) = \log\left(2\cosh d_2-2\cos d_1\right)$, where the right-hand side is maximized only when $|d_2|=2\mu$ and $\cos d_1=-1$, i.e. $d_1 = \pi$ in $\bbT$.
			The lemma is proved.
		\end{proof}
		
		\begin{proposition}[Relative position]\label{prop:rela_pos}
			Let $\varepsilon_j \to 0$ and let $\omega_j = \varepsilon_j^{-2} \left( \mathbf 1_{\Omega_{j,+}} - \mathbf 1_{\Omega_{j,-}} \right) \in \calS_{\mu,\varepsilon_j}$.
			There exist $p_\pm \in \calC$ such that $\omega_{j,\pm} \rightharpoonup \frac{1}{2} \delta_{p_\pm}$, where $p_{\pm} = (p_{\pm, 1}, p_{\pm, 2})$ satisfies $p_{+,1}-p_{-,1} = \pi$ in $\bbT$ and $|p_{+,2} - p_{-,2}| = 2\mu$, up to a subsequence.
			Moreover, if $A_{j,\pm}$ are the rescaled patches \eqref{eq:res_patch} in Lemma~\ref{lem:ren_en_decom_2}, then $\calE(A_{j,\pm}) \to \calE(B_*)$, where $\calE(\cdot)$ and $B_*$ are defined as in Lemma~\ref{lem:ren_en_decom_2}.
		\end{proposition}

		\begin{proof}
			By Lemma~\ref{lem:asy_pt}, there exist $p_\pm\in\mathcal C$ such that $\omega_{j,\pm} \rightharpoonup \frac{1}{2}\delta_{p_\pm}$, up to a subsequence.
			We first identify the limiting constraint.
			Since $\|x_2\omega_j\|_1 = \mu$ and $\supp \omega_j \subset \{ |x_2| \le \bar{R}(\mu) \}$ by Lemma~\ref{lem:supp}, for a test function $\varphi \in C_c(\calC)$ such that $\varphi(x) = |x_2|$ in $\{ |x_2| \le \bar{R} \}$, we see
			\[
			\mu = \int_{\mathcal C} \varphi(x) \left(\omega_{j,+}+\omega_{j,-}\right)(x)\,dx \to \frac{ \varphi(p_+) + \varphi(p_-)}{2} = \frac{|p_{+,2}| + |p_{-,2}|}{2},
			\]
			Thus $(p_+,p_-) \in \calP_\mu$, where $\calP_{\mu}$ is defined in Lemma~\ref{lem:max_pv}.
			By Lemma~\ref{lem:ren_en_decom_2}, we also have
			\begin{equation}\label{eq:en_asy_2}
				E[\omega_j] = \frac{1}{8\pi}|\log\varepsilon_j| + \calE(A_{j,+}) + \calE(A_{j,-}) - \frac{1}{4} G(p_+-p_-) + o(1).
			\end{equation}
			We now construct a disk competitor.
			Let $q_+ = (0,\mu)$ and $q_- = (\pi,-\mu)$, and define
			\[
			\rho_j := \varepsilon_j^{-2} \left( \mathbf 1_{B_{\varepsilon_j/\sqrt{2\pi}}(q_+)} - \mathbf 1_{B_{\varepsilon_j/\sqrt{2\pi}}(q_-)} \right).
			\]
			For all sufficiently large $j$, we have $x_2\rho_j \ge 0$ in $\calC$.
			Each disk has mass $1/2$, and $\|x_2\rho_j\|_1 = \mu$.
			Thus $\rho_j\in\mathcal F_{\mu,\varepsilon_j}$.
			Moreover, since the sequence $\{ \rho_j \}$ satisfies the assumptions of Remark~\ref{rem:ren_en_gen} with $p_{j, \pm} = q_{\pm}$ and $\calD_{j, \pm} = B_{\varepsilon_j/\sqrt{2\pi}}(q_{\pm})$ in there, we find
			\begin{equation}\label{eq:en_disk}
				E[\rho_j] = \frac{1}{8\pi}|\log\varepsilon_j| + 2\calE(B_*) - \frac{1}{4}G(\pi,2\mu) + o(1),
			\end{equation}
			where $B_* = B_{1/\sqrt{2\pi}}(0)$ is the rescaled disk of $B_{\varepsilon_j/\sqrt{2\pi}}(q_{\pm})$.
			
			Since $\omega_j$ is a maximizer, $E[\omega_j] \ge E[\rho_j]$.
			Combining \eqref{eq:en_asy_2} and \eqref{eq:en_disk}, we obtain
			\begin{equation}\label{eq:comp_asy}
				\begin{aligned}
					0 &\le \left( -\frac{1}{4}G(\pi,2\mu) + \frac{1}{4}G(p_+-p_-) \right) \\
					&\le \left( 2\calE(B_*) - \calE(A_{j,+}) - \calE(A_{j,-}) \right) 
					+
					\left( -\frac{1}{4}G(\pi,2\mu) + \frac{1}{4}G(p_+-p_-) \right) \\
					&= E[\rho_j] - E[\omega_j] +  o(1) \le o(1),
				\end{aligned}
			\end{equation}
			where the first and second inequality is due to Lemma~\ref{lem:max_pv} and the last assertion in Lemma~\ref{lem:ren_en_decom_2}, respectively.
			
			Consequently, sending $j \to \infty$ gives $- G(p_+ - p_-) = - G(\pi, 2\mu)$, and hence the desired properties for $p_{\pm}$ follow from Lemma~\ref{lem:max_pv}.
			The last assertion $\calE(A_{j,\pm})\to\calE(B_*)$ also follows from \eqref{eq:comp_asy} with the fact that $\calE(B_*)-\calE(A_{j,\pm}) \ge 0$ obtained in Lemma~\ref{lem:ren_en_decom_2}.
		\end{proof}
		
		\begin{remark}\label{rmk:vert_order}
			In particular, there exists $\varepsilon_0 = \varepsilon_0(\mu) > 0$ such that the following holds: for any $\varepsilon \in (0, \varepsilon_0)$ and  $\omega = \varepsilon^{-2} \left(\mathbf 1_{\Omega_+}-\mathbf 1_{\Omega_-}\right) \in\calS_{\mu,\varepsilon}$, one core lies strictly above the other, with vertical gap at least $\mu$, that is,
			\[
			\operatorname*{ess\,inf}_{\Omega_+}x_2 - \operatorname*{ess\,sup}_{\Omega_-}x_2 \ge \mu
			\]
			or the same inequality with $\Omega_+$ and $\Omega_-$ interchanged.
			Indeed, Proposition~\ref{prop:rela_pos} shows that the limiting vertical separation is $2\mu$, while Lemma~\ref{lem:core} gives
			$\operatorname{diam}(\Omega_\pm) \lesssim_{\mu} \varepsilon$.
			We emphasize that this assertion is uniform over all maximizers in $\calS_{\mu, \varep}$.
		\end{remark}

		\begin{corollary}\label{cor:res_core}
			Under the assumptions in Proposition~\ref{prop:rela_pos}, there exist bounded sequences $\{ a_{j, \pm} \}_j$ such that  $|(A_{j, \pm} - a_{j, \pm}) \triangle B_*| \to 0$.
		\end{corollary}
		
		\begin{proof}
			For simplicity, we write $A_j := A_{j, \pm}$.
			Note that, since $A_j$ is a rescaled patch, we may consider $\bbR^2$ rather than $\calC$.
			Suppose, to the contrary, that for some $\delta > 0$, up to a subsequence,
			\[
			\inf_{a \in \bbR^2} |(A_{j} - a) \triangle B_*| \ge \delta.
			\]
			Since $0 \in A_j$ and $\operatorname{diam} A_{j, \pm} \le C(\mu)$ by Lemma~\ref{lem:ren_en_decom_2}, $A_j$ are contained in a fixed ball $B$.
			Thus, up to a subsequence, $\mathbf{1}_{A_j} \rightharpoonup^* f$ in $L^{\infty}(\bbR^2)$, where $0 \le f \le 1$, $\supp f \subset B$, and $\|f\|_1 = 1/2$.
			Since $\log|x-y| \in L^1(B \times B)$, by the weak-$^*$ convergence and Proposition~\ref{prop:rela_pos} implies
			\[
			\calE(A_j) \to -\frac{1}{4\pi} \iint_{\bbR^2 \times \bbR^2} \log |x-y| f(x) f(y) \,dxdy = \calE(B_*).
			\]
			The equality case in the Riesz rearrangement inequality implies that $f(x + a) = \mathbf{1}_{B_*}$ for some $a \in \bbR^2$.
			Therefore, $|(A_j - a)\triangle B_*| \to 0$ which is a contradiction.
			Finally, the existence of bounded sequence $\{a_j\}_j$ satisfying the assertion follows directly.
		\end{proof}

		\subsection{Asymptotic localization to von K\'arm\'an vortex pair}\label{subsec:local}
		
		\begin{definition}
			For $(a,b) \in \calC$, define the translation operator
			\[
			(\tau_{a,b}\omega) (x_1,x_2) := \omega(x_1-a,x_2-b).
			\]
			We also define the vertical reflection
			\[
			(\mathcal R \omega)(x_1,x_2) := \omega(x_1,-x_2).
			\]
			Note that  $E[\tau_{a,b}\omega]=E[\omega]$, $E[\mathcal R \omega] = E[\omega]$, $I(\tau_{a,b}\omega) = I(\omega)$, and $I(\mathcal R\omega) = - I(\omega)$.
			
			For $\omega = \varepsilon^{-2}(\mathbf{1}_{\Omega_+} - \mathbf{1}_{\Omega_-}) \in \calS_{\mu,\varepsilon}$, the vertical centers of the two patches are defined by
			\[
			z_\pm (\omega) := 2\int_{\calC}x_2\omega_\pm(x)\,dx.
			\]
			We also define the common vertical center by
			\[
			c_2(\omega) := \int_{\calC}x_2|\omega(x)|\,dx = \frac{z_+(\omega)+z_-(\omega)}{2}.
			\]
			The impulse satisfies $2I(\omega) = z_+(\omega) - z_-(\omega)$.
			Moreover, $z_{\pm}(\omega) = c_2(\omega) \pm I(\omega)$, $c_2(\tau_{a,b}\omega) = c_2(\omega)+b$, and $c_2(\mathcal R\omega) = - c_2(\omega)$.
		\end{definition}
		
		\begin{proposition}[Impulse saturation]\label{prop:orientation}
			There exists $\varepsilon_0 = \varepsilon_0(\mu) > 0$ such that, for any $\varepsilon \in (0,\varepsilon_0)$ and $\omega\in\calS_{\mu,\varepsilon}$, we have $|I(\omega)|=\mu$.
			In particular, one of the following alternatives holds: $x_2 \omega(x) \ge 0$ or $x_2\omega(x) \le 0$ for a.e. $x \in \calC$.
		\end{proposition}
		
		\begin{proof}
			Applying the vertical reflection $\mathcal R$ if necessary, Remark~\ref{rmk:vert_order} implies that we may assume $\operatorname*{ess\,inf}_{x\in\Omega_+}x_2 - \operatorname*{ess\,sup}_{y\in\Omega_-}y_2 \ge \mu$.
			In particular, since $\|\omega_+\|_1 = \|\omega_-\|_1 = 1/2$, we have $I(\omega) \ge \mu/2 > 0$.
			Set $\bar\omega := \tau_{0,-c_2(\omega)}\omega$.
			Then $z_+(\bar\omega) = I(\omega)$ and  $z_-(\bar\omega)=-I(\omega)$.
			Also, for all small $\varepsilon$, we have $x_2 \bar{\omega} \ge 0$ in $\calC$, and therefore $\|x_2\bar\omega\|_1 = I(\omega) \le \mu$.
			Since the energy is translation invariant, $\bar\omega$ is also a maximizer in $\calF_{\mu,\varepsilon}$.
			
			Suppose that $I(\omega) < \mu$.
			Choose $0 < h < \mu - I(\omega)$ and define 
			\[
			\bar\omega_h := \tau_{0,h}\bar\omega_+ - \tau_{0,-h}\bar\omega_-.
			\]
			Clearly, $\|\bar\omega_h\|_1 = 1$,  $\|\bar\omega_h\|_{\infty} = \varepsilon^{-2}$, $\int \bar\omega_h = 0$, and $\|x_2\bar\omega_h\|_1 = I(\omega) + h < \mu$.
			Thus $\bar\omega_h \in \calF_{\mu,\varepsilon}$.
			Note that $E[\bar{\omega}_{\pm}] = E[\tau_{0, \pm h}\bar{\omega}_{\pm}]$.
			On the other hand,
			for $x\in\supp\bar\omega_+$ and $y\in\supp\bar\omega_-$, we have
			$x_2-y_2>0$, and
			\[
			\partial_{x_2}G(x_1,x_2) = -\frac{1}{4\pi} \frac{\sinh x_2}{\cosh x_2-\cos x_1} < 0 \quad \text{ for } \, x_2 > 0.
			\]
			Hence $E_{\rm int} \left[ \tau_{0,h}\bar\omega_+, \, \tau_{0,-h}\bar\omega_- \right] < E_{\rm int}[\bar\omega_+, \bar\omega_-]$, and then $E[\bar\omega_h] > E[\bar\omega]$, contradicting maximality.
			Therefore $I(\omega) = \mu$.
			Finally, since $\|x_2\omega\|_1=\mu$, we have $0 = \|x_2\omega\|_1 - I(\omega) = \int \left( |x_2\omega| - x_2\omega \right)$.
			Since the integrand is nonnegative, we have $x_2 \omega \ge 0$.
		\end{proof}

		\begin{definition}[Centered maximizers] \label{def:cen_max}
			For sufficiently small $\varepsilon$, we define the class of centerized maximizers as
			\[
			\calS^{\rm cen}_{\mu,\varepsilon} := \left\{ \omega \in \calS_{\mu,\varepsilon} : I(\omega) = \|x_2\omega\|_1 = \mu, \,\, c_2(\omega)=0 \right\} \subsetneq \calS_{\mu, \varepsilon}.
			\]
			Equivalently,
			\[
			\calS^{\rm cen}_{\mu,\varepsilon} = \left\{ \omega \in \calS_{\mu,\varepsilon}: x_2\omega \ge 0, \,\,
			\int_{\calC}x_2|\omega(x)|\,dx = 0 \right\}.
			\]
		\end{definition}
		
		\begin{remark}\label{rmk:cen_max}
			For all sufficiently small $\varepsilon > 0$, $\calS^{\rm cen}_{\mu,\varepsilon} \neq \emptyset$.
			In fact, any $\omega \in \calS_{\mu,\varepsilon}$ can be mapped into $\calS^{\rm cen}_{\mu,\varepsilon}$ by applying the vertical reflection $\calR$ (if $I(\omega) = -\mu$) followed by the vertical translation $\tau_{0, -c_2(\cdot)}$.
			Indeed, by Proposition~\ref{prop:orientation}, we have $I(\omega) \in \{\mu, -\mu\}$ and the sign of $ x_2\omega$ is strictly determined.
			Since $I(\calR \omega) = -I(\omega)$ and vertical translations shift the center of mass $c_2(\omega)$ while preserving $I(\omega)$, we see
			\[
			\omega  \longmapsto   \tld{\omega} := \begin{cases}
				\omega  &\textrm{if }\,\, I(\omega) = \mu \\
				\calR\omega  &\textrm{if }\,\, I(\omega) = -\mu
			\end{cases}
			\longmapsto  \tau_{0, -c_2(\tilde{\omega})} \tld{\omega} \in \calS^{\rm cen}_{\mu,\varepsilon}.
			\]
		\end{remark}
		
		\begin{proposition}[Properties of $\calS^{\rm cen}_{\mu,\varepsilon}$] \label{prop:cen_prop}
			For any sequence $\omega_j \in \calS^{\rm cen}_{\mu,\varepsilon_j}$ with $\varepsilon_j \to 0$, the following hold:
			\begin{enumerate}
				\item \textit{Localization}: The vertical centers satisfy $z_+(\omega_j) = \mu$ and $z_-(\omega_j) = -\mu$.
				Consequently, 
				\begin{equation}\label{eq:cen_ver}
					\Omega_{j,\pm} \subset \bbT \times [\pm \mu - C\varepsilon_j, \pm \mu + C\varepsilon_j]
				\end{equation}
				for some $C = C(\mu) > 0$.
				
				\item \textit{Asymptotic position}: Up to a subsequence and horizontal translations,
				\[
				\omega_{j,+} \rightharpoonup \frac{1}{2}\delta_{(0,\mu)}, \quad \omega_{j,-} \rightharpoonup \frac{1}{2}\delta_{(\pi,-\mu)}.
				\]
				
				\item \textit{Energy expansion}: 
				\begin{equation}\label{eq:cen_exp}
					E[\omega_j] = \frac{1}{8\pi}|\log\varepsilon_j|  - \frac{1}{2\pi} \iint_{B_* \times B_*} \log |x-y|\,dxdy - \frac{1}{4}G(\pi,2\mu) + o(1),
				\end{equation}
				where $B_* = B_{1/\sqrt{2\pi}}(0)$.
			\end{enumerate}
		\end{proposition}
		
		\begin{proof}
			\textit{(1)}
			For a sequence $\omega_j \in \calS^{\rm cen}_{\mu,\varepsilon_j}$, $z_{\pm}(\omega_j) = c_2(\omega_j) \pm I(\omega_j) = \pm \mu$ since $I(\omega_j) = \mu$ and $c_2(\omega_j) = 0$.
			Then \eqref{eq:cen_ver} directly follows from this with Lemma~\ref{lem:core}.
			
			\smallskip
			
			\noindent \textit{(2)}
			Proposition~\ref{prop:rela_pos} with \eqref{eq:cen_ver} imply, up to a subsequence, $\omega_{j,\pm} \rightharpoonup \frac{1}{2}\delta_{p_\pm}$ such that $p_{\pm, 2} = \pm \mu$ and $p_{+,1} - p_{-,1} = \pi$.
			A suitable horizontal translation maps these limits to $(0,\mu)$ and $(\pi,-\mu)$.
			
			\smallskip
			
			\noindent \textit{(3)}
			Finally, we substitute $p_+ - p_- = (\pi, 2\mu)$ into \eqref{eq:renorm} of Lemma~\ref{lem:ren_en_decom_2}.
			The convergence $\calE(A_{j,\pm}) \to \calE(B_*)$ in Proposition~\ref{prop:rela_pos} then yields \eqref{eq:cen_exp}.
		\end{proof}
		
		Proposition~\ref{prop:cen_prop} allows us to rewrite the Euler--Lagrange equation \eqref{eq:EL} for $\omega  \in \calS^{\rm cen}_{\mu,\varepsilon}$ when $\varepsilon$ is sufficiently small:
		\begin{equation}\label{eq:EL2}
			\omega = \varepsilon^{-2} \left(\mathbf{1}_{ \{ \psi[\omega] - Wx_2 - \gamma_+ > 0 \} \cap \{ x_2 > \mu/2 \} } - \mathbf{1}_{ \{ -\psi[\omega] + Wx_2 - \gamma_- > 0 \} \cap \{ x_2 < -\mu/2 \} } \right)
		\end{equation}
		for uniquely determined constants $\gamma_{\pm} = \gamma_{\pm}(\omega)$ and $W = W(\omega)$ such that $\gamma_+ + \gamma_-, W \ge 0$.
		In particular, since $\Omega_+ \Subset \{ x_2 > 0 \}$ and $\Omega_- \Subset \{ x_2 < 0 \}$,
		\begin{equation}\label{eq:EL3}
			\psi[\omega] - Wx_2 =  \pm \gamma_{\pm}  \quad\text{on }\,\, \partial\Omega_{\pm}.
		\end{equation}
		Thus the velocity $u - We_1 = \nabla^{\perp}\psi - W e_1$ is tangent to both patch boundaries, and hence $\omega(x_1 - Wt, x_2)$ is a \emph{traveling-wave} solution.
		In fact, Lemma~\ref{lem:local_poho} below shows that $W > 0$.

		\subsection{Asymptotic translation speed}\label{subsec:speed}
		
		\begin{lemma} \label{lem:ext_field}
			Let $\varepsilon_j \to 0$ and $\omega_j \in \calS^{\rm cen}_{\mu,\varepsilon_j}$.
			Suppose that $\omega_{j,\pm} \rightharpoonup \frac{1}{2}\delta_{p_\pm}$ where $p_+ = (0,\mu)$ and $p_- =  (\pi,-\mu)$.
			Choose $r \in (0, |p_+ - p_-|/3)$.
			Then, for every integer $k \ge 0$,
			\[
			\psi[\omega_{j,-}] \to \frac{1}{2}G(\cdot - p_-) \quad \text{in}\,\, C^k(\overline{B_r(p_+)}) \quad \textrm{and} \quad \psi[\omega_{j,+}] \to \frac{1}{2}G( \cdot - p_+) \quad \text{in}\,\,C^k(\overline{B_r(p_-)}).
			\]
		\end{lemma}
		
		\begin{proof}
			We prove the first convergence; the second is identical.
			By Lemmas~\ref{lem:core} and~\ref{lem:asy_pt}, $\Omega_{j,\pm}\subset B_r(p_{\pm})$ for all sufficiently large $j$.
			Hence $G(x-y)$ is smooth on $\overline{B_r(p_+)} \times \overline{B_r(p_-)}$.
			Let $\alpha$ be a multi-index with $|\alpha|\le k$.
			Since $\|\omega_{j,-}\|_1=1/2$, for $x \in \overline{B_r(p_+)}$ we have
			\begin{align*}
				\partial_x^\alpha \psi[\omega_{j,-}](x)  -  \frac{1}{2} \partial^\alpha G(x-p_-)
				= \int_{\calC} \left( \partial^\alpha G(x-y) - \partial^\alpha G(x-p_-) \right) \omega_{j,-}(y)\,dy.
			\end{align*}
			Therefore,
			\begin{align*}
				\sup_{x\in\overline{B_r(p_+)}} \left| \partial_x^\alpha\psi[\omega_{j,-}](x) - \frac{1}{2} \partial^\alpha G(x-p_-) \right|
				\le \frac{1}{2} \sup_{\substack{x\in\overline{B_r(p_+)}\\ y\in\Omega_{j,-}}} \left| \partial^\alpha G(x-y) - \partial^\alpha G(x-p_-) \right|.
			\end{align*}
			Since $\sup_{y\in\Omega_{j,-}} |y - p_-| \to 0$ by Lemmas~\ref{lem:core} and~\ref{lem:asy_pt}, the right-hand side tends to zero.
			This proves the $C^k$ convergence.
		\end{proof}

		\begin{lemma}[Traveling speed] \label{lem:local_poho}
			Let $\omega = \omega_+ - \omega_- \in\calS^{\rm cen}_{\mu,\varepsilon}$, and let $W$ be the multiplier in \eqref{eq:EL3}.
			Then
			\begin{equation}\label{eq:260801}
				\int_{\calC} \omega_+(x) \partial_2 \psi[\omega_-](x)\,dx = -\frac{W}{2} \quad \textrm{and} \quad \int_{\calC} \omega_-(x) \partial_2 \psi[\omega_+](x)\,dx  =  \frac{W}{2} .
			\end{equation}
			Furthermore, every $\omega \in \calS^{\rm cen}_{\mu,\varepsilon}$ is a traveling wave solution with a traveling speed $W = W(\omega)$ determined by
			\begin{equation} \label{eq:W_formula}
				W = \frac{1}{2\pi} \iint_{\calC \times \calC} \frac{\sinh(x_2 - y_2)}{\cosh(x_2 - y_2) - \cos(x_1 - y_1)} \omega_+(x)\omega_-(y)\,dxdy > 0.
			\end{equation}
		\end{lemma}
		
		\begin{proof}
			Let $e_2 = (0, 1)$ and set $\Psi_+ :=  \left( \psi[\omega]-W x_2-\gamma_+ \right)_+ \mathbf{1}_{x_2 \ge 0}$.
			By \eqref{eq:EL2}, we see
			\[
			\nabla \Psi_+ = \mathbf 1_{\Omega_+} \left( \nabla\psi[\omega]-We_2 \right).
			\]
			Since $\Psi_+$ is compactly supported (Lemma~\ref{lem:core}),
			\[
			\int_{\Omega_+}\nabla\Psi_+\,dx =\int_{\calC}\nabla\Psi_+\,dx =  0.
			\]
			Multiplying by $\varepsilon^{-2}$ gives
			\[
			\int_{\calC} \omega_+(x)\nabla  \psi[\omega_-](x) \,dx = - \int_{\calC}\omega_+(x) \nabla\psi[\omega](x)\,dx = -\frac{W}{2} e_2,
			\]
			where the first equality is due to the oddness of $\nabla G$.
			This proves the first identity in \eqref{eq:260801}.
			The second identity in \eqref{eq:260801} is obtained from $\Psi_- := \left(-\psi[\omega] + Wx_2-\gamma_-\right)_+ \mathbf{1}_{x_2 \le 0}$ by the same way.
			
			To derive \eqref{eq:W_formula}, we substitute $\partial_2 \psi[\omega_-](x) = \int_{\calC} \partial_2 G(x-y)\omega_-(y)\,dy$ into the first identity.
			For $x \in \supp \omega_+$ and $y \in \supp \omega_-$, we see that
			\[
			\partial_2 G(x-y) = -\frac{1}{4\pi} \frac{\sinh(x_2-y_2)}{\cosh(x_2-y_2) - \cos(x_1-y_1)} <0
			\]
			by Proposition~\ref{prop:cen_prop}, and then we obtain \eqref{eq:W_formula}.
			%\begin{align*}
			%-\frac{W}{2} &= \iint_{\calC \times \calC} \partial_2 G(x-y) \omega_+(x) \omega_-(y) \,dx dy \\
			%&= -\frac{1}{4\pi} \iint_{\calC \times \calC} \frac{\sinh(x_2 - y_2)}{\cosh(x_2 - y_2) - \cos(x_1 - y_1)} \omega_+(x)\omega_-(y)\,dxdy.
			%\end{align*}
			%This yields \eqref{eq:W_formula}.
		\end{proof}
		
		\begin{proposition}[Convergence of the traveling speed] \label{prop:W_limit}
			Let $\omega \in \calS_{\mu, \varepsilon}^{\rm cen}$ and $W = W(\omega)$ be the translation speed given by \eqref{eq:W_formula}.
			As $\varepsilon\to0$, we have
			\begin{equation}\label{eq:W_conv}
				\sup_{\omega\in\calS^{\rm cen}_{\mu,\varepsilon}} \left| W - \frac{1}{8\pi} \tanh \mu \right| \to 0.
			\end{equation}
		\end{proposition}
		
		\begin{proof}
			Let $\varepsilon_j \to 0$ and $\omega_j \in \calS^{\rm cen}_{\mu,\varepsilon_j}$.
			Since $W$ in \eqref{eq:W_formula} is invariant under horizontal translations, we may assume, by passing a subsequence and applying horizontal translations (Proposition~\ref{prop:cen_prop}), that $\omega_{j, \pm} \rightharpoonup \frac{1}{2}\delta_{p_{\pm}}$ where $p_+ = (0, \mu)$ and $p_- = (\pi, -\mu)$.
			By Lemmas~\ref{lem:ext_field} and \ref{lem:local_poho},
			\[
			-\frac{W_j}{2}  = \int_{\calC} \omega_{j,+}\partial_2 \psi[\omega_{j,-}]\,dx \,\, \to \,\, \frac{1}{4} \partial_2 G(\pi,2\mu) = -\frac{1}{16\pi}\tanh\mu.
			\]
			Since this holds for arbitrary sequences $\{\varepsilon_j\}$, we obtain \eqref{eq:W_conv}.
		\end{proof}
		
		\subsection{Simply connected vortex cores}\label{subsec:simp-conn}
		
		We clarify the recentering notation used below.
		The points $p_{j,\pm} \in \Omega_{j,\pm}$ are arbitrary reference points used to define the rescaled patch:
		\[
		A_{j,\pm} = \frac{\Omega_{j,\pm} - p_{j,\pm}}{\varepsilon_j}.
		\]
		The bounded sequences $a_{j,+}$ and $a_{j, -}$ in Corollary~\ref{cor:res_core} recenter the two rescaled patches separately so that
		\[
		\left| (A_{j,\pm} - a_{j,\pm}) \triangle B_* \right| \to 0.
		\]
		The corresponding recentering points in the original variables are
		\[
		q_{j,\pm} := p_{j,\pm} + \varepsilon_j a_{j,\pm},
		\]
		and hence
		\[
		A_{j,\pm} - a_{j,\pm} = \frac{\Omega_{j,\pm} - q_{j,\pm}}{\varepsilon_j}.
		\]
		Thus the two cores are viewed in separate rescaled coordinate systems centered at $q_{j,+}$ and $q_{j,-}$, respectively, and the points $q_{j,\pm}$ need not belong to $\Omega_{j,\pm}$.
		By \eqref{eq:cen_ver}, the vertical components of $q_{j, \pm}$ satisfy
		\[
		q_{j,+,2} = \mu + O(\varepsilon_j), \quad q_{j,-,2} = -\mu + O(\varepsilon_j).
		\]
		Moreover, after a common horizontal translation and passing to a subsequence,
		\[
		q_{j,+} \to (0,\mu), \quad q_{j,-} \to (\pi,-\mu).
		\]
		Here and below, balls in the rescaled variables without a specified center are centered at the origin.
		
		To prove that the cores $\Omega_{\pm}$ are simply connected, we use Turkington's classical rescaling argument for concentrated vortices in bounded planar domains \cite{Tur1983}; see also \cite{CWZ2020}.
		
		\begin{lemma}\label{lem:res_str}
			Let $\varepsilon_j \to 0$ and $\omega_j \in \calS^{\rm cen}_{\mu,\varepsilon_j}$.
			Let $p_{j,\pm}$ and $A_{j,\pm}$ be as in Lemma~\ref{lem:ren_en_decom_2}, and let $\{ a_{j,\pm} \}_j$ be as in Corollary~\ref{cor:res_core}.
			Set $q_{j, \pm} = p_{j, \pm} + \varepsilon_j a_{j, \pm}$ and $D_{j, \pm} := A_{j, \pm} - a_{j, \pm}$.
			%\begin{align*}
			%    q_{j, \pm} &:= p_{j, \pm} + \varepsilon_j a_{j, \pm}, \\
			%    D_{j, \pm} &:= A_{j, \pm} - a_{j, \pm} = \frac{\Omega_{j, \pm} - q_{j, \pm}}{\varepsilon_j}.
			%\end{align*}
			Let $W_j := W(\omega_j)$ and $\gamma_{j,\pm} := \gamma_\pm(\omega_j)$, and define
			\begin{align*}
				\Psi_{j,+}(x) &:= \psi[\omega_j](q_{j,+} + \varepsilon_j x)
				- W_j(q_{j,+,2} + \varepsilon_j x_2) - \gamma_{j,+}, \\
				\Psi_{j,-}(x) &:= -\psi[\omega_j](q_{j,-} + \varepsilon_j x)
				+ W_j(q_{j,-,2} + \varepsilon_j x_2) - \gamma_{j,-}.
			\end{align*}
			Then $\Psi_{j,\pm} \to \Psi_*$ in $C^1_{\rm loc}(\bbR^2)$, where
			\[
			\Psi_*(x) := -\frac{1}{2\pi}\int_{B_*}\log |x-y|\,dy + \frac{1}{4\pi}\log\frac{1}{\sqrt{2\pi}}.
			\]
		\end{lemma}
		
		\begin{proof}
			For simplicity, we prove the assertion for $\Psi_j := \Psi_{j,+}$.
			Write $D_j := D_{j, +}$, $q_j := q_{j, +}$, and define
			\[
			U_j(x) := -\frac{1}{2\pi}\int_{D_j}\log |x-y|\,dy, \quad U_*(x) := -\frac{1}{2\pi}\int_{B_*}\log |x-y|\,dy.
			\]
			By Corollary~\ref{cor:res_core}, $\mathbf{1}_{D_j} \to \mathbf{1}_{B_*} = \mathbf{1}_{B_{1/\sqrt{2\pi}}(0)}$ in $L^1(\bbR^2)$.
			Moreover, we may assume that $D_j \subset  B_r$ for some $r > 1/\sqrt{2\pi}$.
			Since $- \Delta(U_j-U_*) = \mathbf{1}_{D_j} - \mathbf{1}_{B_*}$, by elliptic regularity and Sobolev embedding, we have $U_j \to U_*$ in $C^1_{\rm loc}(\bbR^2)$.
			Setting $R_*:=1/\sqrt{2\pi}$, a direct computation gives
			\begin{equation}\label{eq:0807_02}
				U_*(x) =U_*(|x|) = 
				\begin{cases}
					-\frac{1}{4\pi}\log R_* + \frac{R_*^2-|x|^2}{4}, & |x|\le R_*,  \\ 
					-\frac{1}{4\pi}\log |x|, & |x|\ge R_*.
				\end{cases}
			\end{equation}
			In particular, $U_*$ is strictly decreasing in $|x|$ and every level set of $U_*$ has measure zero.
			
			We now write $G(z) = -\frac{1}{2\pi}\log |z| + R(z)$ as in \eqref{eq:0806_01} and recall that $R$ is smooth near the origin and $R(0)=0$.
			By a change of variable, we see
			\[
			\psi[\omega_{j,+}](q_j + \varepsilon_j x) = \frac{1}{4\pi}|\log\varepsilon_j| + U_j(x) + \underbrace{\int_{D_j} R\left(\varepsilon_j(x-y)\right)\,dy}_{=: V_j(x)}.
			\]
			Since $D_j \subset B_r$, it is clear that $V_j \to 0$ in $C^1_{\rm loc}(\bbR^2)$.
			We then decompose $\Psi_{j}$ by
			\begin{equation}\label{eq:res_str_decomp}
				\Psi_{j}(x) = U_j(x) - c_j + r_j(x),
			\end{equation}
			where
			\[
			c_j := \gamma_{j,+} - \frac{1}{4\pi}|\log\varepsilon_j| + \psi[\omega_{j,-}](q_j) + W_jq_{j,2}
			\]
			and
			\[
			r_j(x) := V_j(x) - \left\{ \psi[\omega_{j,-}](q_j+\varepsilon_jx) - \psi[\omega_{j,-}](q_j) \right\} - \varepsilon_j W_j x_2.
			\]
			By \eqref{eq:cen_ver}, the two cores are separated by a fixed positive distance.
			Hence, for every $R>0$,
			\[
			\left\| \psi[\omega_{j,-}](q_j + \varepsilon_j \cdot) - \psi[\omega_{j,-}](q_j) \right\|_{C^1(B_R)} \le C(R,\mu) \varepsilon_j, \quad q_j = q_{j, +}.
			\]
			Moreover, $\{W_j\}_j$ is bounded by Proposition~\ref{prop:W_limit}.
			Consequently, $r_j \to 0$ in $C^1_{\rm loc}(\bbR^2)$.
			
			We may assume that $D_j \cup B_* \subset B_r$.
			Since $q_{j,2} = \mu + O(\varepsilon_j)$, we also have $q_j + \varepsilon_j B_r = B_{\varepsilon_j r}(q_j) \subset \{x_2 > \mu/2\}$ for all sufficiently large $j$.
			Thus the Euler--Lagrange equation \eqref{eq:EL2} and \eqref{eq:res_str_decomp} imply
			\begin{equation}\label{eq:0807_01}
				D_j = \left\{ x \in B_r: \Psi_j > 0 \right\} = \left\{ x \in B_r: U_j(x) + r_j(x) > c_j \right\}.
			\end{equation}
			Then we see $\inf_{B_r}(U_j + r_j) \le c_j \le \sup_{B_r}(U_j + r_j)$, and therefore, $\{c_j\}_j$ is bounded.
			Suppose that $c_j \to c$ along a subsequence.
			Since $U_j + r_j \to U_*$ uniformly on $B_r$, we have $\mathbf{1}_{D_j} \to \mathbf{1}_{\{ U_* > c \}}$ in $L^1(B_r)$.
			Together with $\mathbf{1}_{D_j}\to\mathbf{1}_{B_*}$ in $L^1(\bbR^2)$, this yields  $\{ x \in B_r: U_* > c \} = B_*$.
			Therefore, $c = U_* \big|_{\partial B_*} = -\frac{1}{4\pi}\log\frac{1}{\sqrt{2\pi}}$.
			Then it follows from \eqref{eq:res_str_decomp} that
			\[
			\Psi_{j,+} \to U_* + \frac{1}{4\pi}\log\frac{1}{\sqrt{2\pi}} = \Psi_* \quad \textrm{in }\,\, C^1_{\rm loc}(\bbR^2).
			\]
			The proof for $\Psi_{j,-}$ is identical.
		\end{proof}
		
		We conclude this section by showing that the vortex cores are simply connected for sufficiently small $\varepsilon$.
		
		\begin{proposition}[Simply connected cores]\label{prop:connect}
			There exists $\varepsilon_0 = \varepsilon_0(\mu) > 0$ such that, for any $\varepsilon \in (0, \varepsilon_0)$ and $\omega = \varepsilon^{-2}(\mathbf{1}_{\Omega_+} - \mathbf{1}_{\Omega_-}) \in \calS_{\mu, \varepsilon}^{\rm cen}$, the cores $\Omega_{\pm}$ (represented by open sets) are simply connected.
		\end{proposition}
		
		\begin{proof}
			Suppose otherwise.
			Then there exist $\varepsilon_j \to 0$ and $\omega_j \in \calS_{\mu,\varepsilon_j}^{\rm cen}$ such that one of the cores is not simply connected.
			Let $D_{j,\pm} = A_{j,\pm} - a_{j,\pm}$ and $\Psi_{j,\pm}$ be as in Lemma~\ref{lem:res_str}.
			Recall that $D_{j,\pm} \subset B_r$ for some $r > 1/\sqrt{2\pi}$.
			Set $R_* := 1/\sqrt{2\pi}$.
			By \eqref{eq:0807_02}, the limiting function $\Psi_*$ in Lemma~\ref{lem:res_str} is given by
			\[
			\Psi_*(x) =
			\begin{cases}
				\frac{R_*^2 - |x|^2}{4}, & |x| \le R_*, \\[5pt]
				\frac{1}{4\pi}\log \frac{R_*}{|x|}, & |x| \ge R_*.
			\end{cases}
			\]
			In particular, for any $R_0 \in (0, R_*)$, $\Psi_* > 0$ in $\overline{B_{R_0}}$, while
			\[
			\frac{x}{|x|} \cdot \nabla\Psi_*(x) = 
			\begin{cases}
				-\frac{|x|}{2}, & 0 < |x| \le R_*, \\[5pt]
				-\frac{1}{4\pi |x|}, & |x| \ge R_*.
			\end{cases}
			\]
			Hence Lemma~\ref{lem:res_str} implies that, for all sufficiently large $j$, $\Psi_{j,\pm} > 0$ in $\overline{B_{R_0}}$ and
			\[
			\frac{x}{|x|} \cdot \nabla\Psi_{j,\pm}(x) < 0 \quad \textrm{for }\,\, |x| \ge R_0.
			\]
			Recall from \eqref{eq:0807_01} that $D_{j,\pm} = \left\{ x \in B_r : \Psi_{j,\pm}(x) > 0 \right\}$.
			Let $x \in D_{j,\pm}$ and $t \in [0,1]$.
			If $|tx| \le R_0$, then $\Psi_{j,\pm}(tx) > 0$.
			If $|tx| > R_0$, then
			\[
			\Psi_{j,\pm}(tx) - \Psi_{j,\pm}(x) = -\int_t^1 x \cdot \nabla\Psi_{j,\pm}(sx)\,ds \ge 0.
			\]
			Thus, $\Psi_{j,\pm}(tx) > 0$ and therefore $tx \in D_{j,\pm}$ for every $t \in [0,1]$.
			Hence $D_{j,\pm}$ are simply connected.
			The same holds for $\Omega_{j,\pm}$, which is a contradiction.
		\end{proof}

		\section{Compactness and Stability}\label{sec:cpt_stab}
		Throughout this section, we fix $\varepsilon_0 = \varepsilon_0(\mu) > 0$ sufficiently small so that all the results in \S \ref{sec:vari}--\ref{sec:karman} hold, possibly decreasing it further whenever necessary, and consider $\varepsilon \in (0, \varepsilon_0)$.
		
		\begin{definition}\label{def:orb}
			We define the positive-impulse orbit by
			\[
			\calO^+_{\mu,\varepsilon} := \left\{ \tau_{a,b}\omega: \omega \in \calS^{\rm cen}_{\mu,\varepsilon}, \quad (a,b)\in\bbT\times\bbR \right\},
			\]
			where $\calS^{\rm cen}_{\mu,\varepsilon}$ is defined in Definition~\ref{def:cen_max}.
			The full orbit is defined by $\calO_{\mu,\varepsilon} := \calO^+_{\mu,\varepsilon} \cup \calR \calO^+_{\mu,\varepsilon}$.
		\end{definition}
		
		\subsection{Compactness in $L^1$}
		
		\begin{lemma}\label{lem:comple}
			Let $\{h_n\}$ be a sequence of nonnegative measurable functions satisfying
			\[
			\limsup_{n\to\infty}\|h_n\|_\infty\le a, \quad h_n\rightharpoonup^* a\mathbf 1_A \quad \text{in} \,\, L^\infty(\calC),
			\]
			where $|A|<\infty$.
			If $\|h_n\|_1 \to a|A|$, then $h_n\to a\mathbf 1_A$ in $L^1(\calC)$.
		\end{lemma}
		
		\begin{proof}
			Since $\mathbf 1_A\in L^1(\calC)$, by the weak-$^*$ convergence we have $\int_A h_n \to a|A|$.
			\begin{equation}\label{eq:comple}
				\int_{A^c} |h_n - a \mathbf{1}_A|\,dx = \int_{A^c}h_n\,dx = \| h_n\|_1 - \int_Ah_n\,dx \to 0.
			\end{equation}
			It remains to show the convergence on $A$.
			Fix $\delta > 0$.
			By the assumption, for all large $n$, we have $0\le h_n \le a + \delta$.
			Then we see $\int_A (h_n - a)_+ \le \delta |A|$ and
			\[
			\int_A (a-h_n)_+\,dx = \int_A (h_n-a)_+\,dx + a|A| - \int_A h_n\,dx \le \delta |A| + a|A| - \int_A h_n\,dx
			\]
			This gives $\limsup_{n\to\infty} \int_A (a-h_n)_+ dx \le \delta |A|$.
			It follows that
			\[
			\limsup_{n\to\infty} \int_A |h_n - a|\,dx \le \limsup_{n\to\infty} \int_A (h_n - a)_+\,dx + \limsup_{n\to\infty} \int_A (a - h_n)_+ \,dx  \le 2\delta |A|.
			\]
			Letting $\delta \to 0$, we obtain $\| \mathbf{1}_A (h_n - a) \|_1 \to 0$.
			Together with \eqref{eq:comple}, this proves the lemma.
		\end{proof}

		\begin{proposition} \label{prop:seq_cpt}
			Let $M_{\mu,\varepsilon} := \sup_{\omega \in \calF_{\mu,\varepsilon}} E[\omega]$ and $\{\omega_n\}_n \subset L^1\cap L^{\infty}(\calC)$ be a sequence satisfying
			\begin{equation}
				\label{eq:seq_cpt}
				\int_{\calC}\omega_n\,dx=0, \quad    \limsup_{n\to\infty}\|\omega_n\|_1\le 1, \quad \limsup_{n\to\infty}\|\omega_n\|_\infty\le \varepsilon^{-2}.
			\end{equation}
			Assume that there is a sequence $\{b_n\} \subset \bbR$ such that $\limsup_{n\to\infty} \|x_2 (\tau_{0,b_n}\omega_n) \|_1\le \mu$.
			If $\liminf_{n\to\infty} E[\omega_n] \ge M_{\mu,\varepsilon}$, then
			\begin{equation}\label{eq:cpt_1}
				\inf_{\omega\in \calO_{\mu,\varepsilon}}\left\|\omega_n-\omega\right\|_1\to0 \quad \textrm{as} \,\, n\to\infty.
			\end{equation}
			If, in addition, $I(\omega_n)\to\mu$, then
			\begin{equation}\label{eq:cpt_2}
				\inf_{\omega \in \calO^+_{\mu,\varepsilon}} \|\omega_n - \omega\|_{1}\to  0 \quad \textrm{as} \,\, n\to\infty.
			\end{equation}
		\end{proposition}
		
		\begin{proof}
			It suffices to prove subsequential convergence. 
			Set $\tld \omega_n(x) :=\tau_{0,b_n}\omega_n (x) = \omega_n(x_1, x_2 - b_n)$.
			By translation invariance, $\tld \omega_n$ satisfies the bounds in \eqref{eq:seq_cpt}, and
			\begin{equation}\label{eq:sec_cpt_2}
				\limsup_{n\to\infty} \|x_2\tld\omega_n\|_1 \le \mu, \quad \liminf_{n\to\infty}E[\tld\omega_n]\ge M_{\mu,\varepsilon}.
			\end{equation}
			Write $\tld\omega_n = \tld\omega_{n, +} - \tld\omega_{n, -}$. 
			Since $\int \tld \omega_n  = 0$, we have $\int \tld\omega_{n, \pm}   = \|\tld \omega_n\|_1/2$.
			By passing to a subsequence if necessary, we find $\tld\omega_{n, \pm} \rightharpoonup^* f_\pm$ in  $L^\infty(\calC)$, where $0 \le f_\pm \le \varepsilon^{-2}$.
			Set $\omega^* := f_+ - f_-$.
			
			By \eqref{eq:seq_cpt} and the first bound in \eqref{eq:sec_cpt_2}, we see $\omega^* \in \calF_{\mu,\varepsilon}$, and then Lemma~\ref{lem:weak_energy_conv} gives $E[\tld \omega_n] \to E[\omega^*]$.
			Since $M_{\mu,\varepsilon}=\sup_{\theta\in\calF_{\mu,\varepsilon}}E[\theta]$, it follows from the second bound in \eqref{eq:sec_cpt_2} that $E[\omega^*] = M_{\mu,\varepsilon}$.
			Hence $\omega^*$ is a maximizer.
			Moreover, we may choose a vertical translation of $\omega^*$, still denoted by $\omega^*$, so that $\|x_2\omega^*\|_1 = \mu$ and $\omega^* \in \calS_{\mu, \varepsilon}$.
			By Proposition~\ref{prop:orientation}, $|I(\omega^*)|  = \mu$.
			By Lemmas~\ref{lem:patch} and~\ref{lem:mass_sat}, $\omega^*$ has the form $\omega^* = \omega^*_+ - \omega^*_-  = \varepsilon^{-2} \left( \mathbf 1_{\Omega_{*,+}} - \mathbf 1_{\Omega_{*,-}} \right)$, where $|\Omega_{*, \pm}| = \varepsilon^2 / 2$.
			
			We now prove the $L^1$-convergence.
			Since $f_+ - f_- = \omega^*$, we have $f_+ \ge \omega^*_+$ and $f_-\ge \omega^*_-$.
			This with the lower semicontinuity, we see
			\[
			1 = \int_{\calC}|\omega^*|\,dx \le \int_{\calC}(f_++f_-)\,dx \le \liminf_{n\to\infty}\|\tld \omega_n\|_1 \le 1.
			\]
			Therefore $f_{\pm} = \omega^*_{\pm}$ a.e. in $\calC$ and $\|\tld\omega_n\|_1 \to 1$.
			Since $\tld\omega_{n, +} \rightharpoonup^* \varepsilon^{-2}\mathbf 1_{\Omega_{*,+}}$ and $\limsup_{n\to\infty}\|\tld\omega_{n, +}\|_\infty \le \varepsilon^{-2}$, 
			Lemma~\ref{lem:comple} gives $\tld\omega_{n, +} \to \varepsilon^{-2}\mathbf 1_{\Omega_{*,+}}$ in $L^1(\calC)$.
			By applying the same argument to $\tld\omega_{n, -}$, we deduce
			\begin{equation}\label{eq:cpt_L1}
				\tld\omega_n \to \omega^* \quad \text{in } \,\, L^1(\calC).
			\end{equation}
			It is clear that $\omega^* \in \calO_{\mu,\varepsilon}$.
			Since $\calO_{\mu,\varepsilon}$ is invariant under translations, $\tau_{0,-b_n}\omega^* \in \calO_{\mu,\varepsilon}$.
			Then by \eqref{eq:cpt_L1}, $\inf_{\omega \in \calO_{\mu,\varepsilon}} \left\|\omega_n - \omega \right\|_1 \le \left\| \omega_n -  \tau_{0,- b_n}  \omega^* \right\|_1 = \left\|  \tld\omega_n - \omega^*  \right\|_1 \to 0$.
			This proves \eqref{eq:cpt_1}.
			
			Now assume in addition that $I(\omega_n) \to \mu$.
			Since $\int \omega_n = 0$, we have $I(\tld\omega_n) = I(\omega_n) \to \mu$.
			Together with the first inequality in \eqref{eq:sec_cpt_2}, this implies $\|x_2\tld\omega_n \|_1 - I(\tld \omega_n) \to 0$.
			Writing $\tld\omega_n = \tld\omega_{n, +} - \tld\omega_{n, -}$, we find
			\[
			\int_{\{x_2<0\}} |x_2| \tld\omega_{n, +}\,dx + \int_{\{x_2>0\}} x_2 \tld \omega_{n, -}\,dx \to 0.
			\]
			It follows that $f_+ = 0$ a.e. on $\{x_2 \le 0\}$ and $f_- = 0$ a.e. on $\{x_2 \ge 0\}$.
			Thus $x_2\omega^*(x) = x_2(f_+ - f_-)(x) \ge 0$, and hence $\omega^* \in \calO^+_{\mu,\varepsilon}$.
			Since  $\tau_{0,-b_n} \omega^* \in \calO^+_{\mu,\varepsilon}$, \eqref{eq:cpt_2} follows from \eqref{eq:cpt_L1}.
		\end{proof}
		
		\subsection{Conditional Lyapunov stability}
		\begin{definition}[Orbital distance]
			We define the distances to the corresponding orbits by
			\begin{align*}
				d_{\calO}(v) := &\inf_{\omega\in\calO_{\mu,\varepsilon}} \left\{ \|v - \omega\|_1 + |I(v)-I(\omega)| \right\},    \\
				d_{\calO^+}(v) := &\inf_{\omega\in\calO^+_{\mu,\varepsilon}} \left\{ \|v - \omega\|_1 + |I(v) - \mu| \right\} .
			\end{align*}
			We also define the translation invariant vertical moment by
			\[
			\calT(v) := \inf_{b\in\bbR} \|x_2(\tau_{0,b}v)\|_1 = \inf_{b\in\bbR} \|x_2 v(x_1, x_2 - b)\|_1.
			\]
		\end{definition}
		
		We first show that $L^1$-closeness to the orbit forces the convergence of energy.
		
		\begin{lemma}\label{lem:en_near_orb}
			Let $\{v_n\}_n\subset L^1\cap L^\infty(\calC)$ satisfy  $\int v_n  = 0$, $\limsup_{n\to\infty} \|v_n\|_\infty \le \varepsilon^{-2}$ and
			\begin{equation}\label{eq:seq_tras}
				\limsup_{n\to\infty} \calT(v_n)\le\mu.
			\end{equation}
			If $\inf_{\omega\in\calO_{\mu,\varepsilon}}\|v_n-\omega\|_1 \to 0$, then $E[v_n] \to M_{\mu,\varepsilon} = \sup_{\omega \in \calF_{\mu, \varep}}E[\omega]$.
		\end{lemma}
		
		\begin{proof}
			It suffices to show that every subsequence has a further subsequence along which $E[v_n] \to M_{\mu,\varepsilon}$.
			Choose $\xi_n\in \calO_{\mu,\varepsilon}$ and $b_n \in \bbR$ such that $\|v_n - \xi_n\|_1 \to 0$ and
			\[
			\|x_2(\tau_{0,b_n}v_n)\|_1 \le \mu + \frac{1}{n}.
			\]
			Since $\|\xi_n\|_1=1$, we have $\|v_n\|_1 \to 1$.
			Set $\tld v_n := \tau_{0,b_n} v_n$ and $\tld{\xi}_n := \tau_{0,b_n} \xi_n$, and note that
			\[
			\limsup_{n\to\infty} \|x_2 \tld v_n \|_1 \le \mu, \quad  \|\tld v_n - \tld{\xi}_n\|_1 \to 0.
			\]
			Since $\tld\xi_n \in \calO_{\mu, \varepsilon}$, by Lemma~\ref{lem:supp}, there is $d_n \in \bbR$ such that $\supp \tld\xi_n \subset \bbT \times [d_n - \bar{R}, d_n + \bar{R}]$, where $\bar{R} = \bar{R}(\mu)$.
			Since $\|\tld\xi_n\|_1 = 1$ (Lemma~\ref{lem:mass_sat}) and $\|\tld v_n - \tld{\xi}_n\|_1 \to 0$, we have
			\[
			\int_{\supp \tld\xi_n}|\tld v_n|\,dx \ge \int_{\supp \tld\xi_n}|\tld\xi_n|\,dx - \|\tld v_n-\tld\xi_n\|_1 = 1 - \|\tld v_n-\tld\xi_n\|_1 \ge \frac{1}{2}
			\]
			for sufficiently large $n$.
			Therefore
			\[
			\mu + 1 \ge \|x_2\tld v_n\|_1 \ge \int_{\supp \tld\xi_n}|x_2|\,|\tld v_n|\,dx \ge \frac{|d_n| - \bar{R}}{2}.
			\]
			Thus $\{d_n\}$ is bounded.
			Consequently, $\sup_n \|x_2\tld\xi_n \|_1 < \infty$.
			
			On the other hand, by Lemma~\ref{lem:weak_energy_conv}, after passing to a further subsequence, there exists $v^* \in \calF_{\mu,\varepsilon}$ such that $\tld{v}_n \rightharpoonup^* v^*$ in $L^{\infty}(\calC)$ and $E[\tld{v}_n]\to E[v^*]$.
			Moreover, $\tld{\xi}_n$ has the same weak-$^*$ limit $v^*$.
			Indeed, set $r_n := \tld{\xi}_n - \tld{v}_n$. 
			Then $\|r_n\|_1 \to 0$ and $\sup_n\|r_n\|_\infty < \infty$.
			For any $\varphi\in L^1(\calC)$ and any $K>0$,
			\[
			\left|\int_{\calC}\varphi r_n\,dx\right| \le K\|r_n\|_1 + \|r_n\|_\infty \int_{\{|\varphi|>K\}}  |\varphi|\,dx.
			\]
			Taking first $n\to\infty$ and then $K\to\infty$, we obtain  $\tld \xi_n - \tld v_n \rightharpoonup^* 0 $ in $L^{\infty}(\calC)$, i.e. $\tld \xi_n\rightharpoonup^* v^*$.
			
			Finally, since $\{\tld\xi_n\}$ has uniform $L^1$, $L^\infty$ bounds and $\sup_n\|x_2\tld\xi_n\|_1<\infty$, by Lemma~\ref{lem:weak_energy_conv} we have $E[\tld\xi_n]\to E[v^*]$.
			Since $E[\tld\xi_n]=E[\xi_n]=M_{\mu,\varepsilon}$, we obtain $E[v_n] = E[\tld v_n]\to M_{\mu,\varepsilon}$.
		\end{proof}
		
		\begin{theorem}[Stability with moment condition]\label{thm:con_stabi}
			For every $\eta > 0$, there exists $\delta = \delta(\eta, \mu, \varepsilon) > 0$ such that the following holds:
			Let $0 \in J \subset \bbR$ be an interval, and let $\theta(t)$ be a global solution to \eqref{eq:Euler_cyl} with initial data $\theta(0) \in L^1 \cap L^{\infty}(\calC)$.
			Assume $\int \theta(0,x) = 0$,
			\[
			\|\theta(0)\|_{\infty} \le \varepsilon^{-2} +  \delta, \quad d_{\calO}(\theta(0)) < \delta, \quad \textrm{and} \quad \sup_{t \in J}\calT(\theta(t)) \le \mu +  \delta.
			\]
			Then $\sup_{t \in J} d_{\calO}(\theta(t)) < \eta$.
			%\begin{equation}\label{eq:sta_res}
			%\sup_{t \in J} d_{\calO}(\theta(t)) < \eta.
			%\end{equation}
			If $d_{\calO^+}(\theta(0))<\delta$ is assumed instead, then $\sup_{t \in J }d_{\calO^+}(\theta(t)) < \eta$.
		\end{theorem}
		
		\begin{proof}
			We argue by contradiction.
			Suppose first that $\calO$-stability fails.
			Then, there exist $\eta_0>0$, intervals $J_n \ni 0$, a sequence of solutions $\theta_n(t)$, and times $t_n \in J_n $ such that
			\begin{equation}\label{eq:sta_ass_4}
				\int_{\calC}\theta_n(0,x)\,dx = 0,\quad \|\theta_n(0)\|_\infty\le\varepsilon^{-2} + \frac{1}{n}, \quad d_{\calO}(\theta_n(0)) \to 0
			\end{equation}
			and
			\begin{equation}\label{eq:sta_ass_3}
				\sup_{t \in J_n}  \calT(\theta_n(t)) \le \mu + \frac{1}{n},
			\end{equation}
			but
			\begin{equation}\label{eq:stab_bad}
				d_{\calO}(\theta_n(t_n))\ge\eta_0.
			\end{equation}
			Set $v_n := \theta_n(t_n)$.
			Since $d_{\calO}(\theta_n(0)) \to 0$ and $\|\omega\|_1 = 1$ for $\omega \in \calO_{\mu,\varepsilon}$, we have $\limsup_{n\to\infty}\|\theta_n(0)\|_1 \le 1$.
			By the conservation property, 
			\[
			\int_{\calC}v_n\,dx = 0, \quad \limsup_{n\to\infty} \|v_n\|_1 \le 1, \quad \limsup_{n\to\infty}    \|v_n\|_\infty \le \varepsilon^{-2}.
			\]
			From \eqref{eq:sta_ass_3}, we also have $\limsup_{n\to\infty} \calT(v_n) \le \mu$, that is, there exists $b_n \in \bbR$ such that
			\[
			\limsup_{n\to\infty} \|x_2(\tau_{0,b_n}v_n)\|_1 = \limsup_{n\to\infty} \|x_2v_n(x_1, x_2 - b_n)\|_1 \le \mu.
			\]
			By Lemma~\ref{lem:en_near_orb} with \eqref{eq:sta_ass_4}--\eqref{eq:sta_ass_3}, we have $E[v_n] =  E[\theta_n(0)] \to  \sup_{\omega \in \calF_{\mu, \varepsilon}}E[\omega]$.
			Hence $\{v_n\}_n$ satisfies the assumptions of Proposition~\ref{prop:seq_cpt}.
			
			On the other hand, since $d_{\calO}(\theta_n(0)) \to 0$ and $|I(\omega)| = \mu$ for every $\omega \in \calO_{\mu,\varepsilon}$, after passing to a subsequence, either $I(v_n) = I(\theta_n(0)) \to \mu$ or $I(v_n) = I(\theta_n(0)) \to -\mu$.
			If $I(v_n)\to\mu$, then \eqref{eq:cpt_2} of Proposition~\ref{prop:seq_cpt} holds, and hence
			\[
			d_{\calO}(v_n) \le \inf_{\omega\in\calO^+_{\mu,\varepsilon}} \left\{ \|v_n-\omega\|_1 + |I(v_n) - \mu| \right\} \to 0,
			\]
			which contradicts \eqref{eq:stab_bad}.
			If $I(v_n)\to-\mu$, we apply Proposition~\ref{prop:seq_cpt} to $\calR v_n$ to obtain $d_{\calO}(v_n) \to 0$, again contradicting \eqref{eq:stab_bad}.
			This proves the $\calO$-stability.
			The $\calO^+$-stability is proved in the same way.
		\end{proof}
		
		\subsubsection{Stability with support condition}
		For a more intuitive version of Theorem~\ref{thm:con_stabi}, we now replace $\calT(\theta(t))$ by a bound on the vertical diameter of $\supp \theta(t)$.
		To be concise, we consider patch perturbations having the same distribution as the maximizers:
		\begin{equation}\label{eq:supp_patch}
			\theta(0) = \varepsilon^{-2} \left( \mathbf 1_{A_0^+}-\mathbf 1_{A_0^-} \right), \quad A_0^+ \cap A_0^- = \emptyset, \quad |A_0^{\pm}| = \frac{\varepsilon^2}{2}.
		\end{equation}

		\begin{theorem}[Stability with support condition]\label{thm:supp_stab}
			There exists $c_0 = c_0(\mu) > 0$ such that the following holds:
			For every $\eta > 0$, there exists $\delta = \delta(\eta, \mu, \varepsilon) > 0$ such that, if $0 \in J \subset \bbR$ is an interval and $\theta(t)$ is the global solution to \eqref{eq:Euler_cyl} whose initial data satisfies \eqref{eq:supp_patch}, $d_{\calO}(\theta(0)) < \delta$, and
			\begin{equation}\label{eq:supp_stab}
				\sup_{t \in J} \sup_{x,y \in \supp\theta(t)} |x_2-y_2| \le c_0|\log\varepsilon|,
			\end{equation}
			then $\sup_{t \in J} d_{\calO}(\theta(t)) < \eta$.
			If $d_{\calO^+}(\theta(0)) < \delta$ is assumed instead, then $\sup_{t \in J} d_{\calO^+}(\theta(t)) < \eta$.
		\end{theorem}
		
		\begin{remark}
			The patch assumption \eqref{eq:supp_patch} is imposed only for simplicity.
			Theorem~\ref{thm:supp_stab} remains valid for general $L^1 \cap L^{\infty}(\calC)$ initial data satisfying
			\[
			\int_{\calC} \theta(0)\,dx = 0, \quad \| \theta(0)\|_{\infty} \le \varepsilon^{-2} + \delta, \quad d_{\calO}(\theta(0)) < \delta,
			\]
			with \eqref{eq:supp_stab} in place of $\sup_{t \in J} \calT(\theta(t)) \le \mu + \delta$ in Theorem~\ref{thm:con_stabi}.
			The only major change in the proof is to replace the set exchange \eqref{eq:0804_03} in Lemma~\ref{lem:supp_en_loss} by a suitable density exchange.
		\end{remark}
		
		We need a local comparison of the maximal energy levels.
		Recall that $M_{\mu, \varepsilon} = \sup_{\omega \in \calF_{\mu,\varepsilon}} E[\omega]$.

		\begin{lemma}\label{lem:M_mu_lip}
			There exist $r = r(\mu) > 0$, $C = C(\mu) > 0$, and $\varepsilon_0 = \varepsilon_0 (\mu) > 0$ such that, for every $\varepsilon\in(0,\varepsilon_0)$ and $\nu \in (\mu - r, \mu + r)$,
			\[
			M_{\nu,\varepsilon} \le M_{\mu,\varepsilon} + C|\nu - \mu|.
			\]
		\end{lemma}
		
		\begin{proof}
			The assertion is trivial when $\nu \le \mu$.
			Let $\mu < \nu < \mu + r$, set $h := \nu - \mu$, and choose $\omega \in \calS_{\nu,\varepsilon}^{\rm cen}$.
			By decreasing $r$ and $\varepsilon_0$ if necessary, Proposition~\ref{prop:cen_prop} gives
			\[
			\tld{\omega} := \tau_{0,-h} \omega_+ - \tau_{0,h}\omega_- \in \calF_{\mu,\varepsilon}.
			\]
			Indeed, $\|\omega_\pm\|_1 = 1/2$ and hence $I(\tld{\omega}) = I(\omega) - h = \mu$, while $x_2 \tld{\omega} \ge 0$ for sufficiently small $r$ and $\varepsilon$.
			Since the self-energies are invariant under translations,
			\[
			\begin{aligned}
				E[\omega] - E[\tld{\omega}]
				&= E_{\rm int} \left[ \tau_{0,-h}\omega_+, \tau_{0,h}\omega_- \right] - E_{\rm int}[\omega_+,\omega_-] \\
				&= -2\int_0^h \iint_{\calC \times \calC} \partial_2G(x-y-2se_2)\omega_+(x)\omega_-(y)\,dxdy\,ds.
			\end{aligned}
			\]
			For $x \in \supp\omega_+$, $y \in \supp\omega_-$, and $s \in [0,h]$, by \eqref{eq:cen_ver} of Proposition~\ref{prop:cen_prop},
			\[
			x_2-y_2-2s \ge 2\nu-2C\varepsilon-2h = 2\mu - 2C\varepsilon \ge \mu.
			\]
			Therefore, $0 < -\partial_{2}G(x-y-2se_2) \lesssim \coth\left(\mu/2 \right)$, and hence $0 \le E[\omega] - E[\tld{\omega}] \le Ch$.
			Since $M_{\nu,\varepsilon} = E[\omega]$ and $E[\tld{\omega}] \le M_{\mu,\varepsilon}$, the lemma is proved.
		\end{proof}
		
		The next lemma controls the stream function in \emph{wrong-signed} regions.
		
		\begin{lemma}\label{lem:wrong_gap}
			There exist $c = c(\mu) > 0$ and $c_1 = c_1(\mu) > 0$ such that the following holds:
			Let $0 < R < c |\log\varepsilon |$.
			If $\omega\in\calS_{\mu,\varepsilon}^{\rm cen}$ and $v\in L^1(\calC)$ satisfy
			\[
			\|v\|_1 = 1, \quad \|v - \omega\|_1<1, \quad \sup_{x,y \in \supp v} |x_2-y_2| \le R,
			\]
			then
			\begin{align*}
				\begin{cases}
					\psi[\omega](x) - W x_2 - \gamma_+ &\le - c_1|\log\varepsilon|    \quad \text{on} \quad \supp v \cap \{x_2 < 0\},\\ 
					-\psi[\omega](x) + Wx_2 - \gamma_- &\le - c_1 |\log\varepsilon| \quad \text{on} \quad \supp v \cap \{x_2 > 0\},
				\end{cases}
			\end{align*}
			where $W, \gamma_{\pm}$ as in \eqref{eq:EL2}.
		\end{lemma}
		
		\begin{proof}
			Due to Proposition~\ref{prop:cen_prop}, we write $\omega = \varepsilon^{-2} (\mathbf1_{\Omega_+} - \mathbf1_{\Omega_-})$, where $\Omega_\pm \subset \bbT \times [\pm \mu - C\varepsilon,  \pm\mu + C\varepsilon]$.
			Since $\|v\|_1=\|\omega\|_1=1$ and $\|v-\omega\|_1<1$, we see $\supp v \cap \supp \omega \neq \emptyset$.
			Hence, 
			\begin{equation}\label{eq:0803_2}
				\supp v \subset \bbT \times [-R - C, R + C].
			\end{equation}
			Since $\int \omega = 0$ and $\|x_2\omega\|_1 = \mu$, we find
			\begin{equation}\label{eq:0803_1}
				\left| \int_{\calC}|x_2-y_2| \omega(y)\,dy\right| = \left|\int_{\calC}\big( |x_2-y_2|-|x_2|\big)\omega(y)\,dy\right| \le \mu.
			\end{equation}
			For $x_2 \le 0$, $\inf_{y \in \Omega_+} | y_2 - x_2 |   \ge \mu/2$.
			Thus, by \eqref{green_sing_2} of Lemma~\ref{lem_green_decom},
			\[
			\int_{\Omega_+} H(x-y)\omega_+(y)\,dy\le C, \quad -\int_{\Omega_-} H(x-y)\omega_-(y)\,dy\le C,
			\]
			and hence, by \eqref{green_lin} in the same lemma and \eqref{eq:0803_1}, $\sup_{ x_2 \le 0 } \psi[\omega]\le C$.
			Similarly, $\sup_{ x_2 \ge 0} (-\psi[\omega])\le C$.
			Then, by \eqref{eq:0803_2} with the help of Lemmas~\ref{lem:W} and \ref{lem:gamma}, we obtain
			\[
			\begin{aligned}
				\psi[\omega](x) - Wx_2 - \gamma_+ &\le C(1+R) - \frac{1}{4\pi}|\log\varepsilon| \quad \textrm{ on } \,\,\supp v\cap\{x_2<0\},\\
				-\psi[\omega](x) + Wx_2 - \gamma_- &\le C(1+R) - \frac{1}{4\pi}|\log\varepsilon| \quad \textrm{ on } \,\, \supp v\cap\{x_2 > 0\}.
			\end{aligned}
			\]
			After decreasing $c$ and then $\varepsilon_0$ if necessary, the assumption $R < c|\log\varepsilon|$ yields the result.
		\end{proof}

		The following lemma quantifies the energy loss caused by \emph{wrong-signed} vorticity.
		
		\begin{lemma}\label{lem:supp_en_loss}
			There exist $c = c(\mu) > 0$ and $C_i = C_i(\mu) > 0$ ($i = 1, 2$) such that the following holds:
			Let $\omega \in \calS_{\mu,\varepsilon}^{\rm cen}$ and
			\begin{equation}\label{eq:0804_1}
				2\mu +  1 <  R <  c|\log\varepsilon|.
			\end{equation}
			Then there exists $\rho = \rho(\mu, \varepsilon, R) >  0$ such that, if
			\[
			v= v_+ - v_- = \varepsilon^{-2} \left(\mathbf1_{A_+} - \mathbf1_{A_-}\right), \quad A_+\cap A_- = \emptyset, \quad |A_\pm| = \frac{\varepsilon^2}{2},
			\]
			satisfy
			\begin{equation}\label{eq:supp_loc_clo}
				\|v - \omega\|_1 +  |I(v) - \mu| <  \rho
			\end{equation}
			and $\sup_{x,y \in\supp v} |x_2 -y_2|  \le R$, then
			\begin{equation}\label{eq:supp_en_loss}
				E[v]\le E[\omega] + C_1|I(v)-\mu| - C_2|\log\varepsilon| \left( \int_{\{x_2<0\}}v_+\,dx + \int_{\{x_2>0\}}  v_-\,dx \right).
			\end{equation}
		\end{lemma}

		\begin{proof}
			Let $c_1 > 0$ be the constant given by Lemma~\ref{lem:wrong_gap}.
			We shall choose $c \le c_1$ sufficiently small below.
			We take $\rho < 1$ such that Lemma~\ref{lem:wrong_gap} applies.
			Write $\omega = \varepsilon^{-2}(\mathbf1_{\Omega_+} - \mathbf1_{\Omega_-})$ and set
			\[
			P := A_+ \cap \{x_2 < 0\}, \quad N := A_- \cap \{x_2 > 0\}, \quad m := \varepsilon^{-2}\left( \left| P \right| + \left| N \right| \right).
			\]
			Thus, $m$ represents the amount of \emph{wrong-signed} mass:
			\[
			m = \int_{\{x_2 < 0\}}v_+\,dx + \int_{\{x_2 > 0\}}v_-\,dx.
			\]
			Since $P \subset A_+ \setminus \Omega_+$, $N \subset A_- \setminus \Omega_-$, and $|A_\pm| = |\Omega_\pm|$, there exist
			\[
			H_+ \subset \Omega_+ \setminus A_+, \quad H_- \subset \Omega_- \setminus A_-, \quad |H_+| = |P|, \quad |H_-| = |N|.
			\]
			Define
			\begin{equation}\label{eq:0804_03}
				\tld v := v + \varepsilon^{-2}\left(\mathbf1_{H_+} - \mathbf1_P + \mathbf1_N - \mathbf1_{H_-}\right).
			\end{equation}
			The patch $\tld v$ has the same distribution as $v$.
			Moreover, since $H_+ \subset \{x_2 > 0\}$ and $H_- \subset \{x_2 < 0\}$, we have
			\[
			x_2\tld v \ge 0, \quad \|\tld v - v\|_1 = 2m, \quad m \le \|v - \omega\|_1 < \rho.
			\]
			Set $\Delta I := I(\tld v) - I(v) \ge 0$.
			Then
			\[
			0 \le \Delta I = \varepsilon^{-2}\left(\int_{H_+}x_2\,dx - \int_Px_2\,dx + \int_Nx_2\,dx - \int_{H_-}x_2\,dx\right) \le C(1 + R)m.
			\]
			By \eqref{eq:EL2} and Lemma~\ref{lem:wrong_gap},
			\[
			\begin{aligned}
				&\int_{\calC} (\psi[\omega] - Wx_2)(\tld v - v)\,dx
				= \varepsilon^{-2}\left( \int_{H_+} - \int_P + \int_N - \int_{H_-} \right) (\psi[\omega] - Wx_2)\,dx\\
				\ge & \varepsilon^{-2}\bigg[ \gamma_+|H_+| -\left(\gamma_+ - c_1 |\log\varepsilon| \right) |P| +\left( -\gamma_- + c_1|\log\varepsilon| \right) |N| + \gamma_- |H_-| \bigg]\\
				= & c_1|\log\varepsilon| \varepsilon^{-2}\left( |P| + |N| \right) = c_1 |\log\varepsilon| m,
			\end{aligned}
			\]
			where we use $|H_+| = |P|$ and $|H_-| = |N|$ in the last line.
			Hence,
			\[
			\int_{\calC}\psi[\omega](\tld v - v)\,dx \ge c_1 |\log\varepsilon| m + W \Delta I \ge c_1 |\log\varepsilon| m.
			\]
			As in the proof of Lemma~\ref{lem:wrong_gap}, $\supp v \subset \bbT \times [-R - C, R + C]$.
			Therefore, by \eqref{eq:supp_loc_clo}, $\|x_2 ( v - \omega) \|_1 \le C(1+R)\rho$, and then by Lemma~\ref{lem_str_bdd},
			\[
			\|\psi[v-\omega]\|_\infty \le C\rho \left(1 + R + \log_+\frac{2\varepsilon^{-2}}{\rho}\right).
			\]
			Since $\|\tld v - v\|_1 = 2m$, by decreasing $\rho = \rho(\mu, \varepsilon, R)$,
			\[
			\begin{aligned}
				\left|E_{\rm int}[v - \omega, \tld v - v]\right| \le \|\psi[v - \omega] \|_\infty\|\tld v - v\|_1 \le \frac{c_1}{2}|\log\varepsilon|  m.
			\end{aligned}
			\]
			Since $E[\tld v-v] \ge 0$ and $E[\tld v] - E[v] = E_{\rm int}[v, \tld v - v] + E[\tld v - v]$, it follows that
			\[
			E[\tld v] - E[v] = \int_{\calC}\psi[\omega](\tld v - v)\,dx + E_{\rm int}[v - \omega, \tld v-v] + E[\tld v - v] \ge \frac{c_1}{2}|\log\varepsilon| m.
			\]
			On the other hand, since $x_2\tld v \ge 0$ and $I(\tld v) = I(v) + \Delta I \ge 0$, we have $\tld v \in \calF_{I(v)+\Delta I,\varepsilon}$.
			By taking smaller $\rho$, we may assume (recall that $\Delta I \le C(1 + R)m$ and $m < \rho$)
			\[
			|I(\tld v) - \mu| = |I(v) + \Delta I - \mu |  \le |I(v) - \mu| + C(1 + R)m < r(\mu),
			\]
			where $r(\mu)$ is given by Lemma~\ref{lem:M_mu_lip}.
			Therefore, by Lemma~\ref{lem:M_mu_lip},
			\[
			E[\tld v] \le M_{I(v)+\Delta I,\varepsilon} \le M_{\mu,\varepsilon} + C|I(v) - \mu| + C(1 + R)m
			\]
			for some $C = C(\mu) > 0$.
			This with the lower bound for $E[\tld v] - E[v]$ implies
			\[
			E[v] \le M_{\mu,\varepsilon} + C|I(v) - \mu| - \left(\frac{c_1}{2}|\log\varepsilon| - C(1 + R)\right)m.
			\]
			By choosing $c > 0$ small and then decreasing $\varepsilon_0$ if necessary, \eqref{eq:supp_en_loss} follows from \eqref{eq:0804_1}.
		\end{proof}
		
		In the proof below, the dependency of parameters is as follows:
		\[
		\mu \longrightarrow c_0 \longrightarrow \varepsilon \in (0, \varepsilon_0) \longrightarrow R = c_0|\log \varepsilon| \longrightarrow \rho \longrightarrow \sigma.
		\]
		
		\begin{proof}[\textbf{Proof of Theorem~\ref{thm:supp_stab}}]
			Let $c > 0$ be the constant in \eqref{eq:0804_1} of Lemma~\ref{lem:supp_en_loss}.
			We choose $c_0 \in (0, c)$ and then decrease $\varepsilon_0$ if necessary so that $2\mu +  1 <  c_0|\log\varepsilon|$ for every $\varepsilon \in (0, \varepsilon_0)$.
			Thus Lemma~\ref{lem:supp_en_loss} applies with $R := c_0|\log\varepsilon|$.

			We only prove the $\calO^+$-stability since $\calO$-stability follows in the similar way.
			Towards contradiction, suppose that $\calO^+$-stability fails.
			There exist $\eta_0 > 0$, intervals $J_n \ni 0$, and solutions $\theta_n(t)$ satisfying \eqref{eq:supp_stab} such that
			\begin{equation}\label{eq:supp_hypo}
				d_{\calO^+}(\theta_n(0)) \to 0, \quad \sup_{t \in J_n} d_{\calO^+}(\theta_n(t)) \ge \eta_0.
			\end{equation}
			We claim that $E[\theta_n(0)] \to M_{\mu,\varepsilon}$.
			Indeed, choose $\omega_n \in \calO^+_{\mu,\varepsilon}$ such that $\|\theta_n(0)-\omega_n\|_1 \to 0$, and apply a common translation to $\theta_n(0)$ and $\omega_n$ so that $\omega_n \in \calS_{\mu,\varepsilon}^{\rm cen}$.
			By Proposition~\ref{prop:cen_prop}, \eqref{eq:supp_stab} and $\|\theta_n(0)- \omega_n\|_1 < 1$, we have
			\[
			\supp\theta_n(0) \subset \bbT \times [-C_R, C_R]
			\]
			for some $C_R > 0$ independent of $n$.
			Hence
			\[
			\calT(\theta_n(0)) \le \|x_2\theta_n(0)\|_1 \le \mu + C_R\|\theta_n(0) - \omega_n\|_1 = \mu + o(1),
			\]
			and the claim follows from Lemma~\ref{lem:en_near_orb}.
			
			Let $\rho = \rho(\mu, \varepsilon, R) > 0$ be given by Lemma~\ref{lem:supp_en_loss}, and choose $\sigma > 0$ such that $2\sigma < \min\{\eta_0, \rho \}$.
			By the continuity of the flow and \eqref{eq:supp_hypo}, for all sufficiently large $n$, there exists $t_n \in J_n$ such that $d_{\calO^+}(\theta_n(t_n)) = \sigma$.
			Set $v_n := \theta_n(t_n)$ and choose $\xi_n \in \calO^+_{\mu,\varepsilon}$ such that
			\[
			\|v_n - \xi_n\|_1 + \left| I(v_n) - \mu \right| \le \sigma + \frac{1}{n} \le \frac{\rho}{2} + \frac{1}{n}.
			\]
			Here, we again assume that $\xi_n \in \calS_{\mu,\varepsilon}^{\rm cen}$ by applying a common translation to $v_n$ and $\xi_n$.
			Let
			\[
			m_n := \int_{\{x_2<0\}}(v_n)_+\,dx + \int_{\{x_2>0\}}(v_n)_-\,dx.
			\]
			By Lemma~\ref{lem:supp_en_loss} and $E[v_n] = E[\theta_n(0)] \to M_{\mu,\varepsilon}$, we have $|\log\varepsilon|m_n \lesssim M_{\mu,\varepsilon} - E[v_n] + C|I(v_n)-\mu| \to 0$, and thus $m_n \to 0$.
			As above, Proposition~\ref{prop:cen_prop}, the assumption \eqref{eq:supp_stab}, and $\|v_n - \xi_n\|_1 < 1$ imply $\supp v_n \subset \bbT \times [-C_R,C_R]$.
			Therefore,
			\begin{align*}
				0 \le \|x_2v_n\|_1-I(v_n) = &2\int_{\{x_2<0\}}|x_2|(v_n)_+\,dx+2\int_{\{x_2>0\}}x_2(v_n)_-\,dx  \lesssim  C_R m_n \to 0.
			\end{align*}
			Since $I(v_n) = I(\theta_n(0)) \to \mu$, this gives $\limsup_{n\to\infty}\|x_2v_n\|_1 \le \mu$.
			Proposition~\ref{prop:seq_cpt} with $b_n = 0$ now yields $\inf_{\omega \in \calO^+_{\mu,\varepsilon}}\|v_n-\omega\|_1 \to 0$ and hence $d_{\calO^+}(v_n) \to 0$, contradicting $d_{\calO^+}(v_n) = \sigma$.
			The theorem is proved.
		\end{proof}
		
		\section{Tilted vortex street}\label{sec:tilted_karman}
		
		In this section, by modifying the construction in \S\ref{sec:vari} and \S\ref{sec:unif}, we obtain vortex streets traveling in a \emph{tilted} direction.
		Since the existence and small-core analysis follow essentially the same arguments, we only emphasize the necessary modifications.

		We use the following conventions:
		\begin{itemize}
			\item \textbf{Orientation.}
			The positive vortex is located above the negative vortex.
			
			\item \textbf{Relative position.}
			The limiting relative position is written as 
			\[
			p_+ - p_- = (-a,h), \quad 0 < a < \pi, \quad h > 0.
			\]
			Here, $a$ is the horizontal distance over $\bbT$ and $h$ is the vertical distance.
			The horizontal component $-a$ is understood modulo $2\pi$.
			
			\item \textbf{Traveling velocity.}
			The traveling velocity is written as $W(1,\kappa)$, with $W > 0$ and $\kappa > 0$.
			
			\item \textbf{Representative in $\bbT$.}
			We denote by $\tld{x}_1 \in  (-\pi, \pi]$ the corresponding representative of $x_1$ in $\bbT$.
			
		\end{itemize}
		For the relative position $(-a,h)$, the signed representative difference
		\[
		d := \tld{p}_{+,1} - \tld{p}_{-,1}
		\]
		is either $-a$ or $2\pi - a$, depending on the placement of the pair relative to the cut in $\bbT$.
		
		\begin{figure}[h]
			\centering
			\includegraphics[scale=0.8]{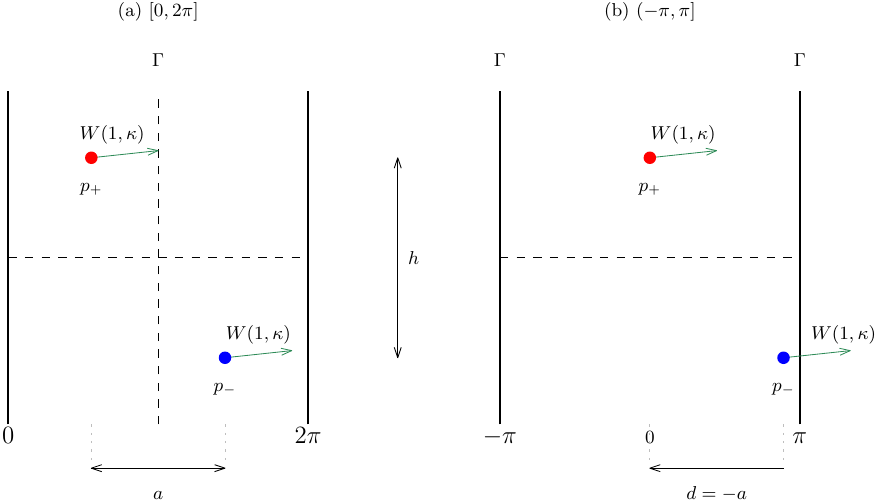}  
			\caption{The same periodic vortex street.
				(a) displays on $[0, 2\pi]$, with the cut $\Gamma$ at $x_1 = \pi$.
				(b) uses representatives in $(-\pi,\pi]$, with $p_+ = (0,h/2)$ and $p_- = (a,-h/2)$.
				Both figures have relative position $(-a,h)$ on $\calC$ and traveling velocity $W(1,\kappa)$. } \label{fig:pv-tilted}
		\end{figure} 
		
		The main goal of this section is to construct slightly tilted vortex streets for all sufficiently small $\kappa > 0$.
		
		\begin{theorem}[Tilted vortex streets]\label{thm:tilt_small}
			Let $\nu > 0$.
			There exists $\kappa_0 = \kappa_0(\nu) > 0$ such that, for each $\kappa \in (0,\kappa_0)$, there are uniquely determined numbers $a_\kappa \in (\pi/2,\pi)$ and $h_\kappa > 0$ satisfying
			\[
			h_\kappa = 2\nu + \kappa(2\pi - a_\kappa), \quad \sin a_\kappa = \kappa \sinh h_\kappa.
			\]
			For all sufficiently small $\varepsilon > 0$, there exist disjoint measurable sets $\Omega_{\pm}^{\varepsilon} \subset \calC$ such that
			\[
			\omega_\varepsilon = \varepsilon^{-2} \left( \mathbf 1_{\Omega_{+}^\varepsilon} - \mathbf 1_{\Omega_{-}^\varepsilon} \right), \quad |\Omega_{\pm}^{\varepsilon}| = \frac{\varepsilon^2}{2},
			\]
			is a traveling vortex patch in $\calC$, with traveling speed $W_\varepsilon > 0$ in direction $(1,\kappa)$.
			After a common translation (depending on $\varepsilon$), we have
			\[
			\varepsilon^{-2}\mathbf 1_{\Omega_{+}^\varepsilon} \rightharpoonup \frac{1}{2}\delta_{p_+}, \quad \varepsilon^{-2}\mathbf 1_{\Omega_{-}^\varepsilon} \rightharpoonup \frac{1}{2}\delta_{p_-}, \quad p_+ - p_- = (-a_\kappa,h_\kappa)
			\]
			as $\varepsilon \to 0$.
			Moreover, $a_\kappa \to \pi$, $h_\kappa \to 2\nu$ as $\kappa \to 0$.
		\end{theorem}
		
		\begin{remark}
			\begin{enumerate}
				\item  An explicit sufficient choice in Theorem~\ref{thm:tilt_small} is
				\[
				\kappa_0(\nu) = \min\left\{ \frac{\nu}{3\pi}, \frac{2}{\pi\sinh(8\nu/3)}
				\right\}.
				\]
				\item Placing the pair across the cut gives the representative difference $2\pi-a$ instead of $-a$.
				This allows a larger vertical separation under the tilted moment constraint; see \eqref{eq:h_tilt} below.
				
				\item For general values of $\kappa$, the limiting point vortex problem need not uniquely select the relative position; see \S \ref{sec:tilt_selection}.
			\end{enumerate}
		\end{remark}
		
		To get tilted streets, we consider the following weight for each $\kappa > 0$:
		\begin{equation}\label{eq:tilt_weight}
			\ell_\kappa(x) := x_2 - \kappa \tld{x}_1.
		\end{equation}
		Note that the case $\kappa = 0$ corresponds to the vertical moment, i.e. $\ell_0(x) = x_2$.
		\begin{remark}
The weight $\ell_\kappa$ is only affine away from the cut $\{ \tld{x}_1 = \pi \}$, that is, $\ell_\kappa$ is discontinuous on $\bbT$.
Thus, in general, $\|\ell_\kappa\omega\|_1$ is not preserved under horizontal translations.
		\end{remark}

		\subsection{Variational problem with tilted moments}\label{sec:tilt_setup}
		
		To prove Theorem~\ref{thm:tilt_small}, we first maximize the kinetic energy at fixed core size $\varepsilon > 0$.
		We fix a bound $\nu > 0$ for the tilted moment, keeping $\kappa > 0$ fixed.
		The point vortex analysis \S \ref{sec:tilt_pv} and \ref{sec:tilt_selection} will then identify the limiting relative position selected by these parameters and determine when the maximizers generate traveling waves.
		
		\begin{definition}
			For $\varepsilon > 0$, we define
			\begin{align*}
				\calF_{\nu,\varepsilon}^{\kappa} := \left\{ \omega \in L^1\cap L^\infty(\calC): \|\omega\|_1 \le 1,\,  \|\omega\|_\infty \le \varepsilon^{-2},\, \|\ell_\kappa \omega\|_1 \le \nu,\, \int_{\calC}\omega(x)\,dx = 0 \right\}.
			\end{align*}
			We also define
			\[
			\tld{\calS}_{\nu,\varepsilon}^{\kappa} := \left\{ \omega \in \calF_{\nu,\varepsilon}^{\kappa}: E[\omega] = \sup_{\theta\in\calF_{\nu,\varepsilon}^{\kappa}}E[\theta] \right\}.
			\]
		\end{definition}
		The tilted moment controls the absolute vertical moment $\|x_2\omega\|_1$.
		Indeed,
		\begin{equation}\label{eq:tilt_vert_bd}
			|x_2| \le |\ell_\kappa(x)| + \pi\kappa,
		\end{equation}
		and hence $\|x_2\omega\|_1 \le \nu + \pi\kappa$ for every $\omega\in\calF_{\nu,\varepsilon}^{\kappa}$.
		It follows from (the proof of) Proposition~\ref{prop:existence} that $\tld{\calS}_{\nu,\varepsilon}^{\kappa} \neq \emptyset$.
		Note that the tilted-moment constraint is preserved under weak-$^*$ limits.
		Moreover, as in \S \ref{sec:vari}, the tilted moment constraint can be saturated by a vertical translation, and then we can define
		\[
		\calS_{\nu,\varepsilon}^{\kappa} := \left\{ \omega \in \tld{\calS}_{\nu,\varepsilon}^{\kappa}: \|\ell_\kappa \omega\|_1 = \nu \right\} \neq \emptyset.
		\]
		
		\begin{lemma}\label{lem:tilt_en_low}
			There exist $\varepsilon_0 = \varepsilon_0(\nu,\kappa)>0$ and $C = C(\nu,\kappa)>0$ such that, for any $\varepsilon \in(0,\varepsilon_0)$ and any $\omega\in\calS_{\nu,\varepsilon}^{\kappa}$,
			\begin{equation}\label{eq:tilt_en_low}
				E[\omega]
				\ge
				\frac{1}{8\pi}|\log\varepsilon|-C.
			\end{equation}
		\end{lemma}
		
		\begin{proof}
			Let $\rho=\rho_+-\rho_-$, where $\rho_\pm$ are uniform disk patches of density $\varepsilon^{-2}$ and mass $1/2$, centered at $(0,\nu)$ and $(0,-\nu)$, respectively.
			For sufficiently small $\varepsilon$, one may check $\rho\in\calF_{\nu,\varepsilon}^{\kappa}$.
			The proof of Lemma~\ref{lem:en_low} then applies without change.
		\end{proof}

		\begin{proposition}\label{prop:tilt_EL}
			There exists $\varepsilon_0=\varepsilon_0(\nu,\kappa)>0$ such that the following holds:
			For any $\varepsilon \in (0, \varepsilon_0)$ and any $\omega\in\calS_{\nu, \varepsilon}^{\kappa}$, there exist measurable sets $\Omega_\pm \subset \calC$ and constants
			\[
			W \ge 0, \quad \gamma_+,\gamma_-\in\bbR, \quad \gamma_+ + \gamma_- \ge 0,
			\]
			such that
			\begin{equation}\label{eq:tilt_patch}
				\omega = \varepsilon^{-2} \left( \mathbf 1_{\Omega_+} - \mathbf 1_{\Omega_-} \right), \quad |\Omega_+| = |\Omega_-| \le \frac{\varepsilon^2}{2},
			\end{equation}
			and
			\begin{equation}\label{eq:tilt_EL}
				\omega = \varepsilon^{-2} \left( \mathbf 1_{\{ \psi[\omega] - W |\ell_\kappa| - \gamma_+ > 0 \}} - \mathbf 1_{\{ -\psi[\omega] - W |\ell_\kappa| - \gamma_- > 0 \}} \right).
			\end{equation}
		\end{proposition}

		\begin{proof}
			The proofs of Lemmas~\ref{lem:patch}, \ref{lem:local_EL}, \ref{lem:impulse_multipl}, and \ref{lem:mass_multip} apply with $|x_2|$ replaced by $|\ell_\kappa|$.
			We only explain the points where the argument requires a modification, and only consider the positive part since the negative part can be treated identically.
			
			We first show \eqref{eq:tilt_patch}.
			Since every level set of $|\ell_\kappa|$ has measure zero, the auxiliary sets in the proof of Lemma~\ref{lem:patch} can be chosen with different averages of $|\ell_\kappa|$.
			Suppose that $K := \{ 0 < \omega < \varepsilon^{-2} \}$ has positive measure.
			The same variational argument in the proof of Lemma~\ref{lem:patch} gives constants $W, \gamma \in \bbR$ such that
			\[
			\Psi := \psi[\omega] - W|\ell_\kappa| - \gamma = 0
			\quad
			\textrm{a.e. on } K.
			\]
			Let
			\[
			\Gamma := \{ \tld{x}_1 = \pi \} \quad \textrm{and} \quad L_\kappa := \{ \ell_\kappa = 0 \}.
			\]
			Note that the set $K \setminus (\Gamma \cup L_\kappa)$ has positive measure, and $|\ell_\kappa|$ is locally affine on this set.
			Hence, we see $0 < \omega = -\Delta \Psi = 0$ a.e. on $K\setminus(\Gamma\cup L_\kappa)$,  which is a contradiction.
			
			We next verify the analogue of Lemma~\ref{lem:no_strip} (no strip lemma) needed to apply Lemma~\ref{lem:aux_basis}.
			Suppose that we cannot choose two points $y,z \notin \Gamma$, each satisfying the (measure-theoretic) boundary condition (for $\Omega_+$) in Lemma~\ref{lem:aux_basis}, such that
			\[
			|\ell_\kappa(z)| > |\ell_\kappa(y)|.
			\]
			Then, at all such boundary points, $|\ell_\kappa| \equiv c$ for some $c \in \bbR$.
			The two lines $\{ \ell_\kappa = \pm c \}$ divide $(-\pi, \pi] \times \bbR$ into three regions.
			Since $\Omega_+$ has finite positive measure, it follows that $\Omega_+ = \{ |\ell_\kappa| < c \}$.
			In particular,
			\[
			4\pi c = |\Omega_+| \le \frac{\varepsilon^2}{2},
			\]
			so $c = O(\varepsilon^2)$.
			Therefore, the proof of Lemma~\ref{lem:no_strip} applies in the similar way, and gives $E[\omega_+] \le C_{\kappa}$ and
			\[
			E[\omega] \le \frac{1}{16\pi}|\log\varepsilon| + C,
			\]
			contradicting \eqref{eq:tilt_en_low} for sufficiently small $\varepsilon$.
			Hence, for each of $\Omega_\pm$, such points $y$ and $z$ can be chosen away from $\Gamma$.
			%Therefore, we may choose boundary points $y, z \notin \Gamma$ such that $|\ell_\kappa(z)| > |\ell_\kappa(y)|$.
			Then the auxiliary function construction (Lemma~\ref{lem:aux_basis}) and the local perturbation argument (Lemma~\ref{lem:local_EL}) apply with $|\ell_\kappa|$ in place of $|x_2|$.
			Hence, modulo null sets,
			\[
			\Omega_\pm = \{ \pm \psi[\omega] - W_\pm |\ell_\kappa| - \gamma_\pm > 0 \}.
			\]
			
			The cross perturbations in Lemma~\ref{lem:impulse_multipl} give $W_+ = W_-$, while the mass exchange argument in Lemma~\ref{lem:mass_multip} gives $\gamma_+ + \gamma_- \ge 0$.
			Finally,
			\[
			|\ell_\kappa(x)| \ge |x_2|-\pi|\kappa| \to  \infty \quad \textrm{as } \,\, |x_2|\to\infty.
			\]
			Since $\psi[\omega]\in L^\infty(\calC)$, $W < 0$ would imply that the positive superlevel set in \eqref{eq:tilt_EL} has infinite measure.
			Therefore, $W_+ = W_- =: W \ge0 $, and \eqref{eq:tilt_EL} follows.
		\end{proof}
		
		\subsubsection{Quantitative estimates}\label{sec:tilt_est}
		From now on, we assume
		\begin{equation}\label{eq:tilt_gap}
			\nu > 3\pi \kappa.
		\end{equation}
		The assumption \eqref{eq:tilt_gap} is required only to obtain \eqref{eq:tilt_W_unif} from \eqref{eq:tilt_gap_ineq}.

		\begin{proposition}\label{prop:tilt_unif}
			There exist $\varepsilon_0 = \varepsilon_0(\nu,\kappa) > 0$, $C = C(\nu,\kappa) > 0$, $R = R(\nu,\kappa) > 0$, and $c = c(\nu, \kappa) > 0$ such that the following holds:
			For any $\varepsilon \in (0,\varepsilon_0)$ and any
			\[
			\omega = \varepsilon^{-2} \left( \mathbf 1_{\Omega_+} - \mathbf 1_{\Omega_-} \right) \in \calS_{\nu,\varepsilon}^{\kappa},
			\]
			we have
			\begin{equation}\label{eq:tilt_mass_asymp}
				\|\omega_+\|_1^2 = \|\omega_-\|_1^2 \ge \frac{1}{4} - \frac{C}{|\log\varepsilon|},
			\end{equation}
			\begin{equation}\label{eq:tilt_W_unif}
				0 \le W \le C,
			\end{equation}
			and
			\begin{equation}\label{eq:tilt_gamma_unif}
				\frac{1}{4\pi}|\log\varepsilon|-C    \le \gamma_\pm \le \frac{1}{4\pi}|\log\varepsilon|+C.
			\end{equation}
			Moreover,
			\begin{equation}\label{eq:tilt_mass_sat}
				\|\omega_+\|_1 = \|\omega_-\|_1 = \frac{1}{2},
			\end{equation}
			\begin{equation}\label{eq:tilt_supp}
				\supp\omega \subset \{ x\in\calC: |x_2| \le R
				\},
			\end{equation}
			\begin{equation}\label{eq:tilt_core}
				\operatorname{diam}(\Omega_\pm) \le C\varepsilon,
			\end{equation}
			and
			\begin{equation}\label{eq:tilt_sep}
				\dist(\Omega_+,\Omega_-)
				\ge
				c.
			\end{equation}
		\end{proposition}
		
		\begin{proof}
			By \eqref{eq:tilt_vert_bd},
			\begin{equation}\label{eq:tilt_vert_bd_2}
				\|x_2\omega\|_1 \le \nu + \pi \kappa.
			\end{equation}
			Thus, the proofs of Lemmas~\ref{lem:mass_asymp} and \ref{lem:weak_bulk} apply with $\mu$ replaced by $\nu + \pi \kappa$.
			This implies \eqref{eq:tilt_mass_asymp} and a constant $R_0=R_0(\nu,\kappa)>0$ such that
			\begin{equation}\label{eq:tilt_weak_bulk}
				\int_{\Omega_\pm\cap\{|x_2|\le R_0\}} \omega_\pm\,dx > \frac{1}{4} - \frac{C}{|\log\varepsilon|}.
			\end{equation}
			
			We next prove \eqref{eq:tilt_W_unif}.
			Let
			\[
			\Psi_+
			:=
			\left(
			\psi-W|\ell_\kappa|-\gamma_+
			\right)_+,
			\quad
			\Psi_-
			:=
			\left(
			-\psi-W|\ell_\kappa|-\gamma_-
			\right)_+.
			\]
			Note that  $\partial_2|\ell_\kappa| = \sgn(\ell_\kappa)$ a.e. in $(-\pi, \pi] \times \bbR$.
			Repeating the calculation leading to \eqref{eq:poho_rhs_2}, we obtain
			\begin{equation}\label{eq:tilt_W_id}
				W \int_{\calC} x_2\sgn(\ell_\kappa) |\omega|\,dx = I_1 + I_2 + I_3,
			\end{equation}
			where
			\[
			I_1 := \int_{\calC} x_2\omega\partial_2 \psi\,dx, \quad I_2 := \int_{\calC} \Psi_+\omega_+\,dx, \quad I_3 := \int_{\calC} \Psi_-\omega_-\,dx.
			\]
			Since $ x_2 \sgn(\ell_\kappa) = |\ell_\kappa| + \kappa\tld{x}_1 \sgn(\ell_\kappa)$ and $\|\ell_{\kappa} \omega\|_1 = \nu $ for $\omega \in \calS_{\nu, \varepsilon}^\kappa$, we have
			\[
			W \left( \nu + \kappa K_\omega \right) = I_1 + I_2 + I_3,
			\]
			where
			\[
			K_\omega := \int_{\calC} \tld{x}_1 \sgn(\ell_\kappa) |\omega|\,dx, \quad |K_\omega | \le \pi,
			\]
			and hence
			\begin{equation}\label{eq:tilt_W_id_2}
				W (\nu - \pi \kappa) \le W (\nu + \kappa K_{\omega}) = I_1 + I_2 + I_3.
			\end{equation}
			The symmetrization argument in the proof of Lemma~\ref{lem:W}, together with \eqref{eq:tilt_vert_bd_2}, we get $|I_1| \le C$.
			To estimate $I_2$, we define
			\[
			D_+ := \left\{ \psi - W|x_2| - \gamma_+ -W \pi \kappa  > 0 \right\}, \quad  \Theta_+ := \left( \psi - W|x_2| - \gamma_+ - W \pi \kappa \right)_+.
			\]
			Since $||\ell_\kappa|(x)-|x_2|| \le \pi \kappa $, we have $D_+ \subset \Omega_+$ and
			\[
			\Psi_+ \le \Theta_+ + 2W \pi\kappa \quad \textrm{on }\,\, \Omega_+.
			\]
			Moreover,
			\[
			-\Delta\Theta_+ = \varepsilon^{-2} + 2W\delta_{\{x_2=0\}} \quad \textrm{in } D_+, \quad |D_+| \le \frac{\varepsilon^2}{2}.
			\]
			Then as in the proof of Lemma~\ref{lem:W}, we get
			\[
			I_2 \le \frac{1}{32\pi} + CW\varepsilon^{1/2} + \pi \kappa W.
			\]
			The same estimate holds for $I_3$.
			Combining these estimates with \eqref{eq:tilt_W_id_2}, we obtain
			\begin{equation}\label{eq:tilt_gap_ineq}
				W \left( \nu - 3\pi \kappa  - C\varepsilon^{1/2} \right) \le C.
			\end{equation}
			By \eqref{eq:tilt_gap}, this proves \eqref{eq:tilt_W_unif} for all sufficiently small $\varepsilon$.
			
			With \eqref{eq:tilt_W_unif} established, the proof of Lemma~\ref{lem:gamma} applies with $|\ell_\kappa|$ in place of $|x_2|$.
			Hence, we have \eqref{eq:tilt_gamma_unif}.
			Indeed, the identity
			\[
			2E[\omega] = M(\gamma_++\gamma_-) + W\nu + I_2 + I_3, \quad M := \|\omega_+\|_1 = \|\omega_-\|_1,
			\]
			gives the lower bound for $\gamma_+ + \gamma_-$, while the individual upper bounds for $\gamma_{\pm}$ follow from Lemma~\ref{lem:en_pos} and \eqref{eq:tilt_vert_bd_2}.
			
			The mass exchange argument in Lemma~\ref{lem:mass_sat}, with $|\ell_\kappa|$ in place of $|x_2|$, implies \eqref{eq:tilt_mass_sat}.
			The proof of Lemma~\ref{lem:supp} then applies with \eqref{eq:tilt_gamma_unif}, \eqref{eq:tilt_vert_bd_2}, and \eqref{eq:tilt_weak_bulk}.
			This implies \eqref{eq:tilt_supp}.
			On this support, $|\ell_\kappa| \le R + \pi \kappa$, so \eqref{eq:tilt_W_unif} and the proof of Lemma~\ref{lem:core} yield \eqref{eq:tilt_core}.
			Finally, we obtain \eqref{eq:tilt_sep} by Lemma~\ref{lem:sep}.
		\end{proof}

		\begin{remark}\label{rmk:tilt_traveling}
			Suppose that a maximizer $\omega \in \calS_{\nu, \varepsilon}^\kappa$ satisfies $\overline{\Omega_+ \cup \Omega_-} \cap \{\tld{x}_1 = \pi \} = \emptyset$ and
			\[
			\Omega_+ \subset \{ \ell_\kappa > 0 \}, \quad \Omega_- \subset \{ \ell_\kappa < 0 \}.
			\]
			Then \eqref{eq:tilt_EL} implies $\psi[\omega] - W\ell_\kappa = \gamma_+$ on $\partial\Omega_+$ and $\psi[\omega] - W\ell_\kappa = -\gamma_-$ on $\partial\Omega_-$.
			Since $\nabla^\perp \ell_\kappa = (1,\kappa)$ in a neighborhood of both cores, the vector field
			\[
			u - W(1,\kappa) = \nabla^\perp \left( \psi[\omega] - W\ell_\kappa \right)
			\]
			is tangent to both patch boundaries.
			Thus, $\omega$ generates a traveling vortex patch with velocity $W(1,\kappa)$.
		\end{remark}

		\subsection{Point-vortex limit}\label{sec:tilt_pv}
		
		\begin{definition}[Admissible point vortices]
			Let
			\[
			\Gamma := \{ x \in \calC : \tld{x}_1 = \pi \}
			\]
			denote the \emph{cut} used to define $\tld{x}_1$, and define
			\[
			\calP_{\nu,\kappa}
			:=
			\left\{
			(p_+,p_-) \in (\calC\setminus\Gamma)^2:
			p_+ \neq p_-,
			\quad
			\frac{|\ell_\kappa(p_+)| + |\ell_\kappa(p_-)|}{2} = \nu
			\right\}.
			\]
		\end{definition}
		
		For a pair $(q_+, q_-) \in \calP_{\nu, \kappa}$, let $a$ be its horizontal distance on $\bbT$ and $h$ be its vertical separation.
		The constraint and energy are symmetric in the two points, we may assume $q_{+,2} \ge q_{-,2}$.
		Let $d := \tld{q}_{+,1} - \tld{q}_{-,1}$.
		By the moment constraint,
		\[
		h - \kappa d = \ell_\kappa(q_+) - \ell_\kappa(q_-) \le |\ell_\kappa(q_+)| + |\ell_\kappa(q_-)| = 2\nu.
		\]
		For the relative position $(-a,h)$, $d$ can be $-a$ or $2\pi - a$.
        Hence
		\begin{equation}\label{eq:h_tilt}
			h \le
			\begin{cases}
				2\nu - \kappa a, & d = -a, \\
				2\nu + \kappa(2\pi - a), & d = 2\pi - a.
			\end{cases}
		\end{equation}
		Note that the latter bound holds for all pairs in $\calP_{\nu,\kappa}$, since $d \le 2\pi - a$.
		In particular, for $0 <a \le \pi$, it is attained by
		\[
		q_{\pm} := \pm \left( \pi - \frac{a}{2}, \,\, \nu + \kappa \left( \pi - \frac{a}{2} \right) \right), \quad \ell_{\kappa}(q_{\pm}) = \pm \nu.
		\]
		This pair is located across the cut and has representative difference $2\pi - a$.
		
		\begin{definition}[Reduced energy]
			Since the point-vortex energy
			\[
			-G(-a,h) = \frac{1}{4\pi} \log\left(2\cosh h - 2\cos a\right)
			\]
			increases in $h$ at fixed $a$, we define the reduced energy by
			\begin{equation}\label{eq:tilt_red}
				\Phi_{\nu,\kappa}(a) := \cosh \left( 2\nu + \kappa(2\pi-a) \right) - \cos a, \quad 0 \le a \le \pi.
			\end{equation}
			When $\kappa=0$, \eqref{eq:tilt_red} reduces to $\Phi_{\nu,0}(a) = \cosh(2\nu)-\cos a$, which is uniquely maximized at $a=\pi$, recovering the result in \S \ref{sec:karman}.
		\end{definition}

		\begin{remark}[The case $a = 0$: cut approaching]
			For any vertically aligned pair $(p_+,p_-) \in \calP_{\nu,\kappa}$, the horizontal representatives coincide.
			Hence its vertical separation satisfies
			\[
			h = |\ell_\kappa(p_+) - \ell_\kappa(p_-)| \le 
			|\ell_\kappa(p_+)| + |\ell_\kappa(p_-)| = 2\nu.
			\]
			For $\kappa > 0$, the admissible pairs $q_\pm$ defined above approach the cut $\Gamma$ from opposite sides as $a \searrow 0$, with vertical separation tending to $2\nu + 2\pi\kappa > 2\nu$.
			In this sense, the endpoint value $\Phi_{\nu,\kappa}(0)$ in \eqref{eq:tilt_red} is interpreted as the limit along \emph{admissible} pairs with $a > 0$, although the limiting pair itself is not admissible.
		\end{remark}

		\begin{figure}[h]
			\centering
			\includegraphics[scale=0.8]{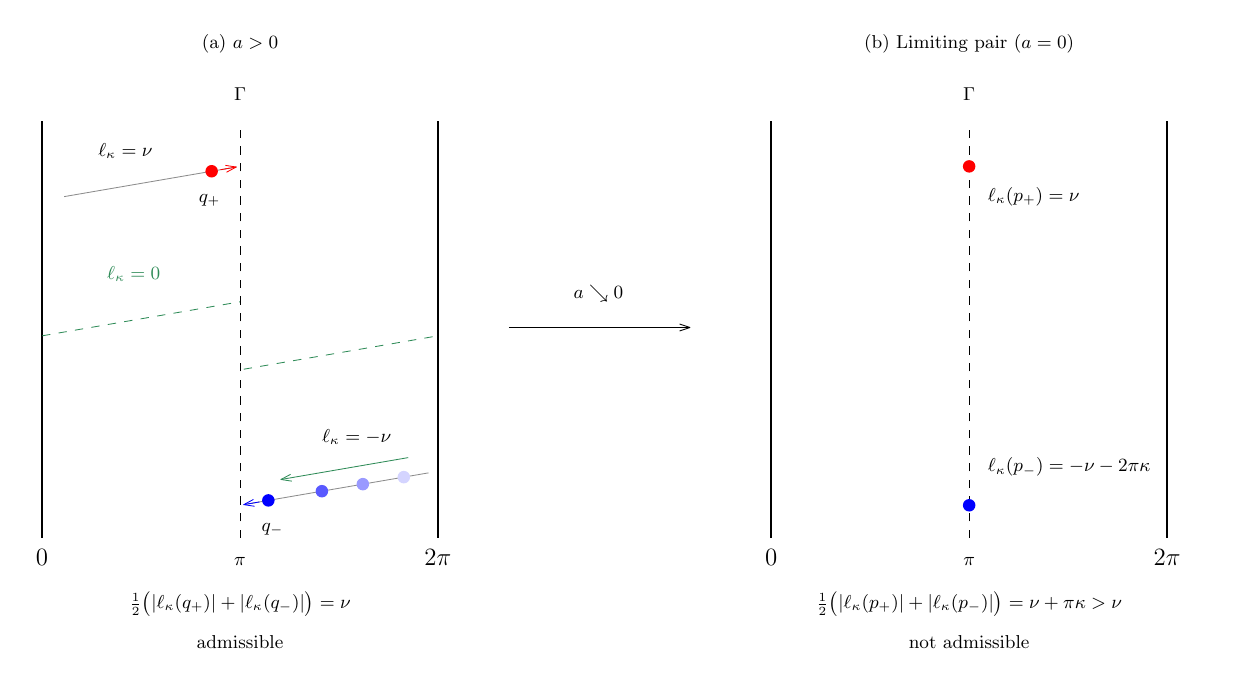}  
			\caption{Cut approaching admissible pairs.
				Two vortices translate along $\ell_\kappa = \pm \nu$, parallel to $\ell_\kappa = 0$, while preserving the tilted moment away from $\Gamma$.
				As both points approach the cut, their interaction energy converges, but the limiting pair is not admissible.} \label{fig:pv-cut-approach}
		\end{figure}

		\begin{proposition}\label{prop:tilt_pv_lim}
			Let $\varepsilon_j \to 0$ and let $\omega_j = \varepsilon_j^{-2} \left( \mathbf 1_{\Omega_{j,+}} - \mathbf 1_{\Omega_{j,-}} \right) \in \calS_{\nu,\varepsilon_j}^{\kappa}$.
			Then, up to a subsequence, there exist distinct points
			$p_\pm \in \calC$ such that
			\[
			\omega_{j,+} \rightharpoonup \frac{1}{2}\delta_{p_+}, \quad \omega_{j,-} \rightharpoonup \frac{1}{2}\delta_{p_-}.
			\]
			Set
			\[
			a := \dist_{\bbT} (p_{+, 1}, p_{-, 1}), \quad h := |p_{+, 2} - p_{-, 2}|.
			\]
			Then
			\begin{equation}
				\label{eq:0905_1}
				h = 2\nu + \kappa (2\pi - a), \quad \Phi_{\nu, \kappa}(a) = \max_{0 \le b \le \pi} \Phi_{\nu, \kappa}(b).
			\end{equation}
			Suppose in addition that $p_\pm \notin \Gamma$.
			Then
			\begin{equation}\label{eq:tilt_pv_moment}
				\frac{|\ell_\kappa(p_+)| + |\ell_\kappa(p_-)|}{2} = \nu,
			\end{equation}
			and $(p_+,p_-)$ maximizes $- G(q_+ - q_-)$ over $(q_+,q_-) \in \calP_{\nu,\kappa}$.
		\end{proposition}
		
		\begin{proof}
			The concentration follows from \eqref{eq:tilt_supp}, \eqref{eq:tilt_core}, and \eqref{eq:tilt_sep}.
			For all $x, y \in \calC$, we have
			\[
			|x_2 - y_2| \le |\ell_\kappa(x)| + |\ell_{\kappa}(y)| + \kappa \left( 2\pi - \dist_{\bbT}(x_1, y_1) \right).
			\]
			Integrating against $\omega_{j, +}(x)\omega_{j, -}(y)$, passing to the limit, and using the mass saturation \eqref{eq:tilt_mass_sat}, we see
			\[
			\frac{h}{4} \le \frac{\nu}{2} + \frac{\kappa}{4}(2\pi - a).
			\]
			In particular, $\cosh h - \cos a \le \Phi_{\nu, \kappa}(a)$.
			
			Fix $b \in (0, \pi]$ and set
			\[
			q_{\pm} := \pm \left( \pi - \frac{b}{2},  \nu + \kappa \left( \pi - \frac{b}{2} \right) \right).
			\]
			We separately translate two cores $\Omega_{j, +}$ and $\Omega_{j, -}$, say $\hat\Omega_{j, \pm}$ so that
			\[
			\frac{1}{|\hat\Omega_{j, \pm}|}\int_{\hat\Omega_{j, \pm}} (\tld{x}_1, x_2) dx = q_{\pm},
			\]
			respectively.
			We denote the resulting vorticity by $\hat\omega_j$.
			It is clear that $\hat\Omega_{j, +}$ and $\hat\Omega_{j, -}$ are disjoint, $ \overline{\hat\Omega_{j, +} \cup \hat\Omega_{j, -}} \cap \Gamma = \emptyset$, and lie in $\{\ell_\kappa > 0\}$ and $\{\ell_\kappa < 0\}$, respectively.
			Moreover, direct computation shows
			\[
			\| \ell_\kappa \hat\omega_j \|_1 = \frac{1}{2}\ell_\kappa(q_+) - \frac{1}{2}\ell_\kappa(q_-) = \nu.
			\]
			Thus $\hat\omega_j \in \calF_{\nu, \varepsilon_j}^\kappa$.
			Then, by the maximality we have
			\[
			0 \le E[\omega_j] - E[\hat\omega_j] \to \frac{1}{16\pi} \log \frac{\cosh h - \cos a}{\Phi_{\nu,\kappa}(b)}.
			\]
			Since $b \in (0, \pi]$ is arbitrary and $\Phi_{\nu, \kappa}$ is continuous, we obtain
			\[
			\max_{0 \le b \le \pi} \Phi_{\nu,\kappa}(b) =  \cosh h - \cos a = \Phi_{\nu, \kappa}(a).
			\]
			This proves \eqref{eq:0905_1}.

			If $p_\pm \notin \Gamma$, $|\ell_\kappa|$ is continuous in fixed neighborhoods of $p_\pm$, and the saturation of the tilted moment $\|\ell_{\kappa}\omega\|_1 = \nu$ yields \eqref{eq:tilt_pv_moment}.
			The maximizing property of $(p_+, p_-)$ follows directly from the above.
		\end{proof}
		
		For a limiting pair satisfying $p_{\pm} \notin \Gamma$, the limiting maximization also determines the separation of the two vortices in the \emph{tilted} coordinate.
		Set $s_\pm := \ell_\kappa(p_\pm)$, and then, by \eqref{eq:tilt_pv_moment}, we have $|s_+| + |s_-| = 2\nu$.
		Fix the horizontal coordinates $\tld{p}_{\pm, 1}$ of $p_\pm$ and then we may write
		\[
		p_{+,2} - p_{-,2} = \kappa \left( \tld{p}_{+,1} - \tld{p}_{-,1} \right) + s_+ - s_-.
		\]
		Under the constraint $|s_+| + |s_-| = 2\nu$, the difference
		$s_+ - s_-$ can take any value in $[-2\nu,2\nu]$.
		Since $-G(p_+-p_-)$ is strictly increasing with
		$|p_{+,2}-p_{-,2}|$ when the horizontal coordinates are fixed, the maximality of $(p_+,p_-)$ implies
		\begin{equation}\label{eq:tilt_opp_lim}
			|s_+ - s_-| = 2\nu.
		\end{equation}
		Indeed, otherwise replacing $s_+-s_-$ by one of the endpoints $\pm2\nu$ would strictly increase $|\kappa \left( \tld{p}_{+,1} - \tld{p}_{-,1} \right) + s_+ - s_-|$ while preserving the tilted-moment constraint.
		
		\begin{corollary}\label{cor:tilt_traveling}
			Under the assumptions of Proposition~\ref{prop:tilt_pv_lim}, we further assume that
			\begin{equation}
				\label{eq:tilt_interior}
				\max_{0 \le b \le \pi} \Phi_{\nu, \kappa}(b) > \Phi_{\nu, \kappa}(0).
			\end{equation}
			Then, for all sufficiently large $j$, the maximizer generates a traveling wave solution in the direction $(1, \kappa)$, up to reversing the direction.
			Moreover, after replacing $\omega_j$ by $-\omega_j$ if necessary, we may write
			\[
			p_+ - p_- = (-a,h), \quad 0 < a < \pi, \quad h > 0,
			\]
			where
			\begin{equation}\label{eq:tilt_dire}
				\kappa = \frac{\sin a}{\sinh h}.
			\end{equation}
			With this choice of sign, the traveling velocities satisfy
			\begin{equation}\label{eq:tilt_speed}
				W_j(1, \kappa) \to \frac{1}{8\pi} \frac{(\sinh h,\sin a)} {\cosh h-\cos a}.
			\end{equation}
		\end{corollary}
		
		\begin{proof}
			By \eqref{eq:0905_1}, \eqref{eq:tilt_interior}, and \eqref{eq:tilt_gap}, we have $a > 0$ and $h > 2\pi\kappa$.
			After replacing $\omega_j$ by $-\omega_j$ if necessary, we may assume $p_{+,2} > p_{-,2}$.
			For all sufficiently large $j$, \eqref{eq:tilt_core} allows a common horizontal translation $t_j$ placing both cores away from $\Gamma$, with the positive core to the right of the negative core in $(-\pi,\pi)$.
			If we apply an appropriate vertical translation further, both cores contained in $\{x_2 > \pi\kappa\}$ and $\{x_2 < -\pi\kappa\}$, respectively.
			We denote this translated vorticity by $\hat\omega_j$.
			
			For $x \in \Omega_{j,+}$ and $y \in \Omega_{j,-}$, let $d \in (0,2\pi)$ be horizontal representative difference of their common translation $x + t_j$ and $y + t_j$.
            Then, since the original difference is $d$ or $d-2\pi$, we have
			\[
			\widetilde{x_1 + t_j} - \tld{y_1 + t_j} \ge \tld{x}_1 - \tld{y}_1.
			\]
			By multiplying $\omega_{j,+}(x)\omega_{j,-}(y)$ and integrating, together with \eqref{eq:tilt_mass_sat}, we find
			\[
			\int_{\calC} \left( \tld{x_1 + t_j} \right) \omega_j\,dx \ge \int_{\calC} \tld{x}_1 \omega_j\,dx,
			\]
			and then
			\begin{align*}
				\| \ell_\kappa \hat\omega_j \|_1 = \int_{\calC} x_2 \omega_j\,dx - \kappa \int_{\calC} \left( \tld{x_1 + t_j} \right) \omega_j\,dx  \le \int_{\calC} \ell_\kappa \omega_j\,dx \le \nu.
			\end{align*}
			Since $E[\hat\omega_j] = E[\omega_j]$, $\hat\omega_j$ is also a maximizer.
			Moreover, $\| \ell_\kappa \hat\omega_j \|_1 = \nu$.
			If not, by slightly increasing the vertical separation of two cores of $\hat\omega_j$, $E[\hat\omega_j] > E[\omega_j]$ which is a contradiction.
			Hence, $\hat\omega_j \in \calS_{\nu,\varepsilon_j}^{\kappa}$.
			By Remark~\ref{rmk:tilt_traveling}, $\hat\omega_j$, and hence $\omega_j$ generates a traveling wave with velocity $W_j(1, \kappa)$.

			Applying the argument of Lemma~\ref{lem:local_poho} to $\hat\omega_j$, we obtain
			\[
			\frac{W_j}{2}(-\kappa,1) = - \iint_{\calC\times\calC} \omega_{j,+}(x)\omega_{j,-}(y) \nabla G(x-y)\,dxdy.
			\]
			Passing to the point-vortex limit in Proposition~\ref{prop:tilt_pv_lim}, we obtain $W_j(-\kappa,1) \to -\frac{1}{2}\nabla G(-a,h)$.
			This with the identity
			\[
			\nabla G(-a,h) = \frac{1}{4\pi} \frac{(\sin a,-\sinh h)}{\cosh h - \cos a} 
			\]
			proves \eqref{eq:tilt_dire} and \eqref{eq:tilt_speed}.
		\end{proof}
		
		\begin{remark}
			In Corollary~\ref{cor:tilt_traveling}, the parameter $\nu$ enters through the constraint on the limiting pair $(a,h)$.
			Once a limiting pair is fixed, its velocity is given by \eqref{eq:tilt_speed}.
		\end{remark}
		
		\subsection{Proof of Theorem~\ref{thm:tilt_small}}\label{sec:tilt_selection}	
		We use Proposition~\ref{prop:tilt_pv_lim} to derive sufficient conditions for selecting a prescribed relative position.
		We then verify these conditions for fixed $\nu > 0$ and sufficiently small $\kappa > 0$.
		
		\medskip
		\noindent\textbf{Choice for critical points.}
		Fix $(-a_0,h_0)$ with $0 < a_0 < \pi$ and $h_0 > 0$.
		Set
		\begin{align}
			\kappa &= \frac{\sin a_0}{\sinh h_0}, \label{eq:tilt_kappa} \\
			\nu &= \frac{1}{2}\left(h_0 - \kappa(2\pi-a_0)\right). \label{eq:tilt_inv}
		\end{align}
		Assume $\nu > 0$.
		These choices imply
		\begin{equation}\label{eq:critical_para}
			\Phi_{\nu,\kappa}'(a_0) = \sin a_0 - \kappa\sinh h_0 = 0,
		\end{equation}
		where we recall \eqref{eq:tilt_red} for the definition of $\Phi_{\nu, \kappa}$.
		
		\medskip
		\noindent\textbf{Local maximality.}
		At this critical point,
		\[
		\Phi_{\nu,\kappa}''(a_0) = \frac{(\cosh h_0-\cos a_0)(1+\cosh h_0\cos a_0)}{\sinh^2 h_0}.
		\]
		Hence $a_0$ is a strict local maximum if
		\begin{equation}\label{eq:tilt_loc}
			1 + \cosh h_0\cos a_0 < 0.
		\end{equation}
		
		\medskip
		\noindent\textbf{Global comparison.}
		Since
		\[
		\Phi_{\nu,\kappa}'''(a) = -\sin a - \kappa^3\sinh\left(2\nu+\kappa(2\pi-a)\right) < 0, \quad 0<a<\pi,
		\]
		and $\Phi_{\nu,\kappa}'$ is concave, there exists at most one interior local maximum of $\Phi_{\nu,\kappa}$.
		Since $\Phi_{\nu,\kappa}'(\pi)<0$, the only remaining competitor is $a=0$.
		This cut-approaching endpoint has value
		\[
		\Phi_{\nu,\kappa}(0) = \cosh(h_0+\kappa a_0)-1.
		\]
		Under \eqref{eq:tilt_loc}, $a_0$ is therefore the unique global maximizer if
		\begin{equation}\label{eq:tilt_uniq}
			\cosh h_0 - \cos a_0 > \cosh(h_0+\kappa a_0)-1.
		\end{equation}
		Note that this implies \eqref{eq:tilt_interior}.
		
		Therefore, under \eqref{eq:tilt_gap}, \eqref{eq:tilt_loc}, and \eqref{eq:tilt_uniq}, the limiting relative position is uniquely given by $(-a_0,h_0)$ in the standard orientation.
		Proposition~\ref{prop:tilt_pv_lim} and Corollary~\ref{cor:tilt_traveling} then yield traveling vortex patches for all sufficiently small $\varepsilon > 0$.

		\begin{proof}[Proof of Theorem~\ref{thm:tilt_small}]
			Fix $\nu > 0$.
			For sufficiently small $\kappa > 0$,
			\[
			\Phi_{\nu,\kappa}'\left( \frac{\pi}{2} \right) =   1 - \kappa\sinh\left( 2\nu + \frac{3\pi}{2}\kappa\right)>  0  >  -\kappa  \sinh (2\nu + \pi\kappa) = \Phi_{\nu,\kappa}'(\pi).
			\]
			Since $\Phi_{\nu,\kappa}' = \sin a -\kappa \sinh(2\nu + \kappa(2\pi - a))$ is strictly concave, there is a unique $a_\kappa \in (\pi/2,\pi)$ such that $\Phi_{\nu,\kappa}'(a_\kappa) = 0$.
			Set $h_\kappa := 2\nu + \kappa(2\pi - a_\kappa)$.
			The critical point condition \eqref{eq:critical_para} reads
			\[
			0 = \Phi_{\nu,\kappa}'(a_\kappa) = \sin a_\kappa - \kappa\sinh h_\kappa.
			\]
			Since $a_\kappa > 0$ and then $h_\kappa < 2\nu + 2\pi\kappa$, we obtain
			\[
			0 < \sin a_\kappa = \kappa\sinh h_\kappa
			\le \kappa\sinh(2\nu+2\pi\kappa) \to 0
			\quad \textrm{as } \kappa \to 0.
			\]
			In particular, $a_\kappa \to \pi$ and $h_\kappa \to 2\nu$ as $\kappa \to 0$.
			Consequently, for all sufficiently small $\kappa$, we have $\nu - 3\pi\kappa > 0$ and
			\[
			\begin{aligned}
				1 + \cosh h_\kappa \cos a_\kappa  &< 0, \\
				\cosh h_\kappa - \cos a_\kappa &>  \cosh(h_\kappa + \kappa a_\kappa) - 1 ,
			\end{aligned}
			\]
			that is, \eqref{eq:tilt_gap}, \eqref{eq:tilt_loc}, and \eqref{eq:tilt_uniq} hold.
			The preceding selection criterion, the mass and core estimates
			\eqref{eq:tilt_mass_sat} and \eqref{eq:tilt_core}, and
			Corollary~\ref{cor:tilt_traveling} yield the desired results
		\end{proof}
		
		\begin{remark}
			For example, to obtain
			\[
			(0,\mu)-(a_0,-\mu)=(-a_0,2\mu),
			\]
			we may choose
			\[
			\kappa = \frac{\sin a_0}{\sinh(2\mu)}, \quad \nu = \mu - \frac{(2\pi-a_0)\sin a_0}{2\sinh(2\mu)}.
			\]
			A common horizontal translation gives representatives $\pm (\pi - a_0/2, \mu)$, for which the tilted coordinates are $\pm\nu$.
			Under the selection conditions above, the limiting relative position is uniquely given by $(-a_0, 2\mu)$.
		\end{remark}  
		
		\bibliographystyle{acm}

	\end{document}